\documentclass[a4paper,dvipsnames,fleqn]{article}
\pdfoutput=1

\usepackage{comment}
\usepackage{amsmath}
\usepackage{amssymb}
\usepackage{amsthm}
\usepackage{bm}
\usepackage[backend=biber,bibencoding=utf8,bibstyle=alphabetic,citestyle=alphabetic,sorting=nyt,maxnames=4]{biblatex}
	\AtEveryBibitem{\clearfield{doi}}
	\AtEveryBibitem{\clearfield{isbn}}
	\AtEveryBibitem{\clearfield{issn}}
	\AtEveryBibitem{\clearfield{url}}
	
	\renewbibmacro{in:}{\ifentrytype{article}{}{\printtext{\bibstring{in}\intitlepunct}}}
	\AtEveryBibitem{\ifentrytype{book}{\clearfield{pages}}{}}
\usepackage[colorlinks=true,allcolors=purple]{hyperref}
\usepackage[strict]{changepage}
\usepackage{cleveref}
	\crefformat{equation}{#2(#1)#3}
\usepackage{dsfont}
\usepackage{enumitem}
\usepackage[OT2,OT1]{fontenc}
\usepackage{index}
\usepackage[cal=cm,scr=boondoxo]{mathalfa}
\usepackage{pifont}
\usepackage{rotating}
\usepackage{setspace}
\usepackage{stmaryrd}
\usepackage{textcomp}
\usepackage{tikz}
	\usetikzlibrary{arrows}
	\usetikzlibrary{patterns}
\usepackage{tikz-cd}
\usepackage[textwidth=3cm,colorinlistoftodos]{todonotes}
\usepackage{xfrac}

\numberwithin{equation}{section}			% Enumeration of equations
\swapnumbers						% Enumeration of theorem-like environments

\newcommand\cyr{%					% cyrillic font
\renewcommand\rmdefault{wncyr}%
\renewcommand\sfdefault{wncyss}%
\renewcommand\encodingdefault{OT2}%
\normalfont
\selectfont}
\DeclareTextFontCommand{\textcyr}{\cyr}

\newindex{sub}{idx-sub}{ind-sub}{Subject index}
\newindex{not}{idx-not}{ind-not}{List of notation}
\newcommand{\IndexN}[1]{\ensuremath{#1}\index[not]{${#1}$}}
\newcommand{\IndexS}[2]{{\emph{#1}}\index[sub]{#2}}
\newcommand{\IndexxS}[1]{\index[sub]{#1}}

\newcounter{Enum}				% Enumerated list
\newenvironment{Enumerate}{\begin{enumerate}[label={\rm({\roman*})}]}{\end{enumerate}}

\newcommand{\descriptionlabelsave}{}		% Itemized list
\newenvironment{Itemize}{%
	\renewcommand{\descriptionlabelsave}{\descriptionlabel}\renewcommand{\descriptionlabel}{$\triangleright$}%
	\begin{description}[leftmargin=15pt,itemindent=-5.2pt]}{%
	\end{description}\renewcommand{\descriptionlabel}{\descriptionlabelsave}}

\newcounter{StepsCount}				% Enumerated list with no indentation (e.g. steps in proof)

\newcounter{StepsRefCount}

\newenvironment{Ilist}{%			% Itemised list without indentation (not-enumerated steps)
	\begin{list}{$\triangleright$}{\leftmargin=0pt \labelwidth=11pt \itemindent=\labelwidth%
	\itemsep=5pt\listparindent=\parindent}}{\end{list}}

\theoremstyle{plain}
	\newtheorem{lemma}{Lemma}[section]
	\newtheorem{proposition}[lemma]{Proposition}
	\newtheorem{theorem}[lemma]{Theorem}
	\newtheorem{corollary}[lemma]{Corollary}
	\newcommand{\GenericTheoremName}{}\newtheorem{generictheorem}[lemma]{\GenericTheoremName}
\theoremstyle{definition}
	\newtheorem{definition}[lemma]{Definition}
	\newcommand{\GenericDefinitionName}{}\newtheorem{genericdefinition}[lemma]{\GenericDefinitionName}
\theoremstyle{remark}
	\newtheorem{remark}[lemma]{Remark}
	\newtheorem{example}[lemma]{Example}
	\newcommand{\GenericRemarkName}{}\newtheorem{genericremark}[lemma]{\GenericRemarkName}
\newenvironment{Lemma}{\begin{lemma}}{\par\noindent\rule{5em}{1pt}\end{lemma}}
\newenvironment{Proposition}{\begin{proposition}}{\par\noindent\rule{5em}{1pt}\end{proposition}}
\newenvironment{Theorem}{\begin{theorem}}{\par\noindent\rule{5em}{1pt}\end{theorem}}
\newenvironment{Corollary}{\begin{corollary}}{\par\noindent\rule{5em}{1pt}\end{corollary}}

\newenvironment{Definition}{\begin{definition}}{\par\noindent\rule{5em}{1pt}\end{definition}}

\newenvironment{Remark}{\begin{remark}}{\par\noindent\rule{5em}{0.5pt}\end{remark}}
\newenvironment{Example}{\begin{example}}{\par\noindent\rule{5em}{0.5pt}\end{example}}

\newcommand{\mc}[1]{{\mathcal{#1}}}			% --- abbreviation ---
\newcommand{\ms}[1]{{\mathscr{#1}}}			% --- abbreviation ---
\newcommand{\mf}[1]{{\mathfrak{#1}}}			% --- abbreviation ---
\newcommand{\bb}[1]{{\mathbb{#1}}}			% --- abbreviation ---
\newcommand{\ov}{\overline}				% --- abbreviation ---
\newcommand{\mr}{\mathring}				% --- abbreviation ---
\newcommand{\wt}{\widetilde}				% --- abbreviation ---
\DeclareMathOperator{\IM}{Im}				% imaginary part
\renewcommand{\Im}{\IM}

\newcommand{\Side}[1]{\hfill{#1}\kern10pt}		% text put on the left side of line with offset
\newcommand{\FD}[5]{%					% definition of function from {#1} to {#2} mapping {#4} to {#5}
	\DF\left\{\begin{array}{rcl}{#1}&\to &{#2}\\[#3pt] {#4}&\mapsto &{#5}\end{array}\right.}

\newcommand{\smmatrix}[4]{\Bigl(			% small matrix for use in textline
\begin{smallmatrix}
\hspace*{-0.2ex} #1 \hspace*{0.2ex} & \hspace*{0.2ex} #2 \hspace*{-0.2ex}
\\[0.5ex]
\hspace*{-0.2ex} #3 \hspace*{0.2ex} & \hspace*{0.2ex} #4 \hspace*{-0.2ex}
\end{smallmatrix}
\Bigr)}
\newcommand{\Dummy}{\text{\textvisiblespace\kern1pt}}	% Platzhaltersymbol fuer Funktionsargumente
\newcommand{\Smallo}{{\rm o}}				% small o
\newcommand{\BigO}{{\rm O}}				% big o
\newcommand{\circleq}{{\mkern5mu\mr{\leq}\mkern5mu}}	% inequalities with circle
\newcommand{\circgeq}{{\mkern5mu\mr{\geq}\mkern5mu}}	% inequalities with circle

\DeclareMathOperator{\Id}{id}				% identity map
\DeclareMathOperator{\Supp}{supp}			% support
\DeclareMathOperator{\Clos}{Clos}			% closure
\DeclareMathOperator{\Ran}{ran}				% range
\DeclareMathOperator{\Tr}{tr}				% trace
\DeclareMathOperator{\SL}{SL}				% special linear group
\DeclareMathOperator{\Loc}{loc}				% "locally"
\DeclareMathOperator{\Ess}{ess}				% "essential"

\newcommand{\DS}{\mid\mkern3mu}				% delimiter for set definition
\newcommand{\DQ}{\mkern6mu}				% distance for successive quantors
\newcommand{\DP}{{:\kern5pt}}				% delimiter for predicate formula
\newcommand{\DF}{\colon}				% delimiter for function domain/codomain
\newcommand{\DE}{\mathrel{\mathop:}=}			% defining equality
\newcommand{\DI}{\mathrel{\mathop:}\Leftrightarrow}	% defining equivalence
\newcommand{\DD}{\mkern4mu\mathrm{d}}			% distance and rm-d for integration differential
\newcommand{\CAS}{&\text{if}\ }				% delimiter in "cases" environment
\newcommand{\CASO}{&\text{otherwise}}			% delimiter in "cases" environment

\DeclareMathOperator{\Ind}{Ind}				% index of regular variation

\DeclareMathOperator{\Arccot}{Arccot}			% Arccot

\newcommand{\PART}[2]{%					% Formatting of heading for "Part"
	\begin{center}
		\rule{100mm}{1pt}
		\\[3mm]
		{\Large\bf PART {#1}}
		\\[2.5mm]
		{\Large\bf {#2}}
		\\[1mm]
		\rule{100mm}{1pt}
		\\[2mm]
		\phantom{}
	\end{center}
	\addtocontents{toc}{\protect\contentsline{section}{\protect\numberline{}\textcolor{Sepia}{\kern-15pt 
	$\triangleright\triangleright$ Part {#1}\,: {#2}}}{}{}}
}

\newcommand{\AUXILIARY}[3]{%					% Formatting of heading for "Auxiliary"
	\begin{center}
		\rule{100mm}{1pt}
		\\[3mm]
		{\Large\bf {#1}}
		\\[2.5mm]
		{\Large\bf {#2}}
		\\[1mm]
		\rule{100mm}{1pt}
		\\[2mm]
		\phantom{}
	\end{center}
	\addtocontents{toc}{\protect\contentsline{section}{\protect\numberline{}\textcolor{Sepia}{\kern-15pt 
	{#3}\,: {#2}}}{}{}}
}

\newcommand{\Intro}[1]{%				% Formatting for introductory text to "Part"
	{\sf\begin{spacing}{1.1}\noindent{#1}\end{spacing}}
	\hspace*{0pt}\\[0mm]
}

\newenvironment{Rlist}{%				% List of remarks and comments 
	\begin{list}{\raisebox{-1pt}{{\large\ding{100}}}}{\leftmargin=0pt \labelwidth=13pt \itemindent=\labelwidth%
	\itemsep=5pt\listparindent=\parindent}}{\end{list}}

\newcommand{\REMARKS}[1]{%
	\medskip\pagebreak[3]
	\begin{center}
		$\bm{\triangleright\triangleright}\quad$\raisebox{-1pt}{{\large\bf Remarks}}$\quad\bm{\triangleleft\triangleleft}$
	\end{center}
	\begin{Rlist}
		{#1}
	\end{Rlist}
}

\begin{document}

\begin{flushleft}
	{\Large\bf Spectral properties of canonical systems:\\[2mm] discreteness and distribution of eigenvalues}
	\\[5mm]
	\textsc{
	Jakob Reiffenstein
	\,\ $\ast$\,\ 
	Harald Woracek
		\hspace*{-14pt}
		\renewcommand{\thefootnote}{\fnsymbol{footnote}}
		\setcounter{footnote}{2}
		\footnote{This work was supported by the project I~4600 of the Austrian Science Fund (FWF),
			and by the Sverker Lerheden foundation.}
		\renewcommand{\thefootnote}{\arabic{footnote}}
		\setcounter{footnote}{0}
	}
	\\[6mm]
	{\small
	\textbf{Abstract:}
		In this survey paper we review classical results and recent progress about a certain topic in the spectral 
		theory of two-dimensional canonical systems. Namely, we consider the questions whether the spectrum $\sigma$ is 
		discrete, and if it is, what is its density. 
		Here we measure density by the growth of the counting function of $\sigma$ in the sense of integrability or 
		$\limsup$-conditions relative to suitable comparison functions. 
		These questions have been around for many years. However, full answers have been obtained only very 
		recently -- in partly still unpublished work. 
		
		\smallskip
		The way we measure density of eigenvalues must be clearly distinguished from spectral asymptotics, where one 
		asks for an asymptotic expansion of eigenvalues. We explicitly and on purpose do not go into this direction and 
		do not present any results about spectral asymptotics.

		\smallskip
		We understand this survey as a focussed 
		presentation of results revolving around the initially stated questions, and not as an exhaustive account on the
		literature which is in the one or other way related to the area. 
		It does not contain any proofs, but we do comment on proof methods.

		\smallskip
		The theorems we present are very diverse. They rely on different methods from operator theory, complex
		analysis, and classical analysis, and beautifully invoke the interplay between these areas. 
		Besides the fundamental theorems solving the problem and several selected additions, we decided to include 
		a number of results devoted to the spectral theory of Jacobi matrices, in other words, 
		to the Hamburger moment problem. The connection between the theories of canonical systems and moment problems
		gives rise to important new insights, yet seems to be less widely known than other connections, e.g., with
		Krein-Feller or Schr\"odinger operators.

		\smallskip
		We compile all necessary prerequisites from the general theory to make the exposition as self-contained as 
		possible. It is our aim that this survey can be read and enjoyed by a broad community, including 
		non-specialist readers.
	\\[3mm]
	{\bf AMS MSC 2020:}
	34\,L\,15,\ 37\,J\,99,\ 30\,D\,15,\ 47\,B\,10,\ 47\,B\,36,\ 44\,A\,60
	}
\end{flushleft}

%---------
%   TEXTBODY
%---------

%**************************************************************************
%***                            Last Change: Mon 31 Mar 2025 12:27
%***   < INTRO >
%***
%**************************************************************************

%
%
%
\section*{Introduction}
\addcontentsline{toc}{section}{Introduction}

The present paper is a survey, not containing any proofs, where we review classical results and recent progress about a certain
topic in the spectral theory of two-dimensional canonical systems. These are first-order systems of the form 
\[
	y'(t)=zJH(t)y(t),\qquad t\in(a,b),
\]
where $-\infty<a<b\leq\infty$, $z\in\bb C$ is the spectral parameter, $J$ is the symplectic matrix $J\DE\smmatrix 0{-1}10$, 
and $H$ is the Hamiltonian of the system. We deal with systems whose Hamiltonian satisfies
$H(t)\in\bb R^{2\times 2}$, $H(t)^T=H(t)$, $H(t)\geq 0$ for a.a.\ $t\in(a,b)$, 
$H\in L^1((a,c),\bb R^{2\times 2})$ for all $c\in(a,b)$, and $H(t)\neq 0$ a.e.

Canonical systems play an important role in many contexts; we mention two:
\begin{Itemize}
\item Several particular equations can be transformed to the form of a canonical systems by sometimes more and sometimes less
	explicit formulae. This applies in particular to Schr\"odinger- or Sturm-Liouville equations, Dirac systems, 
	Jacobi matrices, and Krein-Feller operators (Krein strings). 
\item Canonical systems are a universal model for self-adjoint operators with simple spectrum. Every such operator is unitarily
	equivalent to the differential operator induced by some canonical system in a natural way. 
\end{Itemize}
Each two self-adjoint realisations of the above differential equation are at most $2$-dimensional perturbations of each
other, and hence share most of their spectral properties. When talking about such properties we may thus pick up any self-adjoint
realisation. 

The theorems we present revolve around the following two questions:
\begin{itemize}
\item[(1)] When is the spectrum discrete ?
\item[(2)] If the spectrum is discrete, how densely is it distributed ?
\end{itemize}
Here we understand ``density'' of a discrete subset $\sigma$ of $\bb R$ as being measured by the growth of the counting function 
\[
	n(r)\DE\#\Big(\sigma\cap(-r,r)\Big)
\]
when $r$ tends to infinity. 
We measure the speed of growth of $n(r)$ in ways familiar from complex analysis like finiteness of 
$\limsup_{r\to\infty}\frac{n(r)}{\ms g(r)}$ or $\int_1^\infty\frac{n(r)}{\ms g(r)}\DD r$, where $\ms g(r)$ is a suitable
comparison function (for a start think of a power $\ms g(r)=r^\rho$ with $\rho>0$). It must be said very clearly that we 
do not discuss results about actual asymptotics of $n(r)$. There is a vast literature about spectral asymptotics of various
types of equations, but it is most explicitly not our aim to go into this direction. 

Question~(1) is answered by a single theorem reading as ``The spectrum is discrete, if and only if...''. Contrasting this, 
question~(2) has many facets. The full answer boils down to two theorems: one settling the case of ``dense'' spectrum, 
and another settling the case of ``sparse'' spectrum. It is a good, though not fully correct, intuition to imagine 
the set $\bb Z$ of integers as the border between those two situations. These two cases are very different in their nature, and so
are the proof methods.

\subsection*{A guide to reading the paper}
\addcontentsline{toc}{subsection}{A guide to reading the paper}

On the top level of structuring, the content is arranged in five parts.
\begin{flushleft}
	\textit{Part~I: Fundamental Theorems}\ \dotfill\quad\pageref{U126}
\end{flushleft}
We present the fundamental theorems that answer questions~(1) and (2). 
In the discreteness characterisation and in the case of dense spectrum, the formulae are explicit in terms of the Hamiltonian 
and can be applied in practice. For sparse spectrum formulae are again explicit, but evaluating them is much more difficult
(we discuss situations where the density of the spectrum can be described more explicitly later on).
\begin{flushleft}
	\textit{Part~II: The limit circle case}\ \dotfill\quad\pageref{U146}
\end{flushleft}
We discuss the situation where $H$ is integrable on the whole interval $(a,b)$. In this case the spectrum is always sparse,
and there is a tight connection between the distribution of the spectrum and the growth of a certain entire function. 
This opens up another way of approach, which in particular leads to an algorithm to evaluate the general formula up to a
small error.
\begin{flushleft}
	\textit{Part~III: Moment problems and Jacobi matrices}\ \dotfill\quad\pageref{U127}
\end{flushleft}
The following three notions are just different ways to view one and the same object: moment sequences, Jacobi
matrices, Hamiltonians of a certain discrete form (which we call Hamburger Hamiltonians). 
The connection between moment sequences and Jacobi matrices is standard, the connection with
canonical systems seems to be less widely known, though it has been around for many years.
This interaction proves to be very fruitful, and we use it to discuss question~(2) in this setting. 
In this part we focus on results stated in terms of the moment sequence, its orthonormal polynomials, and its Jacobi parameters. 
Some of them are obtained making a detour via Hamburger Hamiltonians. 
\begin{flushleft}
	\textit{Part~IV: Hamburger Hamiltonians}\ \dotfill\quad\pageref{U147}
\end{flushleft}
In this part we focus on results about Hamburger Hamiltonians. When the Hamiltonian does not vary too wildly a fairly complete
picture can be given. Translating the results presented in this part to Jacobi parameters and the moment sequence is usually hard 
(if not impossible).
\begin{flushleft}
	\textit{Part~V: Additions and Examples}\ \dotfill\quad\pageref{U128}
\end{flushleft}
We present a selection of theorems and examples about various topics. 
Among them a construction of Hamiltonians whose monodromy matrix has prescribed growth, and a result where in one very concrete 
situation not only finiteness of $\limsup_{r\to\infty}\frac{n(r)}{\ms g(r)}$ (for appropriate $\ms g(r)$) is shown, but the 
actual value of this limit is determined. Further, we give results that illuminate the ways in which the growth of $n(r)$ 
depends on the Hamiltonian, and illustrate the pecularities of sparse spectrum.
Last but not least we give a brief account on the literature about Krein strings going back to the pioneering work of M.G.~Krein
in the 1950's. The results about the string equation were a great source of motivation and inspiration when developing the
general theory. 

\bigskip
Stepping down one level in our structuring hierarchy, each of Parts~I--V consists of several sections, in each of which we
focus on one specific result or one specific circle of ideas. 
A brief description guiding through the sections of each single part is provided on the first page of the respective part.

Moving on to one single section, each of them is manufactured as
follows: we start with a few introductory words, then present the results to which the section is dedicated, and close with
remarks about original sources, proof methods, connections to other results, and similar. 

Preceeding Parts~I-V we provide a prologue, where we place some subjects specific for the theory of canonical
systems. The aim is to prepare a common ground for all readers. We close the paper with an appendix, where we collect a few
subjects that we need from general theory and that might not be familiar to every reader. We include a list of notation, a 
subject index, and the bibliography at the end of the paper.

Concerning dependencies between the different sections, we tried to implement the following idea:
\begin{Itemize}
\item All readers should browse through \Cref{U113} in the prologue. After having done that, each reader is welcome to proceed
	to any section in Parts~I--V according to his or her interest, but should be prepared to return to \Cref{U114,U115} of
	the prologue or to visit the appendix if necessary.
\item As a rough guideline, the sections in Parts~I--V can be read independently of each other. 
\end{Itemize}

\subsection*{Table of contents}
\addcontentsline{toc}{subsection}{Table of contents}

\vspace*{-10mm}
\tableofcontents

%**************************************************************************
%***                            Last Change: Mon 31 Mar 2025 11:28
%***   < PROLOGUE >
%***
%**************************************************************************

\clearpage
\AUXILIARY{PROLOGUE}{Two-dimensional canonical systems}{Prologue}
\label{U125}

%%%%%%%%%%%%%%%%%%%%%%%%%%%%%%%%%%%%%%%%%%%%%%%%%%%%%%%%%%%%%%%%%%%%%%%%%%%%%%%%%%%%%%%%%

\Intro{%
	In this part we collect some basic notions and results from the theory of canonical systems and their spectral theory.
	Furthermore, we discuss two particular types of Hamiltonians. 
	This is done in order to provide prerequisites necessary for reading and
	understanding the content of the core Parts~I-V. It is by no means an account of the whole theory; the compilation is
	made specifically for the purpose of this paper and merely covers the present needs. 

	In \Cref{U113} we present what is needed in all further sections. This includes Weyl's limit circle/limit point
	distinction, the monodromy matrix, and the Weyl coefficient. 
	In \Cref{U114} we briefly discuss the operator model associated with a canonical system, which is needed in Part~I. 
	The operator theoretic viewpoint is also a major source of motivation to deal with growth properties of the monodromy
	matrix, which are discussed specifically in Part~II. 
	\Cref{U115} is devoted to two particular classes of Hamiltonians. First, we introduce 
	Hamburger Hamiltonians and explain their relationship to Jacobi matrices and the power moment problem; 
	this class is further studied in Parts~III and IV. Second, we introduce Pontryagin-type Hamiltonians, a class 
	that forms a modest extension of limit circle case; these occur in Part~II. 
}
\begin{center}
	{\large\bf Table of contents}
\end{center}
\begin{flushleft}
	\S\,\ref{U113}.\ Canonical systems\ \dotfill\quad\pageref{U113}
	\\[1mm]
	\quad \S\S\,\ref{U137}.\ The differential equation\ \dotfill\quad\pageref{U137}
	\\[1mm]
	\quad \S\S\,\ref{U139}.\ Weyl's method of nested disks\ \dotfill\quad\pageref{U139}
	\\[1mm]
	\quad \S\S\,\ref{U138}.\ The monodromy matrix\ \dotfill\quad\pageref{U138}
	\\[1mm]
	\quad \S\S\,\ref{U116}.\ The Weyl coefficient\ \dotfill\quad\pageref{U116}
	\\[1mm]
	\S\,\ref{U114}.\ The operator model\ \dotfill\quad\pageref{U114}
	\\[1mm]
	\S\,\ref{U115}.\ Two particular classes of systems\ \dotfill\quad\pageref{U115}
	\\[1mm]
	\quad \S\S\,\ref{U117}.\ Hamburger Hamiltonians: the power moment problem\ \dotfill\quad\pageref{U117}
	\\[1mm]
	\quad \S\S\,\ref{U118}.\ Pontryagin-type Hamiltonians: generalisation of limit circle case\ \dotfill\quad\pageref{U118}
	\\[1mm]
\end{flushleft}
\makeatother
\renewcommand{\thesection}{\Alph{section}}
\renewcommand{\thelemma}{\Alph{section}.\arabic{lemma}}
\makeatletter
\setcounter{section}{0}
\clearpage

%%%%%%%%%%%%%%%%%%%%%%%%%%%%%%%%%%%%%%%%%%%%%%%%%%%%%%%%%%%%%%%%%%%%%%%%%%%%%%%%%%%%%%%%%

%
%
%
\section[{Canonical systems}]{Canonical systems}
\label{U113}

References for the material presented in this section are 
\cite{behrndt.hassi.snoo:2020,remling:2018,romanov:1408.6022v1,winkler:1995,dym.mckean:1976,gohberg.krein:1970,
debranges:1968,atkinson:1964} and, historically, \cite{weyl:1910,titchmarsh:1946}.

\subsection[{The differential equation}]{The differential equation}
\label{U137}

A two dimensional \IndexS{canonical system}{canonical system} is a differential equation of the form 
\begin{equation}
\label{U84}
	y'(t)=zJH(t)y(t)
\end{equation}
on some interval $(a,b)$ where $-\infty<a<b\leq\infty$. Here $z\in\bb C$ is the spectral parameter, 
\IndexN{J} is the symplectic matrix $J\DE\smmatrix 0{-1}10$, 
and $H\in L^1_{\Loc}((a,b),\bb C^{2\times 2})$ is the \IndexS{Hamiltonian}{Hamiltonian} of the system. 
A solution $y:(a,b)\to\bb C^2$ is supposed to be locally absolutely continuous and satisfy \cref{U84} for a.a.\ $t\in(a,b)$. 

We consider a class of canonical systems whose Hamiltonian has certain analytic and algebraic properties. 

\begin{Definition}
\label{U52}
	Let $-\infty<a<b\leq\infty$. We denote by $\IndexN{\bb H_{a,b}}$ the set of all measurable functions 
	$H\DF(a,b)\to\bb R^{2\times 2}$, such that 
	\begin{Itemize}
	\item for a.a.\ $t\in(a,b)$ it holds that $H(t)=H(t)^T$ and $H(t)\geq 0$;
	\item for each $c\in(a,b)$ we have $H\in L^1\big((a,c),\bb R^{2\times 2}\big)$;
	\item the set $\{t\in(a,b)\DS H(t)=0\}$ has measure zero.
	\end{Itemize}
\end{Definition}

\noindent
Given $H\in\bb H_{a,b}$, the initial value problem 
\begin{equation}
\label{U291}
	\left\{
	\begin{array}{l}
		\frac{\partial}{\partial t}W_H(t,z)J=zW_H(t,z)H(t),\qquad t\in(a,b)\text{ a.e.},
		\\[2mm]
		W_H(a,z)=I,
	\end{array}
	\right.
\end{equation}
has a unique absolutely continuous solution $W_H\DF[a,b)\times\bb C\to\bb C^{2\times 2}$. 
We refer to $\IndexN{W_H(t,z)}$ as the \IndexS{fundamental solution}{fundamental solution} 
of the canonical system \cref{U84}. The entries of $W_H(t,z)$ are denoted as $\IndexN{w_{H,ij}(t,z)}$, $i,j\in\{1,2\}$. 

\begin{Remark}
\label{U292}
	Compared to the equation \cref{U84} we passed in \cref{U291} to transposes, so that now the transposed rows of 
	$W_H(t,z)$ are solutions of \cref{U84}. This is practical for technical reasons but has no intrinsic meaning. 
\end{Remark}

\noindent
Intuitively speaking a transformation of the independent variable in \cref{U84} will not change any essential properties 
of the equation \cref{U84} and its solutions. This is formalised by the following notion. 

\begin{Definition}
\label{U85}
	Let $-\infty<a<b\leq\infty$ and $-\infty<\tilde a<\tilde b\leq\infty$, and let $H\in\bb H_{a,b}$ and
	$\tilde H\in\bb H_{\tilde a,\tilde b}$. We say that $\tilde H$ is a
	\IndexS{reparameterisation}{reparameterisation} of $H$, if there exists an increasing bijection 
	$\varphi\DF[\tilde a,\tilde b)\to[a,b)$ such that $\varphi$ and $\varphi^{-1}$ are locally absolutely continuous and 
	$\tilde H(s)=H(\varphi(s))\varphi'(s)$, $s\in(\tilde a,\tilde b)$ a.e.
\end{Definition}

\noindent
If $\tilde H(s)=H(\varphi(s))\varphi'(s)$, the fundamental solution transforms accordingly as
\[
	W_{\tilde H}(s,z)=W_H(\varphi(s),z),\qquad s\in(\tilde a,\tilde b)
	.
\]

\begin{Remark}
\label{U290}
	Let $H\in\bb H_{a,b}$. The function $\IndexN{\mf t(t)}\DE\int_a^t\Tr H(s)\DD s$ is increasing, absolutely continuous,
	and has absolutely continuous inverse. 
	Hence, we may use $\varphi\DE\mf t^{-1}$ to produce a reparameterisation $\tilde H$ of $H$. 
	This reparameterisation satisfies $\Tr\tilde H(s)=1$ a.e.
	Note here that a positive semidefinite matrix $A$ is nonzero if and only if $\Tr A>0$.
\end{Remark}

\noindent
Hamiltonians whose trace is (a.e.) identically equal to $1$ are called 
\IndexS{trace-normalised}{trace-normalised}\index[sub]{Hamiltonian!trace-normalised}.
Due to the above remark one can often restrict considerations to trace-normalised Hamiltonians. 
However, it is not always a good idea to restrict generality.

We observe that reparameterisation induces an equivalence relation on the set 
\[
	\bigcup_{-\infty<a<b\leq\infty}\bb H_{a,b}
	.
\]
In the following we use the notation
\[
	\IndexN{\xi_\phi}\DE\binom{\cos\phi}{\sin\phi},\qquad\phi\in\bb R
	.
\]
There may exist intervals where a Hamiltonian is of a particularly simple form. 

\begin{Definition}
\label{U67}
	Let $H\in\bb H_{a,b}$. A nonempty interval $(c,d)\subseteq(a,b)$ is called \IndexS{indivisible}{indivisible interval} 
	for $H$, if there exists $\phi\in\bb R$ such that 
	\[
		H(t)=\Tr H(t)\cdot\xi_\phi\xi_\phi^T\qquad\text{for }t\in(c,d)\text{ a.e.}
	\]
	The number $\phi$, which is determined up to integer multiples of $\pi$, is called the 
	\IndexS{type}{indivisible interval!type} of the indivisible interval $(c,d)$.

	We call $H\in\bb H_{a,b}$ \IndexS{definite}{definite}\index[sub]{Hamiltonian!definite} if $(a,b)$ is not indivisible.
\end{Definition}

\noindent
We note that $H$ is definite if and only if the matrix $\int_a^b H(t)\DD t$ is positive definite (and not only positive
semidefinite, which it always is). 

Sometimes indivisible intervals are considered exceptional, but this is not at all the case; on the contrary, see \Cref{U117}.

\subsection[{Weyl's method of nested disks}]{Weyl's method of nested disks}
\label{U139}

We denote by $\bb C^+$ the open upper half-plane 
\[
	\IndexN{\bb C^+}\DE\big\{z\in\bb C\DS \Im z>0\big\}
	.
\]
Let $H\in\bb H_{a,b}$ and $(t,z)\in[a,b)\times\bb C$, and consider the fractional linear transformation 
\begin{equation}
\label{U90}
	\zeta\mapsto\frac{w_{H,11}(t,z)\zeta+w_{H,12}(t,z)}{w_{H,21}(t,z)\zeta+w_{H,22}(t,z)}
\end{equation}
as a map of the Riemann sphere $\bb C\cup\{\infty\}$ onto itself. 
Moreover, denote by $\ov{\bb C^+}$ the closure of the upper half-plane in the sphere, explicitly
\[
	\IndexN{\ov{\bb C^+}}\DE\bb C^+\cup\bb R\cup\{\infty\}
	.
\]
For each $t\in[a,b)$ and $z\in\bb C^+$ the transformation \cref{U90} maps $\ov{\bb C^+}$ onto some closed disk 
$\IndexN{\Omega_{t,z}}$ contained in $\ov{\bb C^+}$. These disks $\IndexN{\Omega_{t,z}}$ are called 
\IndexS{Weyl disks}{Weyl disk}. 
For each fixed $z\in\bb C^+$ the Weyl disks $(\Omega_{t,z})_{t\in[a,b)}$, are nested in the sense that 
\[
	\forall s,t\in[a,b)\DP s\leq t\Rightarrow\Omega_{s,z}\supseteq\Omega_{t,z}
\]
and hence there are two possible scenarios:
\begin{Itemize}
\item The intersection $\bigcap_{t\in[a,b)}\Omega_{t,z}$ is a disk with positive radius;
\item The intersection $\bigcap_{t\in[a,b)}\Omega_{t,z}$ consists of a single point.
\end{Itemize}
The chordal radius of $\Omega_{t,z}$ splits as $\rho(z)\cdot(\int_a^t\Tr H(s)\DD s)^{-1}$ with some continuous 
function $\rho$ depending only on $z$. Hence, if the limit disk has positive radius for one $z$ then this is true for every $z$.
Moreover, it has positive radius if and only if $\int_a^b\Tr H(s)\DD s<\infty$,
and in turn if and only if $H\in L^1((a,b),\bb R^{2\times 2})$.

This motivates the following terminology.

\begin{Definition}
\label{U88}
	Let $H\in\bb H_{a,b}$. We say that $H$ is in 
	\IndexS{limit circle case}{limit circle case}\index[sub]{Hamiltonian!limit circle}, if 
	$H\in L^1((a,b),\bb R^{2\times 2})$. Otherwise, $H$ is in 
	\IndexS{limit point case}{limit point case}\index[sub]{Hamiltonian!limit point}.
\end{Definition}

\noindent
We note that two Hamiltonians that are reparameterisations of each other are together in limit circle or limit point case. 

\subsection[{The monodromy matrix}]{The monodromy matrix}
\label{U138}

Assume that $H$ is in limit circle case, i.e., $H\in\bb H_{a,b}\cap L^1((a,b),\bb R^{2\times 2})$. 
Then the fundamental solution can be extended continuously to the right endpoint $b$. 
We call $W_H(b,z)$ the \IndexS{monodromy matrix}{monodromy matrix} of $H$, and denote it as 
$\IndexN{W_H(z)}$ and its entries as $\IndexN{w_{H,ij}(z)}$, $i,j\in\{1,2\}$. 
The fractional linear transformations \cref{U90} depend continuously on $t$, and hence 
the limit disk $\bigcap_{t\in[a,b)}\Omega_{t,z}$ is the image of $\ov{\bb C^+}$ under
\[
	\zeta\mapsto\frac{w_{H,11}(z)\zeta+w_{H,12}(z)}{w_{H,21}(z)\zeta+w_{H,22}(z)}
	.
\]
The monodromy matrix has the power series expansion
\[
	W_H(z)=\sum_{n=0}^\infty W_{H,n}z^n
\]
where the coefficients $W_{H,n}$ are defined by the recurrence 
\[
	W_{H,0}=I,\qquad \forall n\in\bb N\DP W_{H,n+1}=\int_a^b W_{H,n}(s)H(s)\DD s
	.
\]
We note that 
\[
	\forall n\in\bb N\DP W_{H,n}\in\bb R^{2\times 2}\ \wedge\ 
	\|W_{H,n}\|\leq\frac 1{n!}\Big(\int_a^b\|H(t)\|\DD t\Big)^n
	.
\]
The matrix $W_H(z)$ is analytic, in fact entire with exponential type not exceeding $\int_a^b\|H(t)\|\DD t$, 
and by the differential equation is has a positivity property.
Here, and always, we use the spectral norm of a matrix.

\begin{Definition}
\label{U293}
	We denote by \IndexN{\mc M_0} the set of all matrix functions $W(z)=(w_{ij}(z))_{i,j=1}^2$ whose entries are entire
	functions and take real values along the real axis, that satisfy $\det W(z)=1$ for all $z\in\bb C$, and 
	\[
		\forall z\in\bb C^+\DP \frac{W(z)JW(z)^*-J}{z-\ov z}\geq 0
		.
	\]
\end{Definition}

\begin{Theorem}
\label{U87}
\IndexxS{Theorem!direct and inverse spectral theorem!limit circle case}
	\phantom{}
	\begin{Enumerate}
	\item Let $H\in\bb H_{a,b}$ be in limit circle case. Then $W_H\in\mc M_0$. 
	\item The map $H\mapsto W_H$ induces a bijection of the quotient set of 
		\[
			\bigcup_{-\infty<a<b\leq+\infty} \Big(\bb H_{a,b}\cap L^1((a,b),\bb R^{2\times 2})\Big)
		\]
		modulo reparameterisation onto $\mc M_0$.
	\end{Enumerate}
\end{Theorem}

\noindent
We remark that, using appropriate topologies, the map $H\mapsto W_H$ becomes a homeomorphism.

\subsection[{The Weyl coefficient}]{The Weyl coefficient}
\label{U116}

Assume that $H$ is in limit point case, i.e., $H\in\bb H_{a,b}\setminus L^1((a,b),\bb R^{2\times 2})$. 
Then the fundamental solution does not extend to the right endpoint $b$. 
A result analogous to \Cref{U87} requires a different analytic object.

\begin{Definition}
\label{U89}
	Let $H\in\bb H_{a,b}$ be in limit point case. Then we define a function 
	$\IndexN{q_H}\DF\bb C^+\to\ov{\bb C^+}$ by the requirement that 
	\[
		\forall z\in\bb C^+\DP \bigcap_{t\in[a,b)}\Omega_{t,z}=\{q_H(z)\}
		.
	\]
	This function is called the \IndexS{Weyl coefficient}{Weyl coefficient} of $H$. 

	It is customary to extend $q_H$ to the lower half-plane by symmetry as
	\[
		q_H(\ov z)\DE \ov{q_H(z)},\qquad z\in\bb C^+
		.
	\]
\end{Definition}

\noindent
In some contexts the function $q_H$ is also called the \IndexS{Titchmarsh-Weyl coefficient}{Titchmarsh-Weyl coefficient} or 
the \IndexS{$m$-function}{$m$-function}. 

For each $z\in\bb C^+$ and $\tau_t\DF\bb C^+\to\ov{\bb C^+}$, $t\in[a,b)$, we have 
\[
	q_H(z)=\lim_{t\to b}\frac{w_{H,11}(t,z)\tau_t(z)+w_{H,12}(t,z)}{w_{H,21}(t,z)\tau_t(z)+w_{H,22}(t,z)}
	,
\]
and this limit is attained uniformly in $(\tau_t)_{t\in[a,b)}$ and locally uniformly in $z$. 

The Weyl coefficient is analytic and has a positivity property.

\begin{Definition}
\label{U294}
	A complex valued function $q$ is called a \IndexS{Nevanlinna function}{Nevanlinna function}, if it is defined and
	analytic in $\bb C^+$ and satisfies $q(\bb C^+)\subseteq\bb C^+\cup\bb R$.
	The \IndexS{Nevanlinna class}{Nevanlinna class} is the set of all Nevanlinna functions, and we denote it as
	\IndexN{\mc N_0}.
\end{Definition}

\noindent
In the literature it is also common to use the name \IndexS{Herglotz function}{Herglotz function} instead of Nevanlinna
function. Note that a function $q\in\mc N_0$ satisfies $q(\bb C^+)\subseteq\bb C^+$ unless it is a real constant. 
Sometimes the function that is identically equal to $\infty$ is also considered an element of $\mc N_0$. 
We shall not do so, and write $\mc N_0\cup\{\infty\}$ when we have to include the constant $\infty$. 

\begin{Theorem}
\label{U91}
\IndexxS{Theorem!direct and inverse spectral theorem!limit point case}
	\phantom{}
	\begin{Enumerate}
	\item Let $H\in\bb H_{a,b}$ be in limit point case. Then $q_H\in\mc N_0\cup\{\infty\}$. 
	\item The map $H\mapsto q_H$ induces a bijection of the quotient set of 
		\[
			\bigcup_{-\infty<a<b\leq+\infty} \Big(\bb H_{a,b}\setminus L^1((a,b),\bb R^{2\times 2})\Big)
		\]
		modulo reparameterisation onto $\mc N_0\cup\{\infty\}$. 
	\end{Enumerate}
\end{Theorem}

\noindent
We remark that, using appropriate topologies, the map $H\mapsto q_H$ becomes a homeomorphism.

\section[{The operator model}]{The operator model}
\label{U114}

References for this section are \cite{behrndt.hassi.snoo:2020,remling:2018,hassi.snoo.winkler:2000} and, historically, 
\cite{orcutt:1969,kac:1983,kac:1984}.

The operator model for a canonical system is defined mostly analogously to that of a Schr\"odinger operator. There is a maximal
and a minimal operator, which are adjoint to each other. The minimal operator is symmetric and has equal deficiency indices.
Self-adjoint realizations are then obtained by imposing boundary conditions at the endpoints of the interval $(a,b)$. Notable
differences to the Schr\"odinger case are that the Hilbert space is a weighted $L^2$-space, and that the operator is
formally defined as a linear relation which might have a multi-valued part (which is hardly problematic). 

First we define the \IndexS{model space}{operator model!model space} $L^2(H)$ associated with a Hamiltonian.

\begin{Definition}
\label{U289}
	Let $H\in\bb H_{a,b}$. Then we set 
	\[
		\IndexN{\mf L^2(H)}\DE\Biggl\{f\DF(a,b)\to\bb C^2\;\Bigg|\; 
		\raisebox{1.5ex}{\parbox[t]{46ex}{\small $f$ measurable, $\int_a^b f(t)^*H(t)f(t)\DD t<\infty$,\\[0.5ex]
		$(c,d)$ indivisible type $\phi$\ $\Rightarrow$\ $\xi_\phi^T f$ constant on $(c,d)$}}
		\Biggr\}
		,
	\]
	and let \IndexN{L^2(H)} be the factor space of $\mf L^2(H)$ modulo the equivalence relation that identifies two functions 
	$f_1,f_2$ when $Hf_1=Hf_2$ a.e.
\end{Definition}

\noindent
With the natural scalar product, $L^2(H)$ becomes a Hilbert space. 

Next we define the maximal and minimal relation associated with $H$.

\begin{Definition}
\label{U288}
	Let $H\in\bb H_{a,b}$ and assume that $H$ is definite. 
	The \IndexS{maximal relation}{operator model!maximal relation} associated with $H$ is 
	\begin{align*}
		\IndexN{T_{\max}(H)}\DE \bigg\{
		(f;g)\in L^2(H)\times L^2(H)\;\Big|\; 
		\raisebox{1.5ex}{\parbox[t]{30ex}{\small $f$ has an absolutely continuous\\[0.5ex] 
		representative $\hat{f}$ with $\hat{f}'=JHg$}}
		\bigg\}.
	\end{align*}
	The \IndexS{minimal relation}{operator model!minimal relation} associated with $H$ is
	\[
		\IndexN{T_{\min}(H)}\DE\Clos\big\{(f;g)\in T_{\max}(H)\;\big|\;\hat f\text{ has compact support in }(a,b)\big\}
		.
	\]
\end{Definition}

\noindent
The assumption that $H$ is definite guarantees that the representative $\hat f$ is unique for each $(f;g)\in T_{\max}(H)$, and
hence $T_{\min}(H)$ is well-defined.

The minimal relation is symmetric and its adjoint is the maximal relation:
\[
T_{\min}(H) \subseteq T_{\min}(H)^*=T_{\max}(H).
\]
The case distinction limit circle/limit point manifests itself in the following operator theoretic alternative:
\begin{Ilist}
\item If $H$ is in limit circle case, $T_{\min}(H)$ has deficiency indices $(2,2)$;
\item If $H$ is in limit point case, $T_{\min}(H)$ has deficiency indices $(1,1)$.
\end{Ilist}
If limit circle case takes place, the self-adjoint extensions of $T_{\min}(H)$ are obtained by imposing boundary conditions at
both endpoints of $a,b$. These boundary conditions may be separated or coupled. If limit point case holds, it is enough to
impose a boundary condition at $a$, and the self-adjoint extensions are precisely the restrictions of $T_{\max}(H)$ by a 
boundary condition of the form $(\cos\alpha,\sin\alpha)^T\hat f(a)=0$ where $\alpha\in[0,\pi)$. 

In any case, each two self-adjoint restrictions of $T_{\max}(H)$ are finite-rank perturbations of each other, and 
therefore we have:
\begin{Ilist}
\item If one self-adjoint realisation has discrete spectrum, so does every other;
\item In limit circle case, the spectrum of every self-adjoint realisation is discrete, and the difference between the number of
	eigenvalues of two self-adjoint realisations in any finite interval is at most 2;  
\item In limit point case, if spectra are discrete, then the eigenvalues of any two self-adjoint realisations interlace each 
	other.
\end{Ilist}
The large-scale distribution of eigenvalues is thus independent of the particular choice of boundary conditions. 
It is convenient to use the following self-adjoint realisation by default.

\begin{Definition}
\label{U210}
	Let $H\in\bb H_{a,b}$ and assume that $H$ is definite. We define the 
	\IndexS{model operator}{operator model!model operator} \IndexN{A_H} as follows:
	\begin{Itemize}
	\item If $H$ is in limit circle case, set 
		\[
			A_H\DE\big\{(f;g)\in T_{\max}(H)\DS (1,0)\hat f(a)=(0,1)\hat f(b)=0\big\}
			;
		\]
	\item If $H$ is in limit point case, set 
		\[
			A_H\DE\big\{(f;g)\in T_{\max}(H)\DS (1,0)\hat f(a)=0\big\}
			.
		\]
	\end{Itemize}
\end{Definition}

\noindent
As a self-adjoint linear relation, $A_H$ decomposes into the orthogonal sum of a self-adjoint
operator and the purely multi-valued relation $A_H\cap(\{0\}\times L^2(H))$. 
We write $\IndexN{\sigma(H)}$ for the spectrum of the operator part of $A_H$. 

We can now give operator theoretic meaning to monodromy matrix and Weyl coefficient. 
In line with the scope of this survey, we do not discuss spectral
representations in detail and confine our interest to the case that $\sigma(H)$ is discrete.

\begin{Theorem}
\label{U287}
\IndexxS{Theorem!spectral properties of $A_H$}
	Let $H \in\bb H_{a,b}$ and assume that $H$ is definite. 
	\begin{Itemize}
	\item If $H$ is in limit circle case, then $\sigma (H)$ is discrete and coincides with the set of zeroes of $w_{H,22}$;
	\item If $H$ is in limit point case, then $\sigma(H)$ is discrete if and only if $q_H$ has a meromorphic continuation 
		to $\bb C$ whose values along the real axis are real; 
	\item If $H$ is in limit point case and $\sigma(H)$ is discrete, then $\sigma(H)$ ($=\sigma_p(H)$) is equal to 
		the set of poles of $q_H$ and all eigenvalues are simple.
	\end{Itemize}
\end{Theorem}

\section[{Two particular classes of systems}]{Two particular classes of systems}
\label{U115}

The connection between canonical systems and the Hamburger power moment problem, presented in \Cref{U117}, is first formulated 
explicitly in \cite{kac:1999}. The formulae have been around for much longer, even in a more general setting e.g.\ as 
in \cite{krein.langer:1979,krein.langer:1980}.
For moment problems themselves there is a variety of textbooks, e.g., \cite{shohat.tamarkin:1943,akhiezer:1965,schmuedgen:2017}.

The class of Hamiltonians presented in \Cref{U118} is a rather recent invention which appeared in the context of indefinite 
inner product spaces, see \cite{kaltenbaeck.woracek:p4db,winkler.woracek:del,langer.woracek:intm,woracek:nass} and is also 
related with a notion from \cite{berg.duran:1995}.

\subsection[{Hamburger Hamiltonians: the power moment problem}]{Hamburger Hamiltonians: the power moment
	problem}
\label{U117}

In this subsection we explain the connection between three -- a priori -- different kinds of objects. 

\begin{Definition}
\label{U286}
	\phantom{}
	\begin{itemize}
	\item[$1^\circ$.] A sequence $(s_n)_{n=0}^\infty$ of real numbers is called \IndexS{positive}{positive sequence}, if 
		\[
			\forall N\in\bb N,(\xi_n)_{n=0}^N\in\bb C^{N+1}\DP
			\sum_{i,j=0}^N s_{i+j}\xi_i\ov{\xi_j}\geq 0
			.
		\]
	\item[$2^\circ$.] Assume we have two sequences, $(a_n)_{n=0}^\infty$ and $(b_n)_{n=0}^\infty$, such that $a_n\in\bb R$ and 
		$b_n>0$ for all $n\in\bb N_0$. Then we define a tridiagonal infinite matrix as 
		\[
			\IndexN{\ms J}=
			\begin{pmatrix}
				a_0 & b_0 & 0 & & & \\
				b_0 & a_1 & b_1 & 0 & &\\
				0   & b_1 & a_2 & b_2 & & \\
				    & \raisebox{5pt}[0pt][5pt]{$0$}   & \raisebox{5pt}[0pt][5pt]{$b_2$} & \ddots & \ddots & \\
				    &     &     & \ddots & \ddots & \\
			\end{pmatrix}
			,
		\]
		and refer to $\ms J$ as the \IndexS{Jacobi matrix}{Jacobi!matrix} with parameters $a_n,b_n$. 
	\item[$3^\circ$.] Assume we have two sequences, $l\DE(l_j)_{j=1}^\infty$ and $\phi\DE(\phi_j)_{j=1}^\infty$, such
		that $l_j>0$ and $\phi_j\in\bb R$ with $\phi_{j+1}-\phi_j\not\equiv 0\mod\pi$ for all $j\in\bb N$, and set 
		\[
			x_0\DE 0,\qquad x_n\DE\sum_{j=1}^n l_j,\ n\in\bb N,\qquad L\DE\sum_{j=1}^\infty l_j\in(0,\infty]
			.
		\]
		Then we define a Hamiltonian $H_{l,\phi}\DF[0,L)\to\bb R^{2\times 2}$ as 
		\[
			H_{l,\phi}(t)\DE\xi_{\phi_j}\xi_{\phi_j}^T\quad\text{for }j\in\bb N\text{ and }
			t\in[x_{j-1},x_j)
			,
		\]
		and refer to $H_{l,\phi}$ as the \IndexS{Hamburger Hamiltonian}{Hamiltonian!Hamburger} with
		\IndexS{lengths}{Hamiltonian!Hamburger!lengths} $l_j$ and \IndexS{angles}{Hamiltonian!Hamburger!angles}
		$\phi_j$. 
	\end{itemize}
\end{Definition}

\noindent
We observe that a Hamburger Hamiltonian consists of a sequence of indivisible intervals accumulating at the right endpoint 
of the interval:
\begin{center}
	\begin{tikzpicture}[x=1.2pt,y=1.2pt,scale=0.8,font=\fontsize{8}{8}]
		\draw[thick] (10,30)--(215,30);
		\draw[dotted, thick] (215,30)--(270,30);
		\draw[thick] (10,25)--(10,35);
		\draw[thick] (70,25)--(70,35);
		\draw[thick] (120,25)--(120,35);
		\draw[thick] (160,25)--(160,35);
		\draw[thick] (190,25)--(190,35);
		\draw[thick] (210,25)--(210,35);
		\draw[thick] (270,25)--(270,35);
		\draw (40,44) node {${\displaystyle \xi_{\phi_1}\xi_{\phi_1}^T}$};
		\draw (95,44) node {${\displaystyle \xi_{\phi_2}\xi_{\phi_2}^T}$};
		\draw (140,44) node {${\displaystyle \xi_{\phi_3}\xi_{\phi_3}^T}$};
		\draw (177,43) node {${\cdots}$};
		\draw (-20,30) node {\large $H_{l,\phi}\!:$};
		\draw (10,18) node {${\displaystyle 0}$};
		\draw[dashed,stealth-stealth] (11,26)--(69,26);
		\draw (40,21) node {${l_1}$};
		\draw (70,18) node {${x_1}$};
		\draw[dashed,stealth-stealth] (71,26)--(119,26);
		\draw (95,21) node {${l_2}$};
		\draw (120,18) node {${x_2}$};
		\draw[dashed,stealth-stealth] (121,26)--(159,26);
		\draw (140,21) node {${l_3}$};
		\draw (160,18) node {${x_3}$};
		\draw (195,18) node {${\cdots}$};
		\draw (270,18) node {${\displaystyle L}$};
	\end{tikzpicture}
\end{center}
Note that the lengths $l_j$ are unique for every Hamburger Hamiltonian, due to the condition that 
$\phi_{j+1}-\phi_j\not\equiv 0\mod\pi$. Each angle $\phi_j$, however, is only determined up to an integer multiple of $\pi$.
The term Hamburger Hamiltonian for this type of Hamiltonian was coined by I.S.~Kac to reference the connection we are
going to explain in the sequel. 

We start by recalling some facts about the moment problem. 
The \IndexS{Hamburger moment problem}{moment problem} is the task of describing, for a sequence 
$(s_n)_{n=0}^\infty$ of real numbers, the set 
\[
	\mc M\big((s_n)_{n=0}^\infty\big)\DE
	\bigg\{\,\mu\mkern7mu\Big|\mkern10mu
		\parbox[c]{49mm}{\small $\mu$ positive measure on $\bb R$\\ 
		$s_n=\int_{\bb R}t^n\DD\mu(t)$ for $n=0,1,2,\ldots$}
	\bigg\}
	.
\]
This problem was treated extensively in work of H.~Hamburger, M.~Riesz, R.~Nevanlinna, and many others. 
The moment problem has a solution, i.e., $\mc M((s_n)_{n=0}^\infty)\neq\emptyset$, 
if and only if the sequence $(s_n)_{n=0}^\infty$ is positive, i.e., 
all quadratic forms $\sum_{i,j=0}^ns_{i+j}\xi_i\overline{\xi_j}$ are positive semidefinite.
If so, there are two possible alternatives for the set $\mc M((s_n)_{n=0}^\infty)$:
\begin{Itemize}
\item $\mc M((s_n)_{n=0}^\infty)$ contains exactly one element -- we say that the moment problem is 
	\IndexS{determinate}{moment problem!determinate};
\item $\mc M((s_n)_{n=0}^\infty)$ has infinitely many elements -- we say that the moment problem is 
	\IndexS{indeterminate}{moment problem!indeterminate}. 
\end{Itemize}
In the indeterminate case, $\mc M((s_n)_{n=0}^\infty)$ can be described via the set of its Cauchy-transforms. 

\begin{Theorem}
\label{U282}
\IndexxS{Theorem!Nevanlinna parameterisation}
	Let $(s_n)_{n=0}^\infty$ be a positive sequence and assume that the moment problem is indeterminate. 
	Then there exist four entire functions $A,B,C,D$, such that the formula 
	\begin{equation}
	\label{B90}
		\int_{\bb R}\frac 1{t-z}\DD\mu(t)=\frac{A(z)\tau(z)-C(z)}{-B(z)\tau(z)+D(z)}
	\end{equation}
	establishes a bijection between $\mc M((s_n)_{n=0}^\infty)$ and $\mc N_0\cup\{\infty\}$.
\end{Theorem}

\noindent
The matrix 
\[
	W(z)\DE\begin{pmatrix} A(z) & C(z) \\ B(z) & D(z)\end{pmatrix}
\]
is called the \IndexS{Nevanlinna matrix}{Nevanlinna matrix} of the sequence $(s_n)_{n=0}^\infty$.

\begin{flushleft}
	\textbf{Relating $1^\circ$ and $2^\circ$.}
\end{flushleft}
Given a positive sequence $(s_n)_{n=0}^\infty$ we obtain an associated sequence $(p_n)_{n=0}^\infty$ of 
\IndexS{orthonormal polynomials}{orthonormal polynomials}. Namely, by applying the Gram-Schmidt orthonormalisation process to the
sequence $(z^n)_{n=0}^\infty$ in a space $L^2(\mu)$ where $\mu\in\mc M((s_n)_{n=0}^\infty)$. The polynomials $p_n$ 
do not depend on the choice of $\mu$. They satisfy a \IndexS{three-term recurrence}{three-term recurrence}: there exist 
$a_n\in\bb R$ and $b_n>0$ for $n\in\bb N_0$, such that 
\begin{equation}
\label{B95}
	\forall n\in\bb N_0\DP zp_n(z)=b_np_{n+1}(z)+a_np_n(z)+b_{n-1}p_{n-1}(z)
	.
\end{equation}
Here we formally set $b_{-1}\DE-1$, $p_{-1}(z)\DE 0$. 
The parameters $(a_n)_{n=0}^\infty$ and $(b_n)_{n=0}^\infty$ occurring in \cref{B95} are uniquely determined by the sequence 
$(s_n)_{n=0}^\infty$ and are called the \IndexS{Jacobi parameters}{Jacobi!parameters} of the moment sequence.

\begin{Theorem}
\label{U285}
\IndexxS{Theorem!moment sequence vs.\ Jacobi parameters}
	The assignment outlined above gives a bijective correspondence between the set of all positive sequences 
	$(s_n)_{n=0}^\infty$ with $s_0=1$ and the set of all pairs of sequences $(a_n)_{n=0}^\infty$ and $(b_n)_{n=0}^\infty$ 
	where $a_n$ are real and $b_n$ are positive.
\end{Theorem}

\noindent
The restriction in the theorem to the case that $s_0=1$ is no loss of generality: rescaling the sequence
$(s_n)_{n=0}^\infty$ by any positive factor corresponds to rescaling solutions of the moment problem by the same factor. 

\begin{Remark}
\label{U281}
	Given a moment sequence we may also define another sequence of polynomials $(q_n)_{n=0}^\infty$, called the 
	\IndexS{orthogonal polynomials of the second kind}{orthogonal polynomials!second kind}. Namely, $(q_n)_{n=-1}^\infty$ 
	is the solution of the recurrence \cref{B95} with the initial conditions (recall that we have set $b_{-1}\DE-1$)
	\[
		q_0(z)=0, \qquad q_{-1}(z)=-1
		.
	\]
\end{Remark}

\noindent
The moment problem is indeterminate if and only if the series $\sum_{n=0}^\infty(p_n(0)^2+q_n(0)^2)$ converges. 
If it is indeterminate, the entries of the Nevanlinna matrix can be expressed as series involving the polynomials 
$p_n$ and $q_n$. 

\begin{flushleft}
	\textbf{Relating $2^\circ$ and $3^\circ$.}
\end{flushleft}
Jacobi- and Hamiltonian parameters are related in a purely algebraic (yet, not simple) way. 
Assume we have sequences $(a_n)_{n=0}^\infty,(b_n)_{n=0}^\infty$ with $a_n$ real and $b_n$ positive. 
Then we define sequences $(l_j)_{j=1}^\infty,(\phi_j)_{j=1}^\infty$ with $l_j$ positive and $\phi_j$ real by 
recursively solving the equations 
\begin{align}
\label{U283}
	&
	l_1=1,\quad \phi_1=\frac\pi2,
	\\
\label{U211}
	&
	a_0=\tan \phi_2,
	\\
\label{U213}
	& \forall k\in\bb N_0\DP
	b_k=\frac{1}{\sqrt{l_{k+1}l_{k+2}}|\sin (\phi_{k+2}-\phi_{k+1})|},
	\\
\label{U212}
	& \forall k\in\bb N\DP
	a_k=\frac{\cot (\phi_{k+1}-\phi_{k+2})+\cot (\phi_{k}-\phi_{k+1})}{l_{k+1}},
\end{align}
where the numbers $\phi_j$ are determined up to integer multiples of $\pi$ (for example we could choose
$\phi_j\in[0,\pi)$). Conversely, given sequences $(l_j)_{j=1}^\infty,(\phi_j)_{j=1}^\infty$ such that
$l_j>0$, $\phi_j\in\bb R$ with $\phi_{j+1}-\phi_j\not\equiv 0\mod\pi$, and $l_1=1,\phi_1=\frac\pi2$,
the equations \cref{U211}--\cref{U212} define $(a_n)_{n=0}^\infty,(b_n)_{n=0}^\infty$ with $a_n$ real and $b_n$
positive. These constructions obviously set up a bijective correspondence between the set of all Jacobi matrices and the set of
all Hamburger Hamiltonians with $l_1=1,\phi_1=\frac\pi2$. 

The restriction to the case that $l_1=1,\phi_1=\frac\pi2$ is no essential loss of generality: appending one indivisible interval
at the initial endpoint of a Hamiltonian or removing one, respectively, are transformations that are easily understood and can
be handled explicitly. 

The correspondence introduced above also relates model operators, which we now introduce.
Given a Jacobi matrix, the \IndexS{Jacobi operator}{Jacobi!operator} is the closure in $\ell^2(\bb N_0)$ of the linear operator 
mapping elements $u$ of $\{u\in\ell^2(\bb N_0)\DS u_n=0\text{ for almost all }n\}$ to $\ms Ju$. Given a Hamburger Hamiltonian 
$H_{l,\phi}$, we define a symmetric extension of $T_{\min}(H_{l,\phi})$ as 
\[
	S_{H_{l,\phi}}\DE\Clos\big\{(f;g)\in T_{\max}(H_{l,\phi})\DS (1,0)\hat f(0)=0,\sup\Supp\hat f<L\big\}
	.
\]
\begin{Theorem}
\label{U284}
\IndexxS{Theorem!Jacobi matrices vs.\ Hamburger Hamiltonians}
	Let $\ms J$ and $H_{l,\phi}$ correspond to each other by means of the above bijection. 
	Then the Jacobi operator is unitarily equivalent to $S_{H_{l,\phi}}$.
\end{Theorem}

\noindent
The extension $S_{H_{l,\phi}}$ (and with it the corresponding Jacobi operator) is self-adjoint if $H_{l,\phi}$ is in 
limit point case, and it has deficiency indices $(1,1)$ if $H_{l,\phi}$ is in limit circle case. 
Since a Hamburger Hamiltonian is trace-normalised, limit point case takes place if and only if $L=\infty$. 

We mention that it is also possible to start with any Hamburger Hamiltonian $H_{l,\phi}$, dropping the restriction that 
$l_1=1,\phi_1=\frac\pi2$, and define Jacobi parameters by \cref{U212} and \cref{U213}. 
However, in order for the Jacobi operator to be unitarily equivalent to $S_{H_{l,\phi}}$, the formula \cref{U211} for $a_0$
has to be modified. Furthermore, the case $\phi_1\equiv 0\mod\pi$ has to be treated separately, because $S_H$ then has a
nontrivial multi-valued part.

\begin{flushleft}
	\textbf{Relating $1^\circ$ and $3^\circ$.}
\end{flushleft}
By composing the above two bijections we obtain a bijective correspondence between the set of all positive sequences
$(s_n)_{n=0}^\infty$ with $s_0=1$ and all Hamburger Hamiltonians $H_{l,\phi}$ with $l_1=1,\phi_1=\frac\pi 2$. 
There is a very interesting direct connection between $(s_n)_{n=0}^\infty$ and $H_{l,\phi}$. 

\begin{Theorem}
\label{U207}
\IndexxS{Theorem!Nevanlinna matrix vs.\ monodromy matrix}
	Let $(s_n)_{n=0}^\infty$ be a positive sequence with $s_0=1$, and let $H_{l,\phi}$ be the Hamburger Hamiltonian with 
	$l_1=1,\phi_1=\frac\pi 2$ that corresponds to the sequence $(s_n)_{n=0}^\infty$ by the above bijection. 
	\begin{Enumerate}
	\item The moment problem for $(s_n)_{n=0}^\infty$ is indeterminate, if and only if $H_{l,\phi}$ is in limit circle case. 
	\item Assume that the moment problem is indeterminate, let $W(z)$ be the Nevanlinna matrix of the moment problem and 
		$W_{H_{l,\phi}}(z)$ the monodromy matrix of $H_{l,\phi}$. Then 
		\[
			W(z)=
			\begin{pmatrix} 1 & 0 \\ 0 & -1 \end{pmatrix}
			W_{H_{l,\phi}}(z)J
			\begin{pmatrix} 1 & 0 \\ 0 & -1 \end{pmatrix}
			.
		\]
	\end{Enumerate}
\end{Theorem}

\noindent
We remark that there exist formulae expressing the parameters of the Hamiltonian corresponding to a positive sequence 
in terms of the orthogonal polynomials of that sequence:
\begin{equation}
\label{U214}
	p_n(0)=\sqrt{l_{n+1}}\sin(\phi_{n+1}),\qquad q_n(0)=-\sqrt{l_{n+1}}\cos(\phi_{n+1})
	.
\end{equation}

\subsection[{Pontryagin-type Hamiltonians: generalisation of limit circle case}]{Pontryagin-type
	Hamiltonians: generalisation of limit circle case}
\label{U118}

In this section we introduce a class of Hamiltonians that are in limit point case, yet behave very similarly to Hamiltonians 
in limit circle case. 

Before we can give the definition of said class, we need to carry out a preliminary discussion. 

\begin{Definition}
\label{U191}
	Let $H\in\bb H_{a,b}$. Then we set 
	\begin{align*}
		& \mc D\DE\big\{f\DF(a,b)\to\bb C^2\DS \forall c\in(a,b)\DP f|_{(a,c)}\in L^2(H|_{(a,c)})\big\}
		,
		\\
		& (\IndexN{V_H}f)(t)\DE\int_a^t JH(s)f(s)\DD s,\quad t\in[a,b),f\in\mc D
		.
	\end{align*}
\end{Definition}

\noindent
Note that $L^2(H)\subseteq\mc D$ and $\bb C^2\subseteq\mc D$, where we understand $\bb C^2$ as set of constant functions. 

In the following we work with functions of the form $\sum_{j=0}^nV_H^ja_j$ where $n\in\bb N$ and $a_0,\ldots,a_n\in\bb C^2$. We
call such functions \IndexS{$H$-polynomials}{$H$-polynomial}. 
This terminology is chosen by analogy with the case when $V$ is the classical Volterra operator $(Vf)(t)\DE\int_0^tf(s)\DD s$. 

\begin{Definition}
\label{U190}
	Let $H\in\bb H_{a,b}$. Then we set
	\begin{align*}
		& 
		C_n\DE\Big\{a_n\in\bb C^2\DS \exists a_0,\ldots,a_{n-1}\in\bb C^2\DP\sum_{j=0}^nV_H^ja_j\in L^2(H)\Big\}
		,\qquad n\in\bb N,
		\\
		& 
		\IndexN{\Delta(H)}\DE\inf\big\{n\in\bb N\DS \dim C_n=2\big\}\in\bb N\cup\{\infty\}
		.
	\end{align*}
\end{Definition}

\noindent
Note that $H$ is in limit circle case if and only if $\Delta(H)=0$. 

\begin{Lemma}
\label{U189}
	Let $H\in\bb H_{a,b}$ and assume that 
	\begin{equation}
	\label{U188}
		\begin{aligned}
			\exists \phi\in[0,\pi)\DP &\, \xi_\phi\in L^2(H)\ \wedge 
			\\
			&\, \limsup_{t\to b}\bigg(
			\int_t^b\xi_\phi^TH(s)\xi_\phi\DD s\cdot\int_a^t\xi_{\phi+\frac\pi 2}^TH(s)\xi_{\phi+\frac\pi 2}\DD s
			\bigg)<\infty
			.
		\end{aligned}
	\end{equation}
	Then $V_H$ maps $L^2(H)$ boundedly into itself. 
\end{Lemma}

\noindent
We remark that the condition \cref{U188} means that $0\notin\sigma_{\Ess}(H)$, cf.\ \Cref{U2}. 

This lemma implies that, under the condition \cref{U188}, we have 
\[
	\forall n\in\bb N\DP C_n\subseteq C_{n+1}
	.
\]

\begin{Definition}
\label{U187}
	Let $H\in\bb H_{a,b}$. We say that $H$ is of 
	\IndexS{Pontryagin type}{Pontryagin type}\index[sub]{Hamiltonian!Pontryagin type}, if 
	\begin{Enumerate}
	\item \raisebox{-16pt}{
		${\displaystyle
		\begin{aligned}
			\exists \phi\in[0,\pi)\DP &\, \xi_\phi\in L^2(H)\ \wedge
			\\
			&\, \lim_{t\to b}\bigg(
			\int_t^b\xi_\phi^TH(s)\xi_\phi\DD s\cdot\int_a^t\xi_{\phi+\frac\pi 2}^TH(s)\xi_{\phi+\frac\pi 2}\DD s
			\bigg)=0,
		\end{aligned}
		}$
		}
	\item $\Delta(H)<\infty$. 
	\end{Enumerate}
\end{Definition}

\noindent
We remark that the condition (i) means that $\sigma_{\Ess}(H)=\emptyset$, cf.\ \Cref{U4}.

The class of Pontryagin type Hamiltonians can be seen as a mild generalisation of limit circle case Hamiltonians. This, and 
our choice of terminology, is motivated by the fact that a Hamiltonian is of Pontryagin type if and only if it is 
a section of an indefinite Hamiltonian in the sense of \cite{kaltenbaeck.woracek:p4db} that is in limit circle case. Such
indefinite Hamiltonians give rise to an operator model in a Pontryagin space instead of a Hilbert space. 

\begin{Example}
\label{U186}
	Let $\alpha\in\bb R$ and consider the Hamiltonian $H\in\bb H_{0,1}$ defined as 
	\[
		H(t)\DE\begin{pmatrix} 1 & 0\\ 0 & \big(\frac 1{1-t}\big)^\alpha\end{pmatrix}, \quad t\in(0,1)
		.
	\]
	Then 
	\begin{Itemize}
	\item $\alpha<1$: limit circle case.
	\item $\alpha\in[2-\frac 1n,2-\frac 1{n+1})$, $n\geq 1$: Pontryagin type with $\Delta(H)=n$.
	\item $\alpha=2$: $0\notin\sigma_{\Ess}(H)\neq\emptyset$ and $\Delta(H)=\infty$.
	\item $\alpha>2$: $0\in\sigma_{\Ess}(H)$. 
	\end{Itemize}
\end{Example}

\noindent
Pontryagin type Hamiltonians can be characterised in terms of their Weyl coefficients. The following is shown in 
\cite[Theorem~5.1]{langer.woracek:intm}.

\begin{Theorem}
\label{U185}
\IndexxS{Theorem!characterisation of Pontryagin type}
	Let $H\in\bb H_{a,b}$ be in limit point case. Then $H$ is of Pontryagin type, if and only if its Weyl coefficient $q_H$
	has the following properties:
	\begin{Enumerate}
	\item $q_H$ has a meromorphic extension to the whole plane $\bb C$.
	\item Denote by $P$ the set of poles of $q_H$, and set 
		\[
			n_+(r)\DE\#\big[P\cap(0,r)\big],\quad n_-(r)\DE\#\big[P\cap(-r,0)\big]
			.
		\]
		Then the limits 
		\[
			\lim_{r\to\infty}\sum_{\substack{w\in P\\ |w|\leq r}}\frac 1w,\qquad
			\lim_{r\to\infty}\frac{n_+(r)}r,\quad
			\lim_{r\to\infty}\frac{n_-(r)}r
			,
		\]
		exist in $\bb R$, and the last two limits are equal. 
	\item Denote by $p$ the entire function defined as 
		\[
			p(z)\DE
			\begin{cases}
				\lim\limits_{r\to\infty}\prod_{\substack{w\in P\\ |w|\leq r}}\big(1-\frac zw\big) 
				\CAS 0\notin P,
				\\[4mm]
				z\lim\limits_{r\to\infty}\prod_{\substack{w\in P\setminus\{0\}\\ |w|\leq r}}
				\big(1-\frac zw\big) 
				\CAS 0\in P.
			\end{cases}
		\]
		Moreover, for $w\in P$, let $c_w$ be minus the residuum of $q_H$ at $w$. Then there exists $\Delta\in\bb N$ such
		that 
		\begin{equation}
		\label{U267}
			\sum_{w\in P}\frac 1{(1+w^{2(1+\Delta)})p'(w)c_w}<\infty
			.
		\end{equation}
	\end{Enumerate}
	If $q_H$ satisfies (i)--(iii), then $\Delta(H)$ is the minimum of all numbers $\Delta\in\bb N$ that satisfy \cref{U267}.
\end{Theorem}

%**************************************************************************
%***                            Last Change: Mon 31 Mar 2025 11:28
%***   < PART I >
%***
%**************************************************************************

\clearpage
\PART{I}{Fundamental Theorems}
\label{U126}

%%%%%%%%%%%%%%%%%%%%%%%%%%%%%%%%%%%%%%%%%%%%%%%%%%%%%%%%%%%%%%%%%%%%%%%%%%%%%%%%%%%%%%%%%

\Intro{%
	In this part we present the fundamental theorems that determine, explicitly in terms of $H$,
	\begin{itemize}
	\item[(1)] whether the spectrum is discrete, and 
	\item[(2)] if it is discrete, how large its density is.
	\end{itemize}
	We understand ``density'' in the sense familiar from complex analysis, namely, we ask for convergence of series or 
	finiteness of limit superior w.r.t.\ to comparison functions $\ms g$ 
	(the model case being $\ms g(r)=r^\rho$ with some $\rho>0$), and also use corresponding wording:
	\begin{Itemize}
	\item \IndexS{Convergence class}{convergence class}:
		$\sum\limits_{\lambda\in\sigma(H)\setminus\{0\}}\frac 1{\ms g(|\lambda|)}<\infty$;
	\item \IndexS{Finite (or minimal-) type}{finite/minimal type}: 
		$\limsup\limits_{r\to\infty}\frac{n_H(r)}{\ms g(r)}<\infty$ (or $=0$, respectively), where 
		\[
			n_H(r)\DE\#\{\lambda\in\sigma(H)\DS|\lambda|<r\}
			.
		\]
	\end{Itemize}
	We note that every discrete subset of $\bb R$ can be realised as $\sigma(H)$ for some $H$.

	Question (1) is answered by \Cref{U2}. The given condition is surprisingly simple. 
	Question (2) has two facets: the case that $\sigma(H)$ has large density and the case that $\sigma(H)$ is sparse. 
	Intuitively, though not fully correct, large density means that $\sigma(H)$ is ``more dense that the integers'', 
	while sparse means roughly ``at most as dense as $\bb Z$''.
	The case of dense spectrum is settled by \Cref{U37,U39}. The given conditions are explicit and accessible from a 
	computational viewpoint. The case of sparse spectrum is settled by \Cref{U55}. The formula is still explicit, yet 
	hard to evaluate.

	A major result on the way to understand spectrum with large density is \Cref{U35}, that contains an astonishing
	independence property and emphasises the role of integer distribution as a
	borderline from an operator theoretic viewpoint.
}
\begin{center}
	{\large\bf Table of contents}
\end{center}
\begin{flushleft}
	\S\,\ref{U108}.\ The discreteness criterion\ \dotfill\quad\pageref{U108}
	\\[1mm]
	\S\,\ref{U100}.\ Independence from off-diagonal via the Matsaev property\ \dotfill\quad\pageref{U100}
	\\[1mm]
	\S\,\ref{U101}.\ Spectrum with large density: approach via operator ideals\ \dotfill\quad\pageref{U101}
	\\[1mm]
	\S\,\ref{U110}.\ Trace class and sparse spectrum: the Weyl coefficient approach\ \dotfill\quad\pageref{U110}
	\\[1mm]
\end{flushleft}
\makeatother
\renewcommand{\thesection}{\arabic{section}}
\renewcommand{\thelemma}{\arabic{section}.\arabic{lemma}}
\makeatletter
\setcounter{section}{0}
\clearpage

%%%%%%%%%%%%%%%%%%%%%%%%%%%%%%%%%%%%%%%%%%%%%%%%%%%%%%%%%%%%%%%%%%%%%%%%%%%%%%%%%%%%%%%%%

%
%
%
\section[{The discreteness criterion}]{The discreteness criterion}
\label{U108}

For a Hamiltonian $H\in\bb H_{a,b}$ we denote the radius of the possible gap around $0$ in the essential spectrum as 
\[
	\IndexN{R_H}\DE\inf\big\{|t|\DS t\in\sigma_{\rm ess}(H)\big\}\in[0,\infty]
	.
\]
Observe that $\sigma(H)$ is discrete if and only if $\sigma_{\rm ess}(H)$ has a gap around zero with infinite radius, i.e., 
$R_H=\infty$.

The below theorem contains a two-sided estimate for $R_H$. 
Before we formulate this result, we have to note the following fact.

\begin{Lemma}
\label{U3}
	Let $H\in\bb H_{a,b}$. If $R_H>0$, then there exists $\phi\in\bb R$ such that 
	\[
		\int_a^b\xi_\phi^TH(t)\xi_\phi\DD t<\infty
		.
	\]
\end{Lemma}

\begin{Theorem}
\label{U2}
\IndexxS{Theorem!essential spectrum}
	There exist constants $\gamma_+,\gamma_->0$ such that the following statement holds. 
	\begin{Itemize}
	\item Let $\phi\in\bb R$ and $H\in\bb H_{a,b}$, assume that $\int_a^b\xi_\phi^TH(t)\xi_\phi\DD t<\infty$, and 
		set 
		\begin{equation}
		\label{U280}
			\alpha\DE\limsup_{t\to b}\bigg(
			\int_t^b\xi_\phi^TH(s)\xi_\phi\DD s\cdot\int_a^t\xi_{\phi+\frac\pi 2}^TH(s)\xi_{\phi+\frac\pi 2}\DD s
			\bigg)\in[0,\infty]
			.
		\end{equation}
		Then%
		\/\footnote{Here we use the usual conventions for algebra in $[0,\infty]$: 
			$\frac 10=\infty$ and $\frac 1\infty=0$.}
		\[
			\gamma_-\frac 1{\sqrt\alpha}\leq R_H\leq\gamma_+\frac 1{\sqrt\alpha}
			.
		\]
	\end{Itemize}
\end{Theorem}

\noindent
In particular, we have the following \IndexS{discreteness criterion}{discreteness criterion}.

\begin{Corollary}
\label{U4}
\IndexxS{Theorem!discreteness criterion}
	Let $H\in\bb H_{a,b}$. Then $\sigma(H)$ is discrete, if and only if there exists $\phi\in\bb R$ such that 
	\begin{align}
	\label{U5}
		& \int_a^b\xi_\phi^TH(t)\xi_\phi\DD t<\infty
		,
		\\
	\label{U6}
		& \lim_{t\to b}\bigg(
		\int_t^b\xi_\phi^TH(s)\xi_\phi\DD s\cdot\int_a^t\xi_{\phi+\frac\pi 2}^TH(s)\xi_{\phi+\frac\pi 2}\DD s
		\bigg)=0
		.
	\end{align}
\end{Corollary}

\noindent
Observe that the conditions \cref{U5} and \cref{U6} depend only on the ``diagonal entries'' 
$\xi_\phi^TH\xi_\phi$ and $\xi_{\phi+\frac\pi 2}^TH\xi_{\phi+\frac\pi 2}$ of $H$, while the ``off-diagonal'' 
$\xi_\phi^TH\xi_{\phi+\frac\pi 2}$ does not occur in the formulae. 

Indeed it is a crucial step in the proof to show this \IndexS{independence theorem}{independence theorem},
where we mean independence from the off-diagonal. We again state a quantitative version.

\begin{Theorem}
\label{U7}
\IndexxS{Theorem!independence theorem}
	There exist constants $\tilde\gamma_+,\tilde\gamma_->0$ such that the following statement holds. 
	\begin{Itemize}
	\item Let $\phi\in\bb R$ and $H\in\bb H_{a,b}$, assume that $\int_a^b\xi_\phi^TH(t)\xi_\phi\DD t<\infty$, and denote 
		\begin{equation}
		\label{U69}
			\IndexN{H_d(t)}\DE
			\begin{pmatrix}
				\xi_\phi^TH(t)\xi_\phi & 0
				\\
				0 & \xi_{\phi+\frac\pi 2}^TH(t)\xi_{\phi+\frac\pi 2}
			\end{pmatrix}
			\quad\text{for }t\in(a,b)
			.
		\end{equation}
		Then
		\[
			\tilde\gamma_-R_{H_d}\leq R_H\leq\tilde\gamma_+ R_{H_d}
			.
		\]
	\end{Itemize}
\end{Theorem}

\noindent
To make it explicit, here is the formulation for discreteness only.

\begin{Corollary}
\label{U8}
\IndexxS{Theorem!independence theorem}
	Let $\phi\in\bb R$ and $H\in\bb H_{a,b}$, and assume that $\int_a^b\xi_\phi^TH(t)\xi_\phi\DD t<\infty$. 
	Moreover, let $H_d$ be as in \cref{U69}. Then $\sigma(H)$ is discrete, if and only if $\sigma(H_d)$ is discrete. 
\end{Corollary}

\REMARKS{%
\item For the case of strings (equivalently, diagonal Hamiltonians, cf.\ \Cref{U124}), 
	the discreteness criterion was shown already at a 
	very early stage by I.S.~Kac and M.G.~Krein in \cite{kac.krein:1958}, see also \cite[11.9$^\circ$]{kac.krein:1968}. 
	After a long period of no progress, a necessary
	condition and a (different) sufficient condition for discreteness was announced in \cite{kac:1995}. Proofs have never
	been published, and the reason may be that there was a flaw in the argument (cf.\ the discussion in 
	\cite[Appendix~B]{romanov.woracek:ideal}). 

	The final solution is given by R.~Romanov and H.~Woracek in \cite{romanov.woracek:ideal}, where \Cref{U8} and \Cref{U4} 
	are established using operator theoretic methods. In that paper constants were not traced. The quantitative versions
	\Cref{U7} and \Cref{U2} are taken from the work of C.~Remling and K.~Scarbrough \cite{remling.scarbrough:2020a} where a
	different approach via oscillation theory is used. 
\item In \cite{remling.scarbrough:2020a} numerical values for the constants $\gamma_\pm,\tilde\gamma_\pm$ are given: 
	\[
		\gamma_-=\frac 14,\ \tilde\gamma_-=\frac 12,\ \gamma_+=\tilde\gamma_+=\frac 2{3-\sqrt 5}
		.
	\]
	The lower bounds in \Cref{U2,U7} with these constants $\gamma_-,\tilde\gamma_-$ are sharp, 
	see the example in \cite[\S5]{remling.scarbrough:2020a} where $\alpha=1,R_H=\frac 14,R_{H_d}=\frac 12$. 
	The stated value for $\gamma_+$ and $\tilde\gamma_+$ is almost certainly not optimal, and the optimal value is not
	known.
\item For diagonal Hamiltonians a more accurate estimate can be given. Assume $H$ is diagonal, let $\alpha$ be the 
	limit superior from \cref{U280}, and denote by $\beta$ the limit inferior of the same expression. Then 
	\[
		\frac 1{2\sqrt\alpha}\leq R_H\leq\frac 1{2\sqrt\beta}
		.
	\]
	In particular, if the expression from \cref{U280} has a limit, the numerical value of $R_H$ is determined. 
\item For the proof of the discreteness criterion one could, after having established the independence theorem, refer to 
	the mentioned connection with strings and \cite{kac.krein:1958}. However, the proofs in \cite{romanov.woracek:ideal} and 
	\cite{remling.scarbrough:2020a} do not do this, and hence yield new proofs for the Kac-Krein theorem. 

	In this context let us mention that the two cases \cite[(0.8), (0.9)]{kac.krein:1958} correspond to different ways of
	rewriting the string to a diagonal Hamiltonian, and to having ``$\phi=0$'' or ``$\phi=\frac\pi 2$'' in the above
	formulations. 
\item The fact that $H$ is dominated by $H_d$ in \Cref{U7} (meaning that $R_H\geq\text{const.}\cdot R_{H_d}$)
	is not hard to see. The difficult, and rather surprising, part is that also a converse inequality holds. 
\item Hamiltonians $H$ for which $\sigma(H)$ is discrete are related with 
	\IndexS{structure Hamiltonians}{structure Hamiltonian}, a notion which occurs in the theory of de~Branges' Hilbert 
	spaces of entire functions. See e.g.\ \cite{woracek:dbgrr} for this notion and \cite{kac:2007} for the relation. 
	The question how to characterise those Hamiltonians that are the structure Hamiltonian of some de~Branges space 
	is posed by L.~de~Branges in \cite{debranges:1968} as an ``important problem''. By means of the relation made explicit by
	I.S.~Kac, it is equivalent to characterise discreteness of $\sigma(H)$. 
\item Based on the discreteness criterion one can -- in theory -- decide for all types of equations that can be reformulated as
	a canonical system whether their spectrum is discrete. However, despite that the condition in terms of $H$ given in 
	\Cref{U2} is so simple, it seems very hard (if not impossible) to express it in terms of the data of other equations
	without imposing restrictions on the data like smoothness or regularity.
	For instance, rewriting it in terms of Jacobi parameters or a moment sequence seems to be out of reach. 
\item \Cref{U3} is folklore; an explicit argument can be found in \cite[Theorem~3.8(b)]{remling:2018} or 
	\cite[Lemma~6.1]{romanov.woracek:ideal}. 

	In the context of this lemma, we remark that for almost all purposes one can
	restrict w.l.o.g.\ to the case that the ``integrable direction'' $\phi$ is one specific direction, e.g., to the case
	$\phi=0$ which amounts to $\int_a^b h_1(t)\DD t<\infty$. This follows since rotating the Hamiltonian by some angle is a
	simple transformation: let $\psi\in\bb R$ and 
	\[
		\tilde H(t)\DE
		\smmatrix{\cos\psi}{\sin\psi}{-\sin\psi}{\cos\psi}H(t)\smmatrix{\cos\psi}{\sin\psi}{-\sin\psi}{\cos\psi}^{-1}
		,
	\]
	then the model operator $A_{\tilde H}$ is unitarily equivalent to a one-dimensional perturbation of $A_H$. 
}

\section[{Independence from off-diagonal via the Matsaev property}]{Independence from off-diagonal via 
	the Matsaev property}
\label{U100}

The \IndexS{Calkin correspondence}{Calkin correspondence} is the map that assigns to a compact operator on some Hilbert 
space $\mc H$ the sequence of its $s$-numbers
\[
	\IndexN{\mc C}
	\FD{\mf S_\infty}{c_0}{0}{T}{(s_n(T))_{n=1}^\infty}
\]
Here we denote by $\IndexN{\mf S_\infty}$ the ideal of all compact operators, and $s_n(T)$ is the 
$n$-th \IndexS{$s$-number}{$s$-number} (in our context equivalently, the $n$-th 
\IndexS{approximation number}{approximation number}) 
\[
	s_n(T)\DE\inf\big\{\|T-F\|\DS F\in\mc B(\mc H),\dim\Ran F<n\big\}
	.
\]
Provided that $\mc H$ is separable, every proper ideal in the ring $\mc B(\mc H)$ of all bounded linear operators on $\mc H$ 
is contained in $\mf S_\infty(\mc H)$, and hence corresponds to a certain space of sequences. 
For a self-adjoint operator $A$ with discrete spectrum, the Calkin correspondence thus translates properties of the distribution 
of the spectrum into membership in operator ideals of resolvents of $A$. 

\begin{Example}
\label{U10}
	Let $A$ be a self-adjoint operator with discrete spectrum, and let $p\in(0,\infty)$. Then (for $z$ being any point of the
	resolvent set)
	\[
		\sum_{\lambda\in\sigma(A)\setminus\{0\}}\frac 1{|\lambda|^p}<\infty
		\ \Longleftrightarrow\ 
		(A-z)^{-1}\in\mf S_p
	\]
	Here $\IndexN{\mf S_p}$ denotes the \IndexS{Schatten--von~Neumann ideal}{Schatten--von~Neumann ideal}, 
	which is defined as the inverse image under $\mc C$ of the sequence space $\ell^p$. 

	We see that the \IndexS{convergence exponent}{convergence exponent}
	\[
		\rho\DE
		\inf\Big\{p>0\DS \sum_{\lambda\in\sigma(A)\setminus\{0\}}\frac 1{|\lambda|^p}<\infty\Big\}\in[0,\infty]
	\]
	of $\sigma(A)$ can be expressed in terms of the Schatten--von~Neumann ideals $\mf S_p$ as 
	\[
		\rho=\inf\big\{p>0\DS (A-z)^{-1}\in\mf S_p\big\}
		.
	\]
\end{Example}

\noindent
For $p\in[1,\infty)$, the ideals $\mf S_p$ carry another structure besides being merely operator ideals: $\mf S_p$ is 
a Banach space with the norm carried over from $\ell^p$, and this norm enjoys additional algebraic properties. 
The following notion formalises this observation. 

\begin{Definition}
\label{U13}
	Let $\mc H$ be a separable Hilbert space. A \IndexS{symmetrically normed ideal}{symmetrically normed ideal} 
	(\IndexS{s.n.-ideal}{s.n.-ideal} for short) in $\mc H$
	is an ideal $\mf J$ of $\mc B(\mc H)$ with $\mf J\notin\{\{0\},\mc B(\mc H)\}$ that is endowed with a norm 
	$\|\cdot\|_{\mf J}$, such that 
	\begin{Enumerate}
	\item $\mf J$ is complete with $\|\cdot\|_{\mf J}$, 
	\item ${\displaystyle
		\forall T\in\mf J,A,B\in\mc B(\mc H)\DP \|ATB\|_{\mf J}\leq\|A\|\cdot\|T\|_{\mf J}\cdot\|B\|
		}$,
	\item $\|T\|_{\mf J}=\|T\|$ for all $T\in\mc B(\mc H)$ with $\dim\Ran T=1$.
	\end{Enumerate}
\end{Definition}

\noindent
A property that many, but not all, s.n.-ideals have is the following.

\begin{Definition}
\label{U33}
	Let $\mf J$ be a s.n.-ideal with $\mf J\neq\mf S_\infty$. We say that $\mf J$ has the 
	\IndexS{Matsaev property}{Matsaev property}, if 
	\begin{equation}
	\label{U34}
		\forall T\in\mf S_\infty,\sigma(T)=\{0\}\DP T+T^*\in\mf J\Longrightarrow T\in\mf J
		.
	\end{equation}
\end{Definition}

\noindent
This choice of terminology is motivated by a classical result due to V.I.~Matsaev, who showed that all ideals $\mf S_p$ with 
$1<p<\infty$ satisfy \cref{U34} (even with a norm estimate). We should point out that the 
\IndexS{trace class ideal}{trace class ideal} $\mf S_1$ does not have the Matsaev property. 

The second part of the following theorem is a crucial -- and quite unexpected -- result. 

\begin{Theorem}
\label{U35}
	Let $\phi\in\bb R$ and $H\in\bb H_{a,b}$, and assume that $\int_a^b\xi_\phi^TH(t)\xi_\phi\DD t<\infty$. 
	Denote again
	\[
		\IndexN{H_d(t)}\DE
		\begin{pmatrix}
			\xi_\phi^TH(t)\xi_\phi & 0
			\\
			0 & \xi_{\phi+\frac\pi 2}^TH(t)\xi_{\phi+\frac\pi 2}
		\end{pmatrix}
		,\quad t\in(a,b)
		.
	\]
	Then the following statements hold.
	\begin{Itemize}
	\item \IndexS{Diagonal dominance:}{diagonal dominance}\IndexxS{Theorem!diagonal dominance}
		Let $\mf J$ be an operator ideal. Then (for $z\notin\sigma(H)\cup\sigma(H_d)$)
		\[
			(A_{H_d}-z)^{-1}\in\mf J\ \Longrightarrow\ (A_H-z)^{-1}\in\mf J
		\]
	\item \IndexS{Independence theorem:}{independence theorem}\IndexxS{Theorem!independence theorem} 
		Let $\mf J$ be an s.n.-ideal with the Matsaev property. Then 
		(for $z\notin\sigma(H)\cup\sigma(H_d)$)
		\[
			(A_{H_d}-z)^{-1}\in\mf J\ \Longleftrightarrow\ (A_H-z)^{-1}\in\mf J
		\]
	\end{Itemize}
\end{Theorem}

\noindent
The following example demonstrates validity and failure of independence from the off-diagonal. 

\begin{Example}
\label{U36}
	Let $\alpha>0$ and set 
	\[
		m(t)\DE \frac 1{1-t}\Big(1+\log\frac 1{1-t}\Big)^{-\alpha},\quad t\in(0,1)
		.
	\]
	Consider the Hamiltonian on $(0,1)$ defined as 
	\[
		H(t)\DE
		\begin{pmatrix}
			1 & -m(t)
			\\
			-m(t) & m(t)^2
		\end{pmatrix}
		.
	\]
	Then $\sigma(H)$ and $\sigma(H_d)$ are discrete, and the respective convergence exponents $\rho_H$ and $\rho_{H_d}$ are 
	\[
		\rho_H=
		\begin{cases}
			\frac 1\alpha \CAS \alpha\in(0,2),
			\\
			\frac 12 \CAS \alpha\geq 2,
		\end{cases}
		\qquad
		\rho_{H_d}=
		\begin{cases}
			\frac 1\alpha \CAS \alpha\in(0,1),
			\\
			1 \CAS \alpha\geq 1.
		\end{cases}
	\]
	Moreover, we always have $(A_{H_d}-z)^{-1}\notin\mf S_1$. 
\end{Example}

\REMARKS{%
\item The Calkin correspondence, identifying ideals with sequence spaces, goes back to J.W.~Calkin \cite{calkin:1941} and 
	D.J.H.~Garling \cite{garling:1967}. An elaborate account on the topic can be found in the textbook 
	\cite{gohberg.krein:1969}, see also \cite{simon:2005}. Those normed sequence spaces that correspond to s.n.-ideals have
	been identified only comparatively recently by N.J.~Kalton and F.A.~Sukochev in \cite{kalton.sukochev:2008}. 
\item Matsaev's theorem for the ideals $\mf S_p$ is shown in \cite{matsaev:1961}, see also 
	\cite[Theorem~6.2]{gohberg.krein:1970}. A characterisation of the Matsaev property for a large class of s.n.-ideals is
	given by G.I.~Russu in \cite{russu:1979}, \cite{russu:1980}, see also \cite[Theorem~2.5]{romanov.woracek:ideal}. 
	A full characterisation in terms of the sequence space
	can be obtained from work of F.A.~Sukochev, K.~Tulenov, D.~Zanin \cite{sukochev.tulenov.zanin:2019}. 

	As a rule of thumb one may say that the Matsaev property holds for an s.n.-ideal $\mf J$ unless $\mf J$ is too close to
	$\mf S_1$ or to $\mf S_\infty$, or the norm $\|\cdot\|_{\mf J}$ behaves in a weird way. 
\item The independence theorem was proved by R.~Romanov and H.~Woracek in \cite[Theorem~3.6]{romanov.woracek:ideal}, 
	even under a slightly weaker assumption on $\mf J$ than stated above. The proof is based on an operator theoretic trick. 
	Diagonal dominance, proved in \cite[Theorem~3.4]{romanov.woracek:ideal}, is an expected result and was implicitly 
	present in earlier literature. 
\item \Cref{U36} is \cite[Example~1.7]{romanov.woracek:ideal}. Its proof requires the results from \Cref{U101}.
\item The independence theorem sheds light on the specialty of integer distribution from an operator theoretic perspective. 
	Namely, Schatten--von~Neumann classes $\mf S_p$ are s.n.-ideals only if $p\geq 1$ and for $p=1$ the Matsaev property
	fails. 

	Compare this with the Krein-de~Branges formula further below, which also emphasises specialty of integer distribution,
	but from a function theoretic perspective.
}

\section[{Spectrum with large density: approach via operator ideals}]{Spectrum with large density: 
	approach via operator ideals}
\label{U101}

Even though the results in this section hold true for Hamiltonians in either limit circle or limit point case, they are
mostly relevant for limit point Hamiltonians. The reason for that is that they rely on \Cref{U35}\,{\rm(ii)}. 
The \IndexS{Krein--de~Branges formula}{Krein--de~Branges formula} (discussed in detail in \Cref{U109}) states that for a
Hamiltonian $H$ in limit circle case the eigenvalue counting function 
\[
	\IndexN{n_H(r)}\DE\#\big\{\lambda\in\sigma(H)\DS|\lambda|<r\big\}
\]
satisfies
\[
	\lim_{r\to\infty}\frac{n_H(r)}r=\frac 2\pi\int_a^b\sqrt{\det H(t)}\DD t<\infty
	.
\]
Thus, if $H$ is in limit circle case, membership in s.n.-ideals with the Matsaev property almost always come for free.

Due to the independence theorem the
task is to decide membership of resolvents in an operator ideal for a diagonal Hamiltonian. This can be done by exploiting a
method that goes back, at least, to the paper \cite{aleksandrov.janson.peller.rochberg:2002} of A.B.~Aleksandrov, S.~Janson, 
V.V.~Peller, R.~Rochberg. In order to obtain meaningful results on the spectral side, we apply the general theory with certain 
\IndexS{Lorentz ideals}{Lorentz ideal} and \IndexS{Orlicz spaces}{Orlicz space}. 

The first theorem deals with convergence class conditions.

\begin{Theorem}
\label{U37}
\IndexxS{Theorem!convergence class!dense spectrum}
	Let $H\in\bb H_{a,b}$ with $\sigma(H)$ discrete, let $\phi\in\bb R$ be such that%
	\/\footnote{Such an angle always exists, cf.\ \Cref{U3}.}
	$\int_a^b\xi_\phi^TH(t)\xi_\phi\DD t<\infty$, and assume%
	\/\footnote{We will explain in the remarks below that this is no loss in generality.}
	that $\xi_\phi^TH(t)\xi_\phi$ does not vanish a.e.\ on any interval $(c,b)$ with $c<b$.
	Moreover, let $\ms g$ be a regularly varying function with index $>1$ (cf.\ \Cref{U123}). 
	Then, setting $\psi\DE\phi+\frac\pi 2$, we have 
	\begin{multline}
	\label{U40}
		\sum_{\lambda\in\sigma(H)\setminus\{0\}}\frac 1{\ms g(|\lambda|)}<\infty\quad \Longleftrightarrow
		\\[2mm]
		\int\limits_a^b\bigg(\ms g\bigg[\Big(
		\int\limits_t^b\xi_\phi^TH(s)\xi_\phi\DD s\cdot\int\limits_a^t\xi_\psi^TH(s)\xi_\psi\DD s
		\Big)^{-\frac 12}\bigg]\bigg)^{-1}
		\frac{\xi_\phi^TH(s)\xi_\phi}{\int\limits_t^b\xi_\phi^TH(s)\xi_\phi\DD s}\DD t
		<\infty
	\end{multline}
\end{Theorem}

\noindent
Usual \IndexS{convergence class w.r.t.\ an order}{convergence class!w.r.t.\ order} $\rho$ is the particular case 
$\ms g(r)=r^\rho$. In this case the integrand in \cref{U40} simplifies. 

\begin{Corollary}
\label{U38}
	Let $H\in\bb H_{a,b}$ with $\sigma(H)$ discrete, let $\phi\in\bb R$ be such that 
	$\int_a^b\xi_\phi^TH(t)\xi_\phi\DD t<\infty$,
	and assume that $\xi_\phi^TH(t)\xi_\phi$ does not vanish a.e.\ on any interval $(c,b)$ with $c<b$. 
	Moreover, let $\rho>1$. 
	Then, setting $\psi\DE\phi+\frac\pi 2$, we have 
	\begin{multline}
	\label{U41}
		\sum_{\lambda\in\sigma(H)\setminus\{0\}}\frac 1{|\lambda|^\rho}<\infty\quad \Longleftrightarrow
		\\[2mm]
		\int\limits_a^b
		\Big(\int\limits_t^b\xi_\phi^TH(s)\xi_\phi\DD s\Big)^{\frac\rho 2-1}
		\Big(\int\limits_a^t\xi_\psi^TH(s)\xi_\psi\DD s\Big)^{\frac\rho 2}\cdot
		\xi_\phi^TH(s)\xi_\phi\DD t
		<\infty
	\end{multline}
\end{Corollary}

\noindent
In the second theorem we deal with finite type conditions and 
minimal type conditions. This case is more complicated concerning
presentation (not concerning content), since the respective conditions cannot anymore be written in integral form as 
\cref{U40} and \cref{U41}. They have a sequential (instead of integral-) form. 

\begin{Theorem}
\label{U39}
\IndexxS{Theorem!finite/minimal type!dense spectrum}
	Let $H\in\bb H_{a,b}$ with $\sigma(H)$ discrete, let $\phi\in\bb R$ be such that 
	$\int_a^b\xi_\phi^TH(t)\xi_\phi\DD t<\infty$, and assume that $\xi_\phi^TH(t)\xi_\phi$ does not vanish a.e.\ on any 
	interval $(c,b)$ with $c<b$. Choose $a=c_0<c_1<c_2<\ldots$ such that 
	\[
		\int_{c_n}^b\xi_\phi^TH(t)\xi_\phi\DD t=2^{-n}\int_a^b\xi_\phi^TH(t)\xi_\phi\DD t
		.
	\]
	Set $\psi\DE\phi+\frac\pi 2$ and 
	\[
		\omega_n\DE 2^{-\frac n2}\bigg(\int_{c_{n-1}}^{c_n}\xi_\psi^TH(t)\xi_\psi\DD t\bigg)^{\frac 12},\quad n\geq 1
		,
	\]
	and let $(\omega_n^*)_{n=1}^\infty$ be the nonincreasing rearrangement%
	\/\footnote{Since $\sigma(H)$ is assumed to be discrete, we have $(\omega_n)_{n=1}^\infty\in c_0$. Hence, this sequence
		can be rearranged so to become nonincreasing.}
	of $(\omega_n)_{n=1}^\infty$. 

	Let $\ms g$ be a regularly varying function with index $>1$, and let $\ms f$ be an asymptotic inverse of $\ms g$ 
	(cf.\ \Cref{U123}). Then
	\begin{align*}
		& n_H(r)=\BigO\big(\ms g(r)\big)\quad\Longleftrightarrow\quad \omega_n^*=\BigO\Big(\frac 1{\ms f(n)}\Big),
		\\
		& n_H(r)=\Smallo\big(\ms g(r)\big)\quad\Longleftrightarrow\quad \omega_n^*=\Smallo\Big(\frac 1{\ms f(n)}\Big).
	\end{align*}
\end{Theorem}

\noindent
Recall \Cref{U36} which, apart from demonstrating validity and failure of the independence theorem, already showed how 
sensitively $\rho_H$ depends on the data. An example of a slightly different kind is the following.

\begin{Example}
\label{U42}
	Let $\alpha_1,\alpha_2\geq 0$ with $(\alpha_1,\alpha_2)\neq(0,0)$, set 
	\[
		h(t)\DE\Big(\frac 1{1-t}\Big)^2\Big(1+\log\frac 1{1-t}\Big)^{-\alpha_1}
		\Big(1+\log^+\log\frac 1{1-t}\Big)^{-\alpha_2},\quad t\in(0,1)
		,
	\]
	and consider the Hamiltonian defined on $(0,1)$ as 
	\[
		H(t)\DE
		\begin{pmatrix}
			1 & 0 
			\\
			0 & h(t)
		\end{pmatrix}
		.
	\]
	Then the spectrum $\sigma(H)$ is discrete, and its convergence exponent is 
	\[
		\begin{cases}
			\infty \CAS \alpha_1=0
			\\
			\frac 2{\alpha_1} \CAS \alpha_1\in(0,2)
			\\
			1 \CAS \alpha_1\geq 2
		\end{cases}
		.
	\]
	In the case that $\alpha_1\in(0,2)$, we moreover know that 
	\[
		0<\limsup_{r\to\infty}\frac{n_H(r)}{r^{\frac 2{\alpha_1}}(\log r)^{-\frac{\alpha_2}{\alpha_1}}}<\infty
		.
	\]
\end{Example}

\REMARKS{%
\item The results presented above are taken from the work \cite{romanov.woracek:ideal} of R.~Romanov and H.~Woracek. 
	The proofs of \Cref{U37} and \Cref{U39} produce sequential conditions, and for the case of convergence class this can be
	rewritten to the integral form stated in \Cref{U37}.
\item The assumption that $\xi_\phi^TH(t)\xi_\phi$ does not vanish a.e.\ on any interval $(c,b)$ with $c<b$ is no loss of
	generality. Because, if it were violated, $\sigma(H)$ would be equal to the spectrum of one self-adjoint realisation of
	the cut-off Hamiltonian $H|_{(a,c)}$ which is in limit circle case.
\item A distribution of $\sigma(H)$ that is more dense than the integers can
	occur only from growth of $H$ towards its limit point endpoint $b$. This becomes nicely visible in the condition from 
	\Cref{U37}: for each $c<b$ the integral \cref{U40} is certainly finite when integration runs only from $a$ to $c$. 
\item Rewriting the sequential conditions in \Cref{U39} to a closed form similar as in \Cref{U37} fails because of the need to
	pass to the nonincreasing rearrangement of $(\omega_n)_{n=1}^\infty$. It would be interesting to have a 
	condition that is more directly formulated in terms of $H$.
\item \Cref{U42} is taken from \cite[Example 1.6]{romanov.woracek:ideal}.
\item Membership of resolvents $(A_H-z)^{-1}$ in the Hilbert-Schmidt ideal $\mf S_2$ was characterised earlier in 
	\cite[Theorem~2.4]{kaltenbaeck.woracek:hskansys}. 
	This case is very simple, since one has available the classical criterion
	that an integral operator is Hilbert-Schmidt if and only if its kernel is $L^2$ (and the proof in the mentioned paper
	exploits this). The independence theorem is not needed in this case. 
\item Using the connection between Krein strings and diagonal Hamiltonians (\Cref{U124}), 
	\Cref{U37} gives an alternative approach to the
	results of I.S.~Kac in \cite{kac:1962} and \cite{kac:1986} for orders between $\frac 12$ and $1$. 
	It must be emphasised in this
	place that the form of the conditions in the mentioned papers, as well as their proof, is very different.
	In fact, we do not have a direct computational argument that shows equivalence of the conditions in Kac's 
	papers and the conditions obtained by rewriting \Cref{U37}. 
}

\section[{Trace class and sparse spectrum: the Weyl coefficient approach}]{Trace class and sparse spectrum: 
	the Weyl coefficient approach}
\label{U110}

Let $H\in\bb H_{a,b}$. We again denote
\[
	\IndexN{n_H(r)}\DE\#\big\{\lambda\in\sigma(H)\DS|\lambda|<r\big\}
	.
\]
For the case that $n_H(r)$ grows slower than $r^2$ we can give a formula in terms of $H$ for the Stieltjes transform 
of $\frac{n_H(\sqrt r)}r$. This result is specific for slow growth, since it relies on a product representation of 
functions with real and symmetric zeroes that does not have exponential factors as a general Weierstra{\ss} product does.

In order to formulate the theorem, we need to introduce some notation. Given a Hamiltonian $H\in\bb H_{a,b}$, we write 
\[
	H(t)=\begin{pmatrix} \IndexN{h_1(t)} & \IndexN{h_3(t)}\\ h_3(t) & \IndexN{h_2(t)}\end{pmatrix},\qquad t\in(a,b)\text{ a.e.},
\]
and set
\[
	\IndexN{\Omega_H(s,t)}\DE\int_s^t H(u)\DD u,\quad 
	\IndexN{\omega_{H,j}(s,t)}\DE\int_s^t h_j(u)\DD u,
\]
for $s,t\in[a,b),s\leq t$. If $H$ is in limit circle case, we can extend the definition of $\Omega_H$ and $\omega_{H,j}$ to 
the range $s,t\in[a,b]$, $s\leq t$. 

\begin{Definition}
\label{U54}
	Let $H\in\bb H_{a,b}$, $r_0\geq 0$, and $c_-,c_+>0$. 
	\begin{Enumerate}
	\item We call $\hat t$ a \IndexS{compatible function}{compatible function} for $H,r_0$ with constants $c_-,c_+$, if 
		\begin{align*}
			& \hat t\DF(r_0,\infty)\to(a,b)
			,
			\\
			& \forall r\in(r_0,\infty)\DP\frac{c_-}{r^2}\leq\det\Omega_H\big(a,\hat t(r)\big)\leq\frac{c_+}{r^2}
			.
		\end{align*}
	\item If $\hat t$ is a compatible function for $H,r_0$, we denote 
		\[
			\Gamma(\hat t)\DE\big\{(t,r)\in(a,b)\times(r_0,\infty)\DS \hat t(r)\leq t\big\}
			.
		\]
	\item We call $(\hat t,\hat s)$ a \IndexS{compatible pair}{compatible pair} for $H,r_0$ with constants $c_-,c_+$, if 
		$\hat t$ is a compatible function for $H,r_0$ with constants $c_-,c_+$, and 
		\begin{align*}
			& \hat s\DF\Gamma(\hat t)\to[a,b)
			,
			\\
			& \forall (t,r)\in\Gamma(\hat t)\DP \hat s(t,r)\leq t,\ 
			\frac{c_-}{r^2}\leq\det\Omega_H\big(\hat s(t,r),t\big)\leq\frac{c_+}{r^2}
			.
		\end{align*}
	\end{Enumerate}
\end{Definition}

\noindent
It is easy to see that compatible pairs always exist provided that $(a,b)$ is not indivisible. 
This follows since $\det\Omega_H$ is a continuous and monotone (in either of its two arguments) function. 

We use the following notation to compare functions up to constants. If $f,g$ are functions defined on some 
set $D$ and taking positive real numbers as values, we write 
\begin{align*}
	& f\IndexN{\lesssim}g\ \DI\ \exists C>0\DQ\forall x\in D\DP f(x)\leq Cg(x)
	,
	\\
	& f\IndexN{\asymp}g\ \DI\ f\lesssim g\wedge g\lesssim f
	.
\end{align*}
In the formulation of the next results recall \Cref{U3} and the discussion in the remarks to \Cref{U108}. 

\begin{Theorem}
\label{U55}
\IndexxS{Theorem!formula for $n_H$ (slow growth)}
	Let $H\in\bb H_{a,b}$ and assume that $H$ is definite and that $\sigma(H)$ is discrete.
	For normalisation assume that $\int_a^b h_1(t)\DD t<\infty$.
	Let $r_0\geq 0$, $c_-,c_+>0$, and let $(\hat t,\hat s)$ be a compatible pair for $H,r_0$ with constants $c_-,c_+$. 
	For $t\in(a,b)$ and $r>r_0$ set
	\begin{align*}
		\IndexN{K_H(t;r)} \DE \mathds{1}_{[a,\hat t(r))}(t)
		\frac{\omega_{H,2}(a,t)h_1(t)}{\frac{c_+}{r^2}+\omega_{H,3}(a,t)^2}
		+ \mathds{1}_{[\hat t(r),b)}(t) \frac{h_1(t)}{\omega_{H,1}\bigl(\hat s(t;r),t\bigr)}
		.
	\end{align*}
	Then
	\begin{align}
	\label{U140}
		r^2\int_0^\infty\frac 1{t+r^2}\cdot\frac{n_H(\sqrt t)}t\DD t\asymp\int_a^b K_H(t;r)\DD t
		,\qquad r>r_0,
	\end{align}
	where this relation includes that one side is finite if and only if the other side is. 
	The constants implicit in this relation depend on $c_-,c_+$ but not on $H,r_0,\hat t,\hat s$.
\end{Theorem}

\noindent
The formula \cref{U140} gives a meaningful result about the speed of growth of $n_H(r)$ only if the integrals on either side are 
finite. The integral on the left-hand side of \cref{U140} is finite for some (equivalently, for every) $r>0$ if and only if 
\[
	\int_1^\infty \frac{n_H(\sqrt t)}{t^2}\DD t<\infty,
\]
which, in turn, is equivalent to $\sum\limits_{\lambda\in\sigma(H)}\frac{1}{\lambda^2}<\infty$. 

If we are willing to accept that constants depend on $H,r_0,\hat t,\hat s$, then the first summand of $K_H(t;r)$ can be
neglected. 

\begin{Proposition}
\label{U72}
	In the situation of \Cref{U55} we have 
	\[
		\int_a^{\hat t(r)}K_H(t;r)\DD t\leq 2\log r+\gamma,\quad r>r_0
		,
	\]
	where $\gamma\in\bb R$ depends on $c_-,c_+,H,r_0,\hat t,\hat s$.
\end{Proposition}

\noindent
Using Abelian and Tauberian theorems we can pass from the Stieltjes transform in \cref{U140} to the function $n_H(r)$ itself. 
This gives rise to the following two theorems, which are the counterparts of \Cref{U37,U39}. 
The reader should observe the particular case of the theorems that $\ms g(r)=r^\rho$ for some $\rho\in(0,2)$, cf.\ 
\Cref{U58,U59}, and also consult the examples given in the remarks below.

The result for \IndexS{convergence class}{convergence class} reads as follows.

\begin{Theorem}
\label{U56}
\IndexxS{Theorem!convergence class!sparse spectrum}
	Let $H\in\bb H_{a,b}$ and assume that $H$ is definite and that $\sigma(H)$ is discrete.
	For normalisation assume that $\int_a^b h_1(t)\DD t<\infty$.
	Let $r_0\geq 1$, $c_-,c_+>0$, and let $(\hat t,\hat s)$ be a compatible pair for $H,r_0$ with constants $c_-,c_+$. 

	Let $\ms g$ be regularly varying with $\Ind\ms g\in[0,2]$ such that $\ms g(r)\ll r^2$ 
	and $\ms g|_Y \asymp 1$ on every compact set $Y \subseteq [1,\infty)$. 
	Moreover:
	\begin{Itemize}
	\item If $\Ind\ms g\in(0,2)$, set $\ms g^*\DE\ms g$.
	\item If $\Ind\ms g=0$, assume that we have a regularly varying function $\ms g^*$ 
		such that $\frac{1}{\ms g^*}$ is locally integrable and
		\[
			\int_t^\infty\int_u^\infty\frac{1}{\ms g^*(s)}\frac{\DD s}{s}\frac{\DD u}{u} 
			\asymp \frac{1}{\ms g(t)}.
		\]
	\item If $\Ind\ms g=2$, assume that we have a regularly varying function $\ms g^*$
		such that $\frac{1}{\ms g^*}$ is locally integrable and
		\[
			\int_1^t\frac{s^2}{\ms g^*(s)}\frac{\DD s}{s} \asymp \frac{t^2}{\ms g(t)}.
		\]
	\end{Itemize}
	Then 
	\[
		\sum_{\lambda\in\sigma(H)}\frac{1}{\ms g(|\lambda|)} < \infty
		\quad\Longleftrightarrow\quad
		\int_{r_0}^\infty\frac 1{r\ms g^*(r)}\int_a^b K_H(t;r)\DD t\DD r<\infty.
	\]
\end{Theorem}

\noindent
In the boundary cases $\Ind\ms g\in\{0,2\}$ we do not known whether a function $\ms g^*$ with the required properties 
can always be found (see the remarks below for a class of functions $\ms g$ where $\ms g^*$ can be found). 
However, if $\ms g^*$ exists, then necessarily $\ms g=\Smallo(\ms g^*)$. 

We have the obvious corollary for the case of 
\IndexS{convergence class w.r.t.\ an order}{convergence class!w.r.t.\ order} $\rho$.

\begin{Corollary}
\label{U58}
	Let $\rho\in(0,2)$, and consider the situation of \Cref{U56}. Then 
	\[
		\sum_{\lambda\in\sigma(H)}\frac 1{|\lambda|^\rho}<\infty
		\quad\Longleftrightarrow\quad
		\int_{r_0'}^\infty\frac 1{r^{\rho+1}}\int_a^b K_H(t;r)\DD t\DD r<\infty
		.
	\]
\end{Corollary}

\noindent
Let us emphasize the case that $\rho=1$ where we have a characterisation of 
\IndexS{trace class resolvents}{trace class resolvents}. 
In addition, we can give a neat formula for the trace of the inverse of the operator $A_H$. 

\begin{Proposition}
\label{U68}
	Consider the situation of \Cref{U56}. Then
	\[
		A_H^{-1} \in \mathfrak{S}_1
		\quad\Longleftrightarrow\quad
		\int_{r_0}^\infty\frac 1{r^2}\int_a^b K_H(t;r)\DD t\DD r<\infty.
	\]
	If $A_H^{-1}\in\mf S_1$, then
	\[
		\Tr\bigl(A_H^{-1}\bigr) = -\lim_{t\to b}\int_a^t h_3(s)\DD s.
	\]
\end{Proposition}

\noindent
In the context of finite type or minimal type we obtain the following result. 

\begin{Theorem}
\label{U57}
\IndexxS{Theorem!finite/minimal type!sparse spectrum}
	Let $H\in\bb H_{a,b}$ and assume that $H$ is definite and that $\sigma(H)$ is discrete.
	For normalisation assume that $\int_a^b h_1(t)\DD t<\infty$.
	Let $r_0\geq 1$, $c_-,c_+>0$, and let $(\hat t,\hat s)$ be a compatible pair for $H,r_0$ with constants $c_-,c_+$. 
	\begin{Enumerate}
	\item We have 
		\[
			n_H(r)\lesssim\int_a^b K_H(t;r)\DD t.
		\]
		The constant implicit in this relation depends on $c_-,c_+$ but not on $H,r_0,\hat t,\hat s$.
	\item Let $\ms g$ be regularly varying with $\Ind\ms g\geq 0$ such that 
		$\int_1^\infty\frac{\ms g(t)}{t^3}\DD t<\infty$ and $\ms g|_Y \asymp 1$ on every compact set 
		$Y\subseteq [1,\infty)$. Set 
		\[
			\ms g_*(r)\DE
			\begin{cases}
				\int_1^r\frac{\ms g(t)}t\DD t \CAS \Ind\ms g=0,
				\\[1ex]
				\ms g(r) \CAS \Ind\ms g\in(0,2),
				\\[1ex]
				r^2\int_r^\infty\frac{\ms g(t)}{t^3}\DD t \CAS \Ind\ms g=2.
			\end{cases}
		\]
		Then there exists $C_->0$, which depends on $c_-,c_+,\Ind\ms g$ but not on $H,r_0,\hat t,\hat s,\ms g$, such that
		\[
			\limsup_{r\to\infty}\frac {C_-}{\ms g_*(r)}\int_a^b K_H(t;r)\DD t
			\leq \limsup_{r\to\infty}\frac{n_H(r)}{\ms g(r)}
			.
		\]
	\end{Enumerate}
\end{Theorem}

\noindent
Note that in the boundary cases $\Ind\ms g\in\{0,2\}$ we have $\ms g=\Smallo(\ms g_*)$.

The obvious corollary for the case of \IndexS{type w.r.t.\ an order}{type!w.r.t.\ order} $\rho$ reads as follows.

\begin{Corollary}
\label{U59}
	Let $\rho\in(0,2)$, and consider the situation of \Cref{U57}. Then 
	\[
	\limsup_{r\to\infty}\frac {C_-}{r^\rho}\int_a^b K_H(t;r)\DD t \leq 
		\limsup_{r\to\infty}\frac{n_H(r)}{r^\rho}\leq
		\limsup_{r\to\infty}\frac {C_+}{r^\rho}\int_a^b K_H(t;r)\DD t
	\]
	where $C_-,C_+$ depend on $c_-,c_+,\Ind\ms g$ but not on $H,r_0,\hat t,\hat s,\ms g$.
\end{Corollary}

\REMARKS{%
\item These results are taken from the recent work \cite{langer.reiffenstein.woracek:kacest-arXiv} 
	of M.~Langer, J.~Reiffenstein, and H.~Woracek. 
	The proof of \Cref{U55} is analytic and relies on a trick. It first uses the canonical differential 
	equation which leads to an integral formula for $\log|w_{H,22}(z)|$. 
	The trick is to realise that the integrand in this expression is, for each single point, the imaginary part of the Weyl
	coefficient of some canonical system depending on the point. This allows to invoke the work about high-energy
	asymptotics of Weyl coefficients from \cite{langer.pruckner.woracek:heniest} and \cite{reiffenstein:imq}. 
\item In contrast to the results in \Cref{U101}, the statements in this section are interesting for Hamiltonians in limit circle
	and limit point case equally as much.
\item Assume that $H$ is in the limit circle case. Then we can use the formula
	\begin{equation}
	\label{U77}
		\frac{r^2}{2} \int_0^\infty \frac{1}{t+r^2} \frac{n_H(\sqrt t)}{t} \DD t=\log |w_{H,22}(ir)|
		,
	\end{equation}
	to obtain knowledge about the growth of $\max_{|z|=r}\|W_H(z)\|$ when $r\to\infty$. 
	The formula \cref{U77}, and the fact that the growth of one entry of $W_H$ along the imaginary axis governs the behaviour 
	of the whole monodromy matrix in the whole plane, holds for function theoretic reasons, cf.\ \Cref{U136} and \Cref{U130}. 
\item It is obvious that the formula \cref{U140} is cumbersome, and in general difficult or even impossible to evaluate. 
	If $H$ is in limit circle case there exists a method that considerably simplifies this problem, namely, an 
	algorithm that evaluates the right-hand side of \cref{U140} up to a small error. This result will be presented 
	in \Cref{U132}.
\item If $\det H(t)=0$ for a.a. $t$ (which is necessary for the spectrum to be sparser than the integers, by the Krein-de~Branges
	formula), we have the intuition that the speed of growth of $n_H(r)$ is related to the maximum local 
	rotation of $\Ran H(t)$. This is quantified by the notion of a compatible pair from \Cref{U54}, which measures the size of 
	$\det\Omega_H(s,t)$. 
	The connection is made by the following fact, cf.\ \cite[Lemma~6.3]{langer.reiffenstein.woracek:kacest-arXiv}:
	for a Hamiltonian of the form $H(t)=\xi_{\phi(t)}\xi_{\phi(t)}^T$, so that $\phi(t)$ gives the rotation of $\Ran H(t)$, 
	it holds that 
	\begin{equation}
	\label{U218}
		\det\Omega_H(s,t)=\frac 12\int_s^t\int_s^t\sin^2\big(\phi(x)-\phi(y)\big)\DD x\DD y
		.
	\end{equation}
	A more direct instantiation of this intuition is the algorithm mentioned in the previous item. 
\item Consider the situation where the regularly varying function $\ms g$ in \Cref{U56,U57} is a 
	\IndexS{Lindel\"of comparison function}{Lindel\"of comparison function} \cref{U83}. First, assume that $\Ind\ms g=0$:
	\begin{equation}
	\label{U70}
		\ms g(r) = \prod_{k=n}^N\bigl(\log^{[k]}r\bigr)^{\beta_k}
	\end{equation}
	for large enough $r$ with $n,N\in\bb N$, $n\leq N$ and $\beta_n>0$.
	Then a function $\ms g^*$ as in \Cref{U56} exists, namely, 
	\[
		\ms g^*(r) \DE \ms g(r)\cdot\log r\prod_{j=1}^n\log^{[j]}r.
	\]
	Second, let $\Ind\ms g=2$:
	\begin{equation}
	\label{U71}
		\ms g(r) = r^2\prod_{k=n}^N\bigl(\log^{[k]}r\bigr)^{\beta_k}
	\end{equation}
	for large enough $r$ with $n,N\in\bb N$, $n\leq N$ and $\beta_n<0$. 
	Again a function $\ms g^*$ as in \Cref{U56} exists, namely, 
	\[
		\ms g^*(r) \DE \ms g(r)\cdot\prod_{j=1}^n\log^{[j]}r.
	\]
	Somewhat surprisingly, the gap between $\ms g$ and $\ms g^*$ is different depending on whether
	$\Ind\ms g=0$ or $\Ind\ms g=2$.  

	In \Cref{U57} we have in both cases (that $\ms g$ is of the form \cref{U70} or \cref{U71}) that 
	\[
		\ms g_*(r)=\ms g(r)\cdot\prod_{j=1}^n\log^{[j]}r.
	\]
\item With the same method as in \Cref{U68} we can compute the trace of powers $A_H^{-p}$ for $p\in\{2,3,\ldots\}$. 
	Recall that membership of $A_H^{-1}$ in $\mf S_p$ for such $p$ is characterised in \Cref{U38}.

	For $p=2$ we reobtain the formula for the Hilbert--Schmidt norm given in \cite{kaltenbaeck.woracek:hskansys}:
	if $A_H^{-1}\in\mf S_2$, then 
	\[
		\sum_{\lambda\in\sigma(H)}\frac{1}{\lambda^2} = 2\int_a^b \omega_{H,2}(a,t)h_1(t) \DD t.
	\]
	For $p=3$ we get: if $A_H^{-1}\in\mf S_3$, then 
	\[
		\sum_{\lambda\in\sigma(H)}\frac{1}{\lambda^3}
		= 12\cdot\lim_{t\to b}\int_a^t\int_a^s\omega_{H,2}(a,x)h_3(x)\DD x\,h_1(s)\DD s.
	\]
	The formulae for higher $p$ are getting more cumbersome, but retain their structure of being iterated integrals over
	entries of $H$. This is expected, thinking of the work \cite{kac:1962} where for a Krein string membership in
	$\mf S_p$ for integers $p\geq 2$ is characterised by finiteness of certain iterated integrals.
\item For orders between $1$ and $2$ we have an overlap with the operator theoretic method from \Cref{U101}. We do not know if 
	one can show equivalence of the respective conditions by direct computation. 
\item For Krein strings, a criterion similar to \Cref{U56} is given by I.S.~Kac \cite{kac:1986}.
	Equivalence of the conditions can, at least under certain smoothness assumptions, be verified by direct computation. 
}

%**************************************************************************
%***                            Last Change: Mon 31 Mar 2025 00:14
%***   < PART II >
%***
%**************************************************************************

\clearpage
\PART{II}{The limit circle case}
\label{U146}

%%%%%%%%%%%%%%%%%%%%%%%%%%%%%%%%%%%%%%%%%%%%%%%%%%%%%%%%%%%%%%%%%%%%%%%%%%%%%%%%%%%%%%%%%

\Intro{%
	Assume that the Hamiltonian $H$ is in limit circle case. Then the fundamental solution exists up to the right endpoint of the 
	interval, and we have the monodromy matrix $W_H(z)$ which is an entire function of exponential type. 
	The spectrum $\sigma(H)$ coincides with the zero set of its entry $w_{H,22}(z)$, in particular $\sigma(H)$ is discrete. 

	These facts add another dimension to the question for density of $\sigma(H)$: we can use the connection between growth
	and zero distribution of entire functions. 
	For example, the convergence exponent of $\sigma(H)$ is equal to the order of the entire function $w_{H,22}$.
	More detailed information can be extracted from the central connection \cref{U77} between $n_H(r)$ and $w_{H,22}(z)$.

	The Krein-de~Branges formula gives a simple expression in terms of $H$ for the exponential type of 
	$w_{H,22}(z)$ and hence for $\lim_{r\to\infty}\frac{n_H(r)}r$, cf.\ \Cref{U53}. It fully settles the case when 
	the type is positive, equivalently, when $\det H(t)$ does not vanish identically. In this case the eigenvalues 
	have integer asymptotics.

	The case that $\det H(t)=0$ is much more involved. One may say that the growth of $W_H(z)$ is
	proportional to the maximum local rotation of $\Ran H(t)$. This vague statement is instanciated in different ways.
	An algorithm that evaluates $\log |w_{H,22}(ir)|$ up to a possible error of size $\log r$ is based on 
	partitioning the interval into parts with equal rotation, cf.\ \Cref{U144}. An upper bound, for
	$\log\max_{|z|=r}\|W_H(z)\|$, is obtained by approximating $H$ with piecewise constant Hamiltonians, cf.\ \Cref{U43}.
	Another result determines the order of $W_H$ using coverings of $(a,b)$ that make certain sums involving integrals of $H$ small, cf.\ \Cref{U76}. 
}
\begin{center}
	{\large\bf Table of contents}
\end{center}
\begin{flushleft}
	\S.\,\ref{U109}\ The Krein--de~Branges formula\ \dotfill\quad\pageref{U109}
	\\[1mm]
	\S.\,\ref{U132}\ Algorithm to evaluate growth\ \dotfill\quad\pageref{U132}
	\\[1mm]
	\S.\,\ref{U111}\ Romanov's Theorem I: bound by discretisation\ \dotfill\quad\pageref{U111}
	\\[1mm]
	\S.\,\ref{U112}\ Romanov's Theorem II: bound by coverings\ \dotfill\quad\pageref{U112}
	\\[1mm]
\end{flushleft}
\makeatother
\renewcommand{\thesection}{\arabic{section}}
\renewcommand{\thelemma}{\arabic{section}.\arabic{lemma}}
\makeatletter
\clearpage

%%%%%%%%%%%%%%%%%%%%%%%%%%%%%%%%%%%%%%%%%%%%%%%%%%%%%%%%%%%%%%%%%%%%%%%%%%%%%%%%%%%%%%%%%

%
%
%
\section[{The Krein--de~Branges formula}]{The Krein--de~Branges formula}
\label{U109}

The \IndexS{Krein--de~Branges formula}{Krein--de~Branges formula} gives precise information on the distribution of 
$\sigma(H)$ compared to the integers. 

\begin{Theorem}
\label{U53}
\IndexxS{Theorem!Krein--de~Branges formula!limit circle}
	Let $H\in\bb H_{a,b}$ be in limit circle case. Then the entries of $W_H(z)$ are entire functions of Cartwright class in $\bb C^+$ and $\bb C^-$. Denote by 
	\[
		0<\lambda^+_1<\lambda^+_2<\ldots\quad\text{and}\quad 0>\lambda^-_1>\lambda^-_2>\ldots
	\]
	the (finite or infinite) sequences of positive and negative, respectively, elements of $\sigma(H)$.
	Then%
	\/\footnote{Here the limit of a finite sequence is understood as being equal to $0$.}
	\[
		\pi\cdot\lim_{n\to\infty}\frac n{\lambda^+_n}=\pi\cdot\lim_{n\to\infty}\frac n{|\lambda^-_n|}
		=\limsup_{y\to\infty}\frac 1y\log^+|w_{H,ij}(iy)|=\int\limits_a^b\sqrt{\det H(t)}\DD t
		.
	\]
\end{Theorem}

\noindent
There are two essential assertions in this theorem: one is the function theoretic part that we have Cartwright class, and the 
other is the computation of the type as integral over $\sqrt{\det H}$. The first is related to a functional analytic property
of de~Branges' Hilbert spaces of entire functions, the latter stems from a differential inequality. 
The stated asymptotic of eigenvalues is then a consequence (recall \Cref{U198} and \Cref{U196}).

The Krein-de~Branges formula continues to hold for a certain class of limit point Hamiltonians. 
The corresponding variant reads as follows. 

\begin{Theorem}
\label{U1}
\IndexxS{Theorem!Krein--de~Branges formula!Pontryagin type}
	Let $H\in\bb H_{a,b}$ be of Pontryagin type (see \Cref{U118}), and let 
	\[
		0<\lambda^+_1<\lambda^+_2<\ldots\quad\text{and}\quad 0>\lambda^-_1>\lambda^-_2>\ldots
	\]
	be the (finite or infinite) sequences of positive and negative, respectively, elements of $\sigma(H)$. 
	Then the product 
	\[
		A(z)\DE\lim_{R\to\infty}\prod_{\substack{\lambda\in\sigma(H)\\ 0<|\lambda|\leq R}}\Big(1-\frac z\lambda\Big)
	\]
	converges locally uniformly on $\bb C$, and is an entire function of Cartwright class in $\bb C^+$ and $\bb C^-$. We have 
	\begin{equation}
	\label{U21}
		\sqrt{\det H}\in L^1(a,b)
		,
	\end{equation}
	and 
	\[
		\pi\cdot\lim_{n\to\infty}\frac n{\lambda^+_n}=\pi\cdot\lim_{n\to\infty}\frac n{|\lambda^-_n|}=
		\limsup_{y\to\infty}\frac 1y\log^+|A(iy)|=\int\limits_a^b\sqrt{\det H(t)}\DD t
		.
	\]
\end{Theorem}

\REMARKS{%
\item Original references for \Cref{U53} are \cite[(3.4)]{krein:1951} and 
	\cite[Theorem~X]{debranges:1961}. The version stated in \Cref{U1} is obtained by combining 
	\cite[Theorem~4.1]{langer.woracek:expty} with \cite[Theorem~4.21]{langer.woracek:gpinf}. 
	A variant for higher-dimensional canonical systems (in limit circle case) is also available, see e.g.\ 
	\cite[Theorem~6.1]{gohberg.krein:1970}.
\item The Krein-de~Branges formula sheds light on the specialty of integer distribution from a function theoretic perspective. 
	Namely, that one has the Cartwright class property of the monodromy matrix and thus obtains at most integer density in the
	limit circle case. 
	Compare this with the origin of the distinction into dense or sparse spectra in \Cref{U101,U110}: the operator
	theoretic facts that Schatten--von~Neumann classes $\mf S_p$ are s.n.-ideals only if $p\geq 1$ and that for $p=1$ the
	Matsaev property fails.
\item Using the connection \Cref{U284} between power moment problems and Hamburger Hamiltonians, 
	the Krein--de~Branges formula gives a one-line proof of the classical result of M.~Riesz from \cite{riesz:1923a} 
	(see also \cite[Theorem~2.4.3]{akhiezer:1965}), that the entries of the Nevanlinna matrix of an indeterminate moment 
	problem are of minimal exponential type. Namely: the determinant of a Hamburger Hamiltonian is identically equal to $0$. 
	The version in \Cref{U1} will lead in the obvious way to a generalisation of Riesz's result; however, this has not been
	carried out yet.
\item Using the connection between strings and diagonal Hamiltonians (laid out in \Cref{U124}), 
	the Krein--de~Branges formula yields a proof
	of M.G.~Krein's result announced in \cite{krein:1951a} (see also \cite[Theorema~1;4]{krein:1952a} or 
	\cite[11.8$^\circ$]{kac.krein:1968}) about the asymptotics of the eigenvalues of a string. 
	The version in \Cref{U1} leads to a generalisation of Krein's result, cf.\ \cite{woracek:asymp}.
}

\section[{Algorithm to evaluate growth}]{Algorithm to evaluate growth}
\label{U132}

The formula given in \Cref{U55} is explicit in $H$, but it is often difficult to evaluate the integral 
on the right-hand side of \cref{U140} in practice. For Hamiltonians in limit circle case there is an algorithm that simplifies evaluation of the formula at the cost of a small loss of precision.

Recall the notation $\IndexN{\Omega_H(s,t)}\DE\int_s^t H(u)\DD u$  and \cref{U218} 
where we explained that $\det\Omega_H(s,t)$ quantifies the maximum local rotation of $\Ran H(t)$.

\begin{Definition}
\label{U143}
	Let $H\in\bb H_{a,b}$ be in limit circle case. For each $r>0$ we define points $\IndexN{\sigma_j^{(r)}}$ and a
	number $\IndexN{\kappa_H(r)}$ by the following procedure.
	\begin{Itemize}
	\item Set $\sigma_0^{(r)}\DE a$.
	\item If $\det\Omega_H\big(\sigma_{j-1}^{(r)},b\big)>\frac 1{r^2}$, let $\sigma_j^{(r)}\in(\sigma_{j-1}^{(r)},b)$ be the 
		unique point such that 
		\[
			\det\Omega_H\big(\sigma_{j-1}^{(r)},\sigma_j^{(r)}\big)=\frac 1{r^2}
			.
		\]
		Otherwise, set $\sigma_j^{(r)}\DE b$ and $\kappa_H(r)\DE j$, and terminate.
	\end{Itemize}
\end{Definition}

\noindent
By an argument involving Minkowski's determinant inequality this algorithm terminates for each $r$ after finitely many steps. 
In fact, we have
\[
	\kappa_H(r) \leq \Big\lfloor r\cdot\sqrt{\det\Omega_H(a,b)}\Big\rfloor+1.
\]
The function $\kappa_H$ is nondecreasing.

\begin{Theorem}
\label{U144}
\IndexxS{Theorem!algorithmic bounds}
	Let $H\in\bb H_{a,b}$ and assume that $H$ is definite and in limit circle case. 
	Then there exists $r_0>0$ such that for $r>r_0$
	\begin{equation}
	\label{U78}
		\kappa_H(r)\lesssim \log\Big(\max_{|z|=r}\|W_H(z)\|\Big) \lesssim \kappa_H \Big(\frac{r}{\log r} \Big)\log r
		.
	\end{equation}
	The number $r_0$ and the constants implicit in \cref{U78} depend on $\Tr\Omega_H(a,b)$. 
\end{Theorem}

\noindent
When comparing $n_H(r)$ to a regularly varying function, the algorithmic character of this theorem becomes more pronounced: in
this case it is enough to compute $\kappa_H(r)$ for countably many values of $r$. 

\begin{Proposition}
\label{U81}
	Let $H\in\bb H_{a,b}$ and assume that $H$ is definite and in limit circle case. 
	Moreover, let $\ms g$ be regularly varying and nondecreasing, and let $(r_n)_{n\in\bb N}$ be an increasing and unbounded 
	sequence of positive numbers with $\sup_{n\in\bb N}\frac{r_{n+1}}{r_n}<\infty$. Then 
	\begin{align*}
		\limsup_{n\to\infty}\frac{\kappa_H\big(\frac{r_n}{\log r_n} \big)\cdot\log r_n}{\ms g(r_n)}<\infty
		& \ \Longrightarrow\ 
		\limsup_{r\to\infty}\frac{\log\big(\max_{|z|=r}\|W_H(z)\|\big)}{\ms g(r)}<\infty
		\\
		& \ \Longrightarrow\ 
		\limsup_{n\to\infty}\frac{\kappa_H(r_n)}{\ms g(r_n)}<\infty
		.
	\end{align*}
\end{Proposition}

\noindent
In particular, we can compute the convergence exponent of $\sigma(H)$ from the sequence $(\kappa_H(r_n))_{n\in\bb N}$. 

\begin{Corollary}
\label{U145}
	Let $H\in\bb H_{a,b}$ and assume that $H$ is definite and in limit circle case. 
	Moreover, let $(r_n)_{n\in\bb N}$ be an increasing and unbounded sequence of positive numbers with 
	$\sup_{n\in\bb N}\frac{r_{n+1}}{r_n}<\infty$. Then the convergence exponent $\rho_H$ of $\sigma(H)$ is given as
	\begin{align*}
		\rho_H=\limsup_{n\to\infty}\frac{\log\kappa_H(r_n)}{\log r_n}
		.
	\end{align*}
\end{Corollary}

\REMARKS{
\item \Cref{U144} was first proved in the slightly weaker form
	\begin{align}
	\label{U220}
		\kappa_H(r)\lesssim \log |w_{H,22}(ir)| \lesssim \kappa_H(r)\log r
	\end{align}
	in a paper of M.~Langer, J.~Reiffenstein, and H.~Woracek, cf.\ 
	\cite[Theorem~5.3]{langer.reiffenstein.woracek:kacest-arXiv}. Its proof is obtained by direct but somewhat
	tricky estimates of the integrand $K_H(t;r)$ in \cref{U140}. 
	The improved upper estimate is proved in \cite{reiffenstein:romanii} using Romanov's Theorem I 
	(cf.\ \Cref{U111}) which is a more directly derived general upper bound.
\item If $\kappa_H$ is well-behaved, e.g.,
	$\kappa_H(r) \asymp r^{\rho_H} (\log r)^\sigma$, the gap between the lower and upper bounds in \Cref{U144} is of size
	$(\log r)^{1-\rho_H}$, i.e., it gets smaller as $\rho_H$ increases and closes at $\rho_H=1$.
\item The results presented in this section are bound to the limit circle case, and we do not know an analogue of \Cref{U144}
	for the limit point case. 
\item Regarding sharpness of the estimates in \Cref{U144}, we distinguish cases depending on $\rho_H$.
	\begin{Itemize}
	\item $\rho_H=0$: In all examples we are aware of, the upper bound in \Cref{U144} is the correct one, i.e.,
		\begin{equation}
		\label{U221}
				\log\Big(\max_{|z|=r}\|W_H(z)\|\Big) \asymp \kappa_H \Big(\frac{r}{\log r} \Big)\log r
				.
		\end{equation}
		Trivial examples are Hamiltonians consisting of a finite number of indivisible intervals only, for which both
		sides of \cref{U221} are $\asymp$ to $\log r$. A concrete and nontrivial example is
		\[
			H(t)
			\DE
			\begin{cases}
				\smmatrix 1000 \CAS e^{-(2n+1)} \le t < e^{-2n},
				\\[1.5ex]
				\smmatrix 0001 \CAS e^{-(2n+2)} \le t < e^{-(2n+1)}
			\end{cases}
			\qquad\text{for }n\in\bb N_0.
		\]
		For this Hamiltonian, both sides of \cref{U221} are $\asymp$ to $(\log r)^2$.
	\item $\rho_H \in (0,1)$: Here we only have examples for which the lower bound in \Cref{U144} is correct:
		\begin{equation}
		\label{U222}
			\log\Big(\max_{|z|=r}\|W_H(z)\|\Big) \asymp \kappa_H (r)
			.
		\end{equation}
		A concrete and nontrivial example is 
		$H(t)\DE\xi_{\phi(t)}\xi_{\phi(t)}^T$ where $\phi(t)$ is the Weierstra{\ss}\ function 
		$\phi(t)\DE\sum_{n=0}^\infty\frac 1{2^n}\cos(8^n\pi t)$. For this Hamiltonian 
		\[
			\log\Big(\max_{|z|=r}\|W_H(z)\|\Big) \asymp \kappa_H(r)\asymp  r^{\frac 34}.
		\]
		Varying the parameters of the Weierstra{\ss}\ function yields a scale of examples covering the range of orders 
		$(\frac 12,1)$, cf.\ \Cref{U184}. Examples with $\rho_H \leq \frac 12$ can be found in \Cref{U134} or, more
		concretely, in \Cref{U149}.
	\item $\rho_H=1$: Typically there is no gap in this case, such as for $\kappa_H(r) \asymp r (\log r)^\sigma$. Trivial
		examples are Hamiltonians with $\det H \not\equiv 0$, for which both sides of \cref{U222} are $\asymp$ to $r$.
		We do not have examples where $\kappa_H$ is of a more complicated form.
	\end{Itemize}
\item We do not have examples with 
	\[
		\kappa_H (r)\ll\log\big(\max_{|z|=r}\|W_H(z)\|\big)\ll\kappa_H \big(\frac r{\log r}\big)\log r
		.
	\]
}

\section[{Romanov's Theorem I: bound by discretisation}]{Romanov's Theorem I: bound by discretisation}
\label{U111}

If $H$ consists of a finite number of indivisible intervals only, then $W_H(z)$ is a polynomial whose degree equals the number
of indivisible intervals. This suggests that $W_H(z)$ will grow slowly if $H$ can be approximated well by such simple
Hamiltonians. The below theorem quantifies this intuition.

\begin{Theorem}
\label{U43}
\IndexxS{Theorem!Romanov's Theorem I (improved version)}
	Let $H\in\bb H_{a,b}$ be in limit circle case, assume that $\det H(t)=0$ a.e., 
	and write $H(t)=\Tr H(t)\cdot\xi_{\phi(t)}\xi_{\phi(t)}^T$ with some measurable function $\phi\DF(a,b)\to\bb R$. 
	Assume we are given a set of parameters $\mf p$ that consists of a number $N\in\bb N$, a partition
	$a=y_0<y_1<\cdots<y_N=b$, rotation parameters $\psi_1,\ldots,\psi_N\in\bb R$, and distortion parameters
	$a_1,\ldots,a_N\in(0,1]$, and set
	\begin{align*}
		A_1(\mf p) \DE &\, 
		\sum_{j=1}^N a_j^2\int_{y_{j-1}}^{y_j}\cos^2\big(\phi(t)-\psi_j\big)\cdot\Tr H(t)\DD t
		,
		\\
		A_2(\mf p) \DE &\, 
		\sum_{j=1}^N \frac 1{a_j^2}\int_{y_{j-1}}^{y_j}\sin^2\big(\phi(t)-\psi_j\big)\cdot\Tr H(t)\DD t
		,
		\\
		A_3(\mf p) \DE &\, 
		\sum_{j=1}^{N-1} \log\bigg(
		\max\Big\{\frac{a_j}{a_{j+1}},\frac{a_{j+1}}{a_j}\Big\}\cdot\big|\cos\big(\psi_j-\psi_{j+1}\big)\big|
		\\
		&\mkern300mu +\frac{|\sin(\psi_j-\psi_{j+1})|}{a_ja_{j+1}}
		\bigg)
		,
		\\
		A_4(\mf p) \DE &\, 
		-\log a_1-\log a_N
		.
	\end{align*}
	Then, for every $r>0$,  
	\begin{equation}
	\label{U44}
		\log\Big(\max_{|z|=r}\|W_H(z)\|\Big)\leq r\cdot\big(A_1(\mf p)+A_2(\mf p)\big)+A_3(\mf p)+A_4(\mf p)
		.
	\end{equation}
%	where $\|\Dummy\|$ denotes the spectral norm on $\bb C^{2\times 2}$.
\end{Theorem}

\begin{Remark}
\label{U209}
	If a parameter set $\mf p$ is fixed, then \cref{U44} does not give any new information about the growth when
	$r\to\infty$. On the contrary, we know even better since $\det H=0$ implies $\log\|W_H(z)\|=\Smallo(r)$.

	The strength of \Cref{U43} lies in the quantitative aspect: the estimate \cref{U44} holds for all $y_j,\psi_j,a_j,r$. 
	The trick to successfully apply the theorem is to make the approximation $\mf p$ dependent on $r$: 
	\begin{Itemize}
	\item Choose a family of parameters $(\mf p(r))_{r>0}$ and apply \Cref{U43} with each of them,
	\item for each $r$ use the bound \cref{U44} obtained from $\mf p(r)$, 
	\item send $r\to\infty$. 
	\end{Itemize}
\end{Remark}

\noindent
As we see from this remark, applying \Cref{U43} usually requires some clever guessing. 
The next result only requires guessing the partitions, and 
then uses a predefined $r$-dependent choice of rotation and distortion parameters. It yields an upper estimate that is almost as good as the one from \Cref{U43} (see the remarks for details).

\begin{Theorem}
\label{U73}
\IndexxS{Theorem!Romanov's Theorem I (partition only)}
	Let $H\in\bb H_{a,b}$ be in limit circle case with $\det H=0$. Assume we have $k\in\bb N$ and a partition 
	\begin{equation}
	\label{U74}
		a=y_0<y_1<\cdots<y_k=b.
	\end{equation}
	Set
	\begin{equation}
	\label{U75}
		M \DE \frac 1k\sum_{j=1}^k\sqrt{\det\Omega(y_{j-1},y_j)}
		.
	\end{equation}
	Then, for every $r>0$, 
	\[
		\log\Big(\max_{|z|=r}\|W_H(z)\|\Big)\leq 
		k\log^+r+3Mr+\gamma\cdot k+\log^+r+\gamma'
		,
	\]
	where $\gamma,\gamma'$ depend only on $\Tr\Omega(a,b)$.
\end{Theorem}

\REMARKS{%
\item A first variant of \Cref{U43} was shown by R.~Romanov in \cite{romanov:2017} as Theorem~1 of that paper. 
	The presently stated version is taken from \cite[Theorem~4.1]{pruckner.woracek:sinqA}. 
	It improves upon Romanov's theorem, even on the rough scale of
	order (for more details see the remarks to \Cref{U106}).
	The proof of \Cref{U43} is carried out by slightly tedious but elementary estimates using Gr\"onwall's Lemma. 
\item \Cref{U73} is taken from work of J.~Reiffenstein \cite{reiffenstein:romanii}. There also explicit values for the constants 
	$\gamma,\gamma'$ are given, namely 
	\[
		\gamma=\log^+\Tr\Omega(a,b)+3+2\log 2,\quad \gamma'=\log^+\Tr\Omega(a,b) +\log 2
		.
	\]
\item The upper bound $\kappa_H (r/\log r )\log r$ in \Cref{U144} is what we obtain when applying \Cref{U73}
	with the $r$-dependent partitions defined in \Cref{U143}. This even yields the more explicit estimate
	\begin{align}
\label{U80}
\log\Big(\max_{|z|=r}\|W_H(b;z)\|\Big)\leq 
		\kappa_H \Big(\frac{r}{\log r} \Big) \big(4\log r+\gamma \big)+\log r+\gamma',
\end{align}
where $\gamma,\gamma'$ are as in the previous item.
	 We point out that in the background lies an application of 
	\Cref{U43} with parameter sets $\mf p(r)$ chosen in a certain deterministic way. 
\item It is unclear whether or not the choice of parameters mentioned in the previous item is optimal, or if making a different
	choice can lead to an asymptotically better upper bound. That is, if $J(r)$ is the infimum of upper bounds produced by
	\Cref{U43}, taken over all parameter sets $\mf p$, we do not know if for all Hamiltonians 
	\[
		\log\big(\max_{|z|=r}\|W_H(z)\|\big) \asymp J(r)
		.
	\]
	However, we have 
\begin{align*}
\log\big(\max_{|z|=r}\|W_H(z)\|\big) \leq J(r)\lesssim\kappa_H \big(\frac{r}{\log r}\big)\log r\lesssim\log r\cdot\log\big(\max_{|z|=r}\|W_H(z)\|\big),
\end{align*}	
where the last inequality follows from \Cref{U144}. Since \eqref{U80} is obtained from \Cref{U73} we see that the optimal upper bounds from \Cref{U43,U73} are both accurate up to a factor of $\log r$.

	Revisiting the context of \Cref{U144}, we do not know if there exists a Hamiltonian $H$ for which
	\[
		\liminf_{r \to \infty} \frac{J(r)}{\kappa_H \big(\frac{r}{\log r} \big)\log r}=0.
	\]
}

\section[{Romanov's Theorem II: bound by coverings}]{Romanov's Theorem II: bound by coverings}
\label{U112}

Minkowski's determinant inequality states that, for any two positive semidefinite matrices $A,B \in \bb C^{2 \times 2}$ the determinant of $A+B$ can be estimated as
\begin{align*}
\sqrt{\det (A+B)} \geq \sqrt{\det A}+\sqrt{\det B},
\end{align*}
where equality holds if and only if $A$ and $B$ are linearly dependent.
For a Hamiltonian $H\in\bb H_{a,b}$ in limit circle case, applying this inequality to matrices of the form $\IndexN{\Omega_H(s,t)}\DE\int_s^t H(x)\DD x$, we get 
\begin{align}
\label{U227}
\sqrt{\det \Omega_H(s,u)} \geq \sqrt{\det \Omega_H(s,t)}+\sqrt{\det \Omega_H(t,u)}, \quad t \in (s,u).
\end{align}
Thinking of $\det \Omega_H$ as a quantifier for the rotation of $\Ran H(x)$ (see \eqref{U218}), the matrices $\Omega_H(s,t)$ and $\Omega_H(t,u)$ should be close to linearly dependent when $\Ran H(x)$ rotates slowly in $[s,u]$, leading to a small relative loss when applying \eqref{U227}. 

The main result of this section determines the order of the monodromy matrix based on how well $\sqrt{\det \Omega}$ can be split up using \eqref{U227}, thus establishing a connection between the rotation of $\Ran H(x)$ and the growth of the monodromy matrix. Recall the notation $\rho_H$ for the common order of all entries of the monodromy matrix (cf.\ \Cref{U201}).

% Due to the Krein-de~Branges formula (cf. \Cref{U53}) we have $\rho_H=1$, provided that $\sqrt{\det H}$ does not vanish identically. Otherwise the Krein-de~Branges formula only yields $\rho_H \in [0,1]$.

%This section is dedicated to a theorem that determines $\rho_H$ in terms of the function $\sqrt{\det \Omega_H}$, where $\IndexN{\Omega_H(s,t)}\DE\int_s^t H(u)\DD u$. That is, instead of the integral of $\sqrt{\det H}$, which is the main input of the Krein-de~Branges formula, we consider $\sqrt{\det (.)}$ of the integral of $H$. 
%\[
%\Omega_H (s,t)=\int_s^t H(x) \, dx
%\]
% (recall the definition of $\Omega_H$ in \Cref{U110}), and that is reminiscent of the Krein-de~Branges formula. As in \eqref{U218} we think of $\det \Omega_H$ as a quantifier for the  rotation of $\Ran H(t)$.

\begin{Theorem}
\label{U223}
\IndexxS{Theorem!Romanov's Theorem II (general version)}
	Let $H\in\bb H_{a,b}$ be in limit circle case. 
	Then the order $\rho_H$ of the monodromy matrix is equal to the infimum of all numbers $\alpha>0$ that satisfy
	\begin{itemize}
	\item[$(\ast )$]
		For every sufficiently large $N\in\bb N$ there exists a covering $\{[c_j,d_j] \DS j=1,\ldots,k(N)\}$ of $[a,b]$,
		such that
		\[
			k(N)\leq N,\qquad 
			\sum_{j=1}^{k(N)}\sqrt{\det \Omega_H (c_j,d_j)}\lesssim N^{1-\frac 1\alpha}
			.
		\]
	\end{itemize}
\end{Theorem}

\noindent
For every Hamiltonian there is an explicit family of coverings testifying to the validity of $(\ast)$ for every $\alpha>\rho_H$. The coverings in this family are partitions of $[a,b]$ that are ``equidistant'' with respect to $\det \Omega_H$.

\begin{lemma}
\label{U76}
Let $H\in\bb H_{a,b}$ be in limit circle case, but not a finite rank Hamiltonian. Then for any $\alpha > \rho_H$ the property $(\ast)$ from \Cref{U223} is satisfied, where the coverings can be chosen as $\{[\sigma_{j-1}^{(r(N))},\sigma_{j}^{(r(N))}] \DS j=1,\ldots,N\}$, with $\IndexN{\sigma_j^{(r)}}$ defined as in \Cref{U143} and $r(N)$ chosen such that $\kappa_H(r(N))=N$.
\end{lemma}

\REMARKS{%
\item These results are taken from work of J.~Reiffenstein \cite{reiffenstein:romanii}. 
\item \Cref{U76}, which leads to the lower bound in \Cref{U223}, is shown using the lower bound in \Cref{U144}. The upper bound in \Cref{U223} is a consequence of \Cref{U73}, which in turn uses Romanov's Theorem I, cf. \Cref{U111}.
\item \Cref{U223} is a striking extension of a theorem of R.~Romanov, namely, of \cite[Theorem~2]{romanov:2017}. 
	In Romanov's original version $H$ is supposed to be a.e.\ diagonal with $\det H \equiv 0$ and $\Tr H \equiv 1$. 
	This means that $H(t)$ can take only the two values
	\[
		H_1 \DE \begin{pmatrix}
		1 & 0 \\
		0& 0
		\end{pmatrix}
		, \qquad
		H_2 \DE \begin{pmatrix}
		0 & 0 \\
		0& 1
		\end{pmatrix}.
	\]
\item Since the coverings provided in \Cref{U76} come from partitions of $[a,b]$, a reformulation of \Cref{U223} using partitions instead of coverings is possible. One reason to stick with coverings is historical, as that is what was used in \cite[Theorem~2]{romanov:2017}.
\item No criterion is known for deciding whether or not $(\ast)$ is satisfied for $\rho_H$ itself.
\item An interesting application of (Romanov's original variant of) \Cref{U223} can be found in \cite[Section~7.2]{romanov:2017}.
}

%**************************************************************************
%***                            Last Change: Mon 31 Mar 2025 11:28
%***   < PART III >
%***
%**************************************************************************

\clearpage
\PART{III}{Moment problems and Jacobi matrices}
\label{U127}

%%%%%%%%%%%%%%%%%%%%%%%%%%%%%%%%%%%%%%%%%%%%%%%%%%%%%%%%%%%%%%%%%%%%%%%%%%%%%%%%%%%%%%%%%

\Intro{%
	In this part we study the context of the Hamburger power moment problem. 
	We focus on the indeterminate case, and mainly present results formulated in the language of Jacobi
	parameters. 

	First we discuss a somewhat isolated theorem that is mainly of theoretical interest, cf.\ \Cref{U28}. 
	In theory it fully determines the growth of the Nevanlinna matrix in terms of the coefficients of the orthogonal 
	polynomials. But it has the drawback of being virtually impossible to apply; we know of only one nontrivial application 
	(for this see \Cref{U105}). 

	An elementary lower bound for the growth of the Nevanlinna matrix in terms of the Jacobi parameters is given in
	\Cref{U24}. The historically probably first result where growth different from exponential type is considered is due to 
	M.S.~Liv\v sic dating back to 1939, cf.\ \Cref{U15}. It gives a lower bound in terms of the moment sequence and appears 
	from the nowadays viewpoint as a consequence of \Cref{U24}.

	Recently, another lower bound in terms of Jacobi parameters was established, cf.\ \Cref{U97}. 
	This bound is most easy to apply, and the fact that it holds in full generality is quite surprising. 
	It is closely related to a classical result of Yu.M.~Berezanskii in 1956. We give a recent
	significant improvement of Berezanskii's result in \Cref{U12}.
	Finally, we discuss Jacobi matrices whose parameters have power asymptotics. In this situation a fairly complete picture
	can be given, cf.\ \Cref{U48} and \Cref{U49}.

	Many of the results presented in this part are proven by translating to a canonical system and referring
	to results formulated in the language of Hamiltonian parameters.

	It must be said that there is a large number of results about the determinate case dealing with eigenvalue asymptotics.
	These are out of the scope of this paper, and we do not touch upon that part of the literature.
	}
\begin{center}
	{\large\bf Table of contents}
\end{center}
\begin{flushleft}
	\S.\,\ref{U103}\ Growth in terms of orthogonal polynomials\ \dotfill\quad\pageref{U103}
	\\[1mm]
	\S.\,\ref{U102}\ An elementary lower bound\ \dotfill\quad\pageref{U102}
	\\[1mm]
	\S.\,\ref{U104}\ The Liv\v sic estimate\ \dotfill\quad\pageref{U104}
	\\[1mm]
	\S.\,\ref{U152}\ Convergence exponents as lower bounds\ \dotfill\quad\pageref{U152}
	\\[1mm]
	\S.\,\ref{U107}\ A theorem of Berezanskii\ \dotfill\quad\pageref{U107}
	\\[1mm]
	\S.\,\ref{U120}\ Growth from power asymptotics\ \dotfill\quad\pageref{U120}
	\\[1mm]
\end{flushleft}
\makeatother
\renewcommand{\thesection}{\arabic{section}}
\renewcommand{\thelemma}{\arabic{section}.\arabic{lemma}}
\makeatletter
\clearpage

%%%%%%%%%%%%%%%%%%%%%%%%%%%%%%%%%%%%%%%%%%%%%%%%%%%%%%%%%%%%%%%%%%%%%%%%%%%%%%%%%%%%%%%%%

%
%
%
\section[{Growth in terms of orthogonal polynomials}]{Growth in terms of orthogonal polynomials}
\label{U103}

For an indeterminate moment sequence, the orthonormal polynomials $p_n$ form an orthonormal basis in a reproducing kernel
Hilbert space generated from the Nevanlinna matrix. The growth of the Nevanlinna matrix $W(z)$ can be expressed in terms
of $p_n$: set 
\[
	\IndexN{\Delta(z)}=\Big(\sum_{n=0}^\infty|p_n(z)|^2\Big)^{\frac 12}
	,
\]
then 
\[
	\log\Big(\max_{|z|=r}\|W(z)\|\Big)=\log\Big(\max_{|z|=r}\Delta(z)\Big)+\BigO(\log r)
	.
\]
This is seen by an elementary argument using the reproducing kernel Hilbert space built from $W(z)$. 

\begin{Theorem}
\label{U28}
\IndexxS{Theorem!Nevanlinna matrix!Berg-Szwarc}
	Let $(s_n)_{n=0}^\infty$ be an indeterminate moment sequence, and write the orthonormal polynomial $p_n$ of degree $n$
	as 
	\begin{equation}
	\label{U32}
		p_n(z)=\sum_{k=0}^n b_{k,n}z^k
		.
	\end{equation}
	Then 
	\[
		\bigg(\sum_{k=0}^\infty\Big(\sum_{n=k}^\infty b_{k,n}^2\Big)r^{2k}\bigg)^{\frac 12}
		\leq\max_{|z|=r}\Delta(z)\leq
		\sum_{k=0}^\infty\Big(\sum_{n=k}^\infty b_{k,n}^2\Big)^{\frac 12}r^k
		.
	\]
\end{Theorem}

\noindent
This two-sided estimate leads to formulae for order and type.

\begin{Corollary}
\label{U29}
	Consider the situation of \Cref{U28}. 
	\begin{Enumerate}
	\item Let $\rho$ and $\tau$ be order and type of the Nevanlinna matrix. Then 
		\[
			\rho=\limsup_{k\to\infty}\frac{2k\log k}{-\log\sum_{n=k}^\infty b_{k,n}^2}
			,
		\]
		and if $\rho>0$
		\[
			\tau=\frac 1{e\rho}\limsup_{k\to\infty}\bigg[k\Big(\sum_{n=k}^\infty b_{k,n}^2\Big)^{\frac\rho{2k}}\bigg]
			.
		\]
	\item Let $\ms f,\ms g$ be regularly varying functions which are asymptotic inverses of each other, and 
		denote by $\tau_{\ms g}$ the type of the Nevanlinna matrix w.r.t.\ $\ms g$.
		Then
		\[
			\tau_{\ms g}=\frac 1{e\rho}\limsup_{k\to\infty}
			\bigg[\ms f(k)^\rho\Big(\sum_{n=k}^\infty b_{k,n}^2\Big)^{\frac\rho{2k}}\bigg]
			.
		\]
	\end{Enumerate}
\end{Corollary}

\REMARKS{%
\item This result goes back to the work \cite{berg.szwarc:2014} of C.~Berg and R.~Szwarc, and the version stated as \Cref{U28} is 
	extracted from that paper. 
	An efficient proof of the theorem runs along the lines indicated before the statement. 

	\Cref{U29} is obtained simply by applying the standard formulae that express order and type of an entire function in
	terms of its Taylor coefficients. See \cite[Theorem~I.2]{levin:1980} for the statement in item (i), and
	\cite[Theorem~2']{levin:1980} (in conjunction with \cite[Proposition~7.4.1]{bingham.goldie.teugels:1989}) for 
	item (ii). Recall also \Cref{U201}.
\item In theory the question for the speed of growth of the Nevanlinna matrix is fully answered by \Cref{U28} and its
	corollary. Unfortunately this answer is of limited practical use. Applying the results requires precise knowledge about 
	all coefficients of all orthonormal polynomials, and additionally the ability to handle the series 
	$\sum_{n=k}^\infty b_{k,n}^2$; this is hardly possible. 
	We know of only one nontrivial situation where the formula for type from \Cref{U29} can be evaluated. 
	This will be presented in \Cref{U105}. 
\item \Cref{U28} can be used to give an alternative proof of an important (but actually much more elementary) lower
	bound for the growth of the Nevanlinna matrix, see the remarks in \Cref{U102}.
}

\section[{An elementary lower bound}]{An elementary lower bound}
\label{U102}

The below theorem contains a lower bound for the growth of the Nevanlinna matrix in terms of the Jacobi parameters $b_n$ 
which is easily applicable, elementary to prove, and sharp. 

\begin{Theorem}
\label{U24}
\IndexxS{Theorem!Nevanlinna matrix!lower bound}
	Let $a_n\in\bb R$ and $b_n>0$, and let $\ms J$ be the Jacobi matrix with these parameters. Assume that $\ms J$ is in limit
	circle case, and denote the corresponding Nevanlinna matrix as $W(z)=(w_{ij}(z))_{i,j=1}^2$. 
	Let $G$ be the function 
	\[
		G(z)\DE\sum_{n=1}^\infty\frac{z^{2n}}{(b_0\cdot\ldots\cdot b_{n-1})^2}
		.
	\]
	Then $G$ is entire, and (for, say, $r\geq 1$)
	\begin{equation}
	\label{U25}
		G(r)^{\frac 12}\lesssim |w_{11}(ir)|
	\end{equation}
	where the constant implicit in $\lesssim$ depends on $\ms J$.
\end{Theorem}

\noindent
If the parameters $b_n$ are dominated by a regularly varying function, the growth of $G(r)$ can be estimated explicitly. 

\begin{Corollary}
\label{U27}
	Let $a_n\in\bb R$ and $b_n>0$, and let $\ms J$ be the Jacobi matrix with these parameters. Assume that $\ms J$ is in limit
	circle case, and denote the corresponding Nevanlinna matrix as $W(z)=(w_{ij}(z))_{i,j=1}^2$.
	Let $\ms f,\ms g$ be regularly varying functions that are asymptotic inverses of each other. 
	Then 
	\[
		b_n\lesssim\ms f(n)\ \Longrightarrow\ \log|w_{11}(ir)|\gtrsim\ms g(r)
		.
	\]
\end{Corollary}

\noindent
We are led to lower bounds for order and type. 

\begin{Corollary}
\label{U26}
	Let $a_n\in\bb R$ and $b_n>0$, and let $\ms J$ be the Jacobi matrix with these parameters. 
	Assume that $\ms J$ is in limit circle case.
	\begin{Enumerate}
	\item Let $\rho$ and $\tau$ be order and type of the Nevanlinna matrix. Then 
		\begin{equation}
		\label{U30}
			\rho\geq\limsup_{n\to\infty}\frac{n\log n}{\sum_{j=0}^{n-1}\log b_j}
			,
		\end{equation}
		and if $\rho>0$ 
		\[
			\tau\geq \frac 1{e\rho}\limsup_{n\to\infty}\frac n{(b_0\cdot\ldots\cdot b_{n-1})^{\frac\rho n}}
			.
		\]
	\item Let $\ms f,\ms g$ be regularly varying functions which are asymptotic inverses of each other, and 
		denote by $\tau_{\ms g}$ the type of the Nevanlinna matrix w.r.t.\ $\ms g$.
		Then
		\[
			\tau_{\ms g}\geq \frac 1{e\rho}
			\limsup_{n\to\infty}\frac{\ms f(n)^\rho}{(b_0\cdot\ldots\cdot b_{n-1})^{\frac\rho n}}
			.
		\]
	\end{Enumerate}
\end{Corollary}

\noindent
\Cref{U27}, and with it \Cref{U24} and \Cref{U26}, is sharp. 
Properties which in conjunction ensure that the lower bound is attained are: 
a regular behaviour of the parameters, sufficently fast growth of the off-diagonal, and relative smallness of the diagonal. 
We state a simple version that can be easily formulated.

\begin{Proposition}
\label{U96}
	Let $\ms f$ be a regularly varying function with $\Ind\ms f>2$, and let $\ms g$ be an asymptotic inverse of $\ms f$. 
	Let $b_n>0$ be such that 
	\[
		b_n\asymp\ms f(n),\qquad \frac{b_n}{\sqrt{b_{n-1}b_{n+1}}}-1 \in \ell^1(\bb N)
		,
	\]
	set $a_n\DE 0$, and let $\ms J$ be the Jacobi matrix with these parameters. 
	Then $\ms J$ is in limit circle case and its Nevanlinna matrix satisfies
	\[
		\log\Big(\max_{|z|=r}\|W(z)\|\Big)\asymp\ms g(r)
		.
	\]
\end{Proposition}

\REMARKS{%
\item \Cref{U24} and \Cref{U27} are taken from \cite[Section~2]{pruckner.woracek:sinqA}. The proof is elementary:
	for the theorem one passes to the associated Hamburger Hamiltonian, and merely observes that a polynomial with real 
	coeffients can be estimated by its leading term; the corollary follows from the common connection between growth of an 
	entire function and its Taylor coefficients, e.g.\ \cite[Theorems~I.2,I.2']{levin:1980} 
	(in conjunction with \cite[Proposition~7.4.1]{bingham.goldie.teugels:1989}), 
\item A generalisation of \Cref{U27}, in terms of the associated Hamburger Hamiltonian, is given in 
	\cite[Theorem 3.6]{reiffenstein:kachamB-arXiv}, cf.\ \Cref{U95}. The statement of \Cref{U27} is the special case $k=2$
	in that result, as \cref{U213} and a short calculation shows.
\item \Cref{U96} follows by passing to Hamiltonian parameters and using
	\cite[Theorem~5.3]{pruckner.reiffenstein.woracek:sinqB-arXiv}. The translation, however, is not trivial. One can use
	\cite[Theorem~4.1]{yafaev:2020} which implies that $l_{n+1}=p_n(0)^2+q_n(0)^2 \asymp \frac{1}{b_n}$.

	Note further that the condition $\frac{b_n}{\sqrt{b_{n-1}b_{n+1}}}-1 \in \ell^1(\bb N)$ cannot be dropped. 
	As an example, take
	\begin{align*}
			a_n\DE 0,\qquad b_n\DE
			\begin{cases}
				(n+1)^{\frac 52} (n+2)^{\frac 12} \CAS n \text{ odd}
				\\
				(n+1)^{\frac 12} (n+2)^{\frac 52} \CAS n \text{ even}.
			\end{cases}
	\end{align*}
	Clearly $b_n \sim n^3$. However, the corresponding Hamiltonian parameters are $\phi_n=n\frac{\pi}{2}$ and
	\[
		l_{2n-1}=n^{-1}, \qquad l_{2n}=n^{-5}
	\]
	which shows that the Jacobi matrix is in limit point case.
\item The estimates in \Cref{U26}\,(i) were first given in \cite[Proposition~7.1]{berg.szwarc:2014} with a proof that 
	uses \Cref{U29} (just drop all terms in the sum $\sum_{n=k}^\infty b_{k,n}^2$ but the first one, 
	and remember \cite[p.86]{akhiezer:1965}).
	Using the stated argument for the proof of \Cref{U24}, one obtains an elementary direct proof of \Cref{U26}, 
	cf.\ \cite[Remark~2.15]{pruckner.romanov.woracek:jaco}. 
\item Despite giving the correct order for a wide range of Jacobi matrices, the bound \cref{U25} may also fail drastically,
	even on the rough scale of order as in \cref{U30}.
	One can say that this is expected; remember the alternative proof by dropping -- a lot of -- summands. 
	The following example is obtained by translating \cite[Theorem~3.1,Corollary~3.6]{pruckner.romanov.woracek:jaco} to the
	language of Jacobi parameters.

	Let $\rho\in(0,1]$ and $r\in(0,\rho)$, and set 
	\[
		a_n\DE 0,\qquad b_n\DE
		\begin{cases}
			1 \CAS n=0,1,
			\\
			(n\log^2n)^{\frac 1\rho} \CAS n\equiv 0,1\!\mod 4,\ n>2,
			\\
			(n\log^2n)^{\frac 1\rho+2(\frac 1r-\frac 1\rho)} \CAS n\equiv 2,3\!\mod 4,\ n>2.
		\end{cases}
	\]
	Then $\ms J$ is in limit circle case and the order of the Nevanlinna matrix is equal to $\rho$, 
	while the limit superior in \cref{U30} is equal to $r$. 
}

\section[{The Liv\v sic estimate}]{The Liv\v sic estimate}
\label{U104}

Already back in 1939 M.S.~Liv\v sic gave a lower bound for the Nevanlinna matrix in terms of the moment sequence 
$(s_n)_{n=0}^\infty$ itself. From a nowadays viewpoint, his theorem appears as a consequence of the results from \Cref{U102}.

\begin{Theorem}
\label{U15}
\IndexxS{Theorem!Nevanlinna matrix!Liv\v sic}
	Let $(s_n)_{n=0}^\infty$ be an indeterminate moment sequence, and denote the corresponding Nevanlinna matrix as
	$W(z)=(w_{ij}(z))_{i,j=1}^2$. Let $F$ be the \IndexS{Liv\v sic function}{Liv\v sic function} 
	\[
		F(z)\DE\sum_{n=0}^\infty\frac{z^{2n}}{s_{2n}}
		.
	\]
	Then $F$ is entire, and (for, say, $r\geq 1$)
	\[
		F(r)^{\frac 12}\lesssim |w_{11}(ir)|
		.
	\]
\end{Theorem}

\noindent
The analogue of \Cref{U26} reads as follows. 

\begin{Corollary}
\label{U93}
	Let $a_n\in\bb R$ and $b_n>0$, and let $\ms J$ be the Jacobi matrix with these parameters. 
	Assume that $\ms J$ is in limit circle case.
	\begin{Enumerate}
	\item Let $\rho$ and $\tau$ be order and type of the Nevanlinna matrix. Then 
		\begin{equation}
		\label{U22}
			\rho\geq\limsup_{n\to\infty}\frac{2n\log n}{\log s_{2n}}
			,
		\end{equation}
		and if $\rho>0$ 
		\[
			\tau\geq \frac 1{e\rho}
			\limsup_{n\to\infty}\frac n{(s_{2n})^{\frac\rho{2n}}}
			.
		\]
	\item Let $\ms f,\ms g$ be regularly varying functions which are asymptotic inverses of each other, and 
		denote by $\tau_{\ms g}$ the type of the Nevanlinna matrix w.r.t.\ $\ms g$.
		Then
		\[
			\tau_{\ms g}\geq \frac 1{e\rho}
			\limsup_{n\to\infty}\frac{\ms f(n)^\rho}{(s_{2n})^{\frac\rho{2n}}}
			.
		\]
	\end{Enumerate}
\end{Corollary}

\REMARKS{%
\item The original reference for \Cref{U15} is \cite[Theorema~2]{livshits:1939}. More accessible sources are
	\cite[\S7]{berg.szwarc:2014} or \cite[p.224]{pruckner.romanov.woracek:jaco}. 
\item The argument for the proof of \Cref{U15} and \Cref{U93} from \cite{pruckner.romanov.woracek:jaco} is by reduction to 
	\Cref{U24} and \Cref{U26}: noting that $(b_0\cdot\ldots\cdot b_{n-1})^{-1}$ is the leading coefficient of the 
	orthonormal polynomial of degree $n$, yields
	\begin{equation}
	\label{U16}
		\frac 1{\sqrt{s_{2n}}}\leq \frac 1{b_0\cdot\ldots\cdot b_{n-1}}
		,
	\end{equation}
	cf.\ \cite[p.86]{akhiezer:1965} and \cite[(3.8)]{pruckner.romanov.woracek:jaco}.
	\Cref{U93} readily follows from \Cref{U26}. To deduce \Cref{U15} from \Cref{U24}, observe that $F(r)\leq G(r)$, 
	where $G$ is as in \Cref{U24}.
\item The Liv\v sic bound \cref{U22} is sharp. In fact, it gives the correct value for $\rho$ in a large class of
	examples; see the remarks to \Cref{U107}.
\item From \cref{U16} and the remarks in \Cref{U102} we obtain that Liv\v sic's bound \cref{U22} may fail drastically. 
	Given $\rho\in(0,1]$ and $r\in(0,\rho)$, there exists an indeterminate moment sequence such that the order of the 
	Nevanlinna matrix is $\rho$ while 
	\begin{equation}
	\label{U17}
		\limsup_{n\to\infty}\frac{2n\log n}{\log s_{2n}}\leq r
		.
	\end{equation}
	We expect, but do not know, that the example can be chosen such that equality holds in \cref{U17}. 
}

\section[{Convergence exponents as lower bounds}]{Convergence exponents as lower bounds}
\label{U152}

Consider a Jacobi matrix in limit circle case, i.e., corresponding to an indeterminate moment problem. We present two lower bounds for the growth of the Nevanlinna matrix in terms of convergence exponents, which are easy to
apply and surprisingly often give the correct order. We also give more quantitative lower bounds, which are less flexible but still widely applicable. Contrasting the bound exhibited in \Cref{U102}, the proofs are not
elementary.

Recall that, since we assume limit circle case, we have
\begin{align}
\sum_{n=1}^\infty\frac 1{b_n}<\infty \quad \text{and} \quad \sum_{n=1}^\infty\frac {|a_{n+1}|}{b_n b_{n+1}}<\infty.
\end{align}
This is because either of the conditions
\begin{align*}
	& \sum_{n=1}^\infty\frac 1{b_n}=\infty 
	&& \text{ {\it(}\IndexS{Carleman's condition}{Carleman's condition}{\it)}}  
	\\
	& \sum_{n=1}^\infty\frac {|a_{n+1}|}{b_n b_{n+1}}=\infty 
	&& \text{ {\it(}\IndexS{Dennis--Wall condition}{Dennis--Wall condition}{\it)}}
\end{align*}
implies limit point case.
Consequently, assuming limit circle case both sequences
\begin{align}
\label{U94}
(b_n)_{n=1}^\infty, \quad \Big(b_n b_{n+1}/\sqrt{a_{n+1}^2+b_n^2+b_{n+1}^2} \Big)_{n=0}^{\infty}
\end{align}
have finite convergence exponent
not exceeding $1$. 

\begin{Theorem}
\label{U97}
\IndexxS{Theorem!Nevanlinna matrix!lower bound}
	Let $a_n\in\bb R$ and $b_n>0$, and let $\ms J$ be the Jacobi matrix with these parameters. 	Assume that $\ms J$ is in limit circle case. Then the order of its Nevanlinna matrix is not less than either of the convergence exponents of the sequences in \eqref{U94}.
\end{Theorem}

\noindent There are also more quantitative versions of this lower bound -- a sequential one as well as one that holds for all 
large $r$, but with a less explicit right-hand side.

\begin{Theorem}
\label{U98}
\IndexxS{Theorem!Nevanlinna matrix!lower bound}
	Let $a_n\in\bb R$ and $b_n>0$, and let $\ms J$ be the Jacobi matrix with these parameters. 
	Assume that $\ms J$ is in limit circle
	case, and denote the corresponding Nevanlinna matrix as $W(z)=(w_{ij})_{i,j=1}^2$. Then
	\begin{Enumerate}
	\item Suppose we have functions $\ms f,\ms g\DF(0,\infty)\to(0,\infty)$ with $\sum_{j=1}^\infty \ms f(b_j)=\infty$
		and $\sum_{j=1}^\infty \ms g(1/j)<\infty$. Then there exists an increasing and unbounded sequence
		$(r_m)_{m=1}^\infty$ of positive numbers with
		\begin{align*}
			\log |w_{21}(ir_m)| \gtrsim \frac{1}{\ms g^- (\ms f(r_m))}, \qquad m \in \bb N,
		\end{align*}
		where $\ms g^-(t) \DE \sup \{r>0\DS\ms g(r)<t\}$.
	\item For all large enough $r>0$ we have
		\[
			\log|w_{21}(ir)|\gtrsim r\sum_{n=h(r)}^{\infty}\frac{1}{b_n},
		\]
		where $h(r) \DE 1+\max \{n \in \bb N \DF b_n < r \}$.
	\end{Enumerate}	
	Both items remain true when $(b_n)_{n=0}^{\infty}$ is replaced with the sequence 
	$\big(b_n b_{n+1}/\sqrt{a_{n+1}^2+b_n^2+b_{n+1}^2}\big)_{n=0}^\infty$.
\end{Theorem}

\REMARKS{%
\item The fact that Carleman's condition implies limit point case is a classical result from \cite{carleman:1926}, see also 
	\cite[p.24(1$^\circ$)]{akhiezer:1965}. A one-line proof can be given by passing to the Hamiltonian parameters:
	we have $\frac 1{b_n}\leq\sqrt{l_nl_{n+1}} \leq \frac 12 (l_n+l_{n+1})$. 
	For the Dennis--Wall condition see \cite[p.25]{akhiezer:1965}

	Let us note in this context that this implies a necessary condition for indeterminacy which is in terms
	of the moment sequence itself: $\sum_{n=0}^\infty\frac 1{\sqrt[2n]{s_{2n}}}<\infty$
	whenever $\ms J$ is in limit circle case, cf.\ \cite[p.85(11$^\circ$)]{akhiezer:1965}.
\item \Cref{U97} and \Cref{U98} are taken from work \cite{reiffenstein:kachamB-arXiv} of J.~Reiffenstein. After
	passing over to the associated Hamburger Hamiltonian, they appear as particular cases of more general theorems that
	rely on the Weyl coefficient approach, cf.\ \Cref{U110,U132}. 
	This makes \Cref{U97,U98} intrinsically different from the lower bound discussed in \Cref{U102} and 
	Liv\v sic's bound in \Cref{U15}, since those results are shown by an elementary estimate of polynomials.
\item The two sequences in \eqref{U94} are just the first and second in an infinite scale of sequences whose convergence
	exponents are all lower bounds for the monodromy matrix, cf. \Cref{U31}. 
\item In many situations the order of the Nevanlinna matrix coincides with the convergence exponent of $(b_n)_{n=1}^\infty$.
	Large classes of Jacobi matrices with this property are presented in \Cref{U107,U120}. In these cases it is clear that
	the convergence exponents of the two sequences in \eqref{U94} coincide. Nonetheless, there are examples where these
	sequences have different convergence exponents and the order of the Nevanlinna matrix is equal to the convergence
	exponent of the second sequence in \eqref{U94}; see the remarks to \Cref{U31}.
\item The convergence exponent of $(b_n)_{n=1}^\infty$ does not at all capture the parameters $a_n$ and neither takes into
	account sparse occurrences of small $b_n$. Hence, it is somewhat obvious that it cannot always give the correct value
	for the order of the Nevanlinna matrix. A simple example that shows this is given by the parameters 
	\[
		b_n\DE n^{\frac 53},\quad a_n\DE 2n^{\frac 53}\big(1-\frac 1n\big)
		,
	\]
	for which limit circle case takes place, the convergence exponent of $(b_n)_{n=0}^\infty$ is $\frac 35$, but the order
	of the Nevanlinna matrix is $\frac 34$. Note that in this example the second sequence in \eqref{U94} 
	also has convergence exponent $\frac 35$.

	This example is just one instance of \Cref{U49}; see the remarks to \Cref{U120}. 
	The key property, which makes this phenomenon possible, is that diagonal and off-diagonal are of almost the same size 
	with quotient $2$.  Intuitively this makes the Jacobi matrix ``close to limit point case'' and allows the Nevanlinna
	matrix to behave in an unusual way. 
\item The order $\rho$ of the Nevanlinna matrix of a Jacobi operator in limit circle case can be expressed via the eigenvalue
	counting function $n(r)$ of one of its self-adjoint extensions as 
	\[
		\rho=\limsup_{r \to \infty} \frac{\log n(r)}{\log r}
		.
	\] 
	If we consider instead a Jacobi operator in limit point case whose spectrum is still discrete, we may ask whether this 
	$\limsup$, defined using the eigenvalue counting function of the Jacobi operator, is again bounded below by the 
	convergence exponent of $(b_n)_{n=0}^\infty$. 

	The following example from \cite[Theorem 2.2]{janas.malejki:2007} shows that the answer to this
	question is negative:
	Let $a_n \DE n^\alpha$, $b_n \DE n^\beta$, where $\beta \geq 0$ and $\alpha \geq 2\beta +1$. Then the eigenvalues
	$\lambda_n$, arranged increasingly, of the Jacobi matrix with parameters $a_n,b_n$ satisfy 
	$\lambda_n = n^\alpha+\BigO(n^{1+2\beta-\alpha})$.
	Intuitively, the Jacobi operator with these parameters is close to a diagonal operator, whose eigenvalues are precisely
	its diagonal entries.

\item Under additional assumptions on the regularity of $b_n$ the bound from item (ii) of \Cref{U98} can be evaluated. For
	example, if $b_n\asymp\ms g(n)$ with some regularly varying function, then the bound is $\gtrsim\ms f(r)$ where $\ms f$
	is an asymptotic inverse of $\ms g$. Thus we reobtain a particular case of \Cref{U27}.
\item Let us compare the lower bounds for order from \Cref{U97} and \cref{U30}. 
	Denote the convergence exponent of $(b_n)_{n=0}^\infty$ by $\gamma$, then
	\[
		\gamma\leq\limsup_{n\to\infty}\frac{\log n}{\log b_n}\geq
		\limsup_{n\to\infty}\frac{n\log n}{\sum_{j=0}^n\log b_j}
		,
	\]
	i.e., the term in the middle is an estimate for both of our lower bounds.
	If $(b_n)_{n=0}^\infty$ is nondecreasing, then the first inequality holds with equality, and hence the bound from 
	\Cref{U97} is potentially better. If, for example, the limit superior in the middle exists as a limit, then equality
	holds throughout, i.e., the bounds from \Cref{U97} and \cref{U30} coincide. 
}

\section[{A theorem of Berezanskii}]{A theorem of Berezanskii}
\label{U107}

We exhibit a particular class of Jacobi matrices for which the general lower bound from \Cref{U97} is attained. 
In short, the theorem says: if the off-diagonal $b_n$ behaves regularly and the diagonal $a_n$ is relatively small, 
then the negation of Carleman's condition is sufficient for $\ms J$ to be in limit circle case and the order of the Nevanlinna matrix to be 
equal to the convergence exponent of $(b_n)_{n=0}^\infty$. 

\begin{Theorem}
\label{U12}
\IndexxS{Theorem!Berezanskii}
	Let $a_n\in\bb R$ and $b_n>0$, and let $\ms J$ be the Jacobi matrix with these parameters. 
	Set 
	\[
		\alpha_n\DE\frac{a_n}{\sqrt{b_{n-1}b_n}}
		,
	\]
	and assume that the following three conditions hold.
	\begin{Itemize}
	\item \IndexS{Carleman's condition}{Carleman's condition} violated: 
		\[
			\sum_{n=0}^\infty\frac 1{b_n}<\infty
			.
		\]
	\item Relative smallness and regularity of the diagonal:
		\[
			\sum_{n=1}^\infty|\alpha_{n+1}-\alpha_n|<\infty,\qquad \lim_{n\to\infty}\alpha_n\in(-2,2)
			.
		\]
	\item Regularity of the off-diagonal:
		\begin{equation}
		\label{U51}
			\sum_{n=0}^\infty\Big|\frac{b_n}{\sqrt{b_{n-1}b_{n+1}}}-1\Big|<\infty
			.
		\end{equation}
	\end{Itemize}
	Then $\ms J$ is in limit circle case, and the order of the Nevanlinna matrix of $\ms J$ is equal to the convergence
	exponent of the sequence $(b_n)_{n=0}^\infty$. 
\end{Theorem}

\REMARKS{%
\item \Cref{U12} is taken from \cite{reiffenstein:kachamB-arXiv}. It is a significant improvement of a classical result 
	due to Yu.M.~Berezanskii in \cite{berezanskii:1956} who assumed a stronger regularity of $b_n$ and that $a_n$ are 
	bounded. A more accessible reference is \cite{berg.szwarc:2014}, where a variant improving upon the original result 
	is given. It states that the conclusion of \Cref{U12} follow under the assumptions that 
	$\sum_{n=1}^\infty\frac{1+|a_n|}{\sqrt{b_nb_{n-1}}}<\infty$ (meaning that Carleman's condition is violated and 
	the diagonal is small), and that the sequence $(b_n)_{n=0}^\infty$ is log-concave or log-convex (meaning that 
	$b_n^2\geq b_{n-1}b_{n+1}$ for all $n$, or $b_n^2\leq b_{n-1}b_{n+1}$ for all $n$, respectively). 
\item Under the hypothesis of Berezanskii's theorem (in the variant from \cite{berg.szwarc:2014}) 
	the Liv\v sic's bound \cref{U22} coincides with the 
	convergence exponent of $(b_n)_{n=0}^\infty$, and hence gives the correct value for the order of the Nevanlinna matrix. 
	This is shown in \cite[Theorem~7.5]{berg.szwarc:2014}.
\item Let $\rho$ be the convergence exponent of $(b_n)_{n=0}^\infty$. If $\sum_{n=0}^\infty\frac 1{b_n^\rho}<\infty$, then the
	Nevanlinna matrix is of finite type w.r.t.\ its order $\rho$. Under the conditions required in \cite{berg.szwarc:2014}
	one also has 
	\[
		\frac 1{b_n}=\BigO\big(n^{-\frac 1\rho}\big)
		.
	\]
}

\section[{Growth from power asymptotics}]{Growth from power asymptotics}
\label{U120}

Consider parameters $a_n\in\bb R$ and $b_n>0$ having power asymptotics of the form
\begin{equation}
\label{U50}
	b_n=n^{\beta_1}\Big(x_0+\frac{x_1}n+\BigO\Big(\frac 1{n^{1+\epsilon}}\Big)\Big),\quad 
	a_n=n^{\beta_2}\Big(y_0+\frac{y_1}n+\BigO\Big(\frac 1{n^{1+\epsilon}}\Big)\Big)
\end{equation}
where 
\[
	\beta_1,\beta_2\in\bb R,\quad x_0>0,y_0\neq 0,\quad x_1,y_1\in\bb R,\quad \epsilon>0
	.
\]
For this class of Jacobi matrices, one can almost fully decide whether $\ms J$ is in limit circle case, and if it is, what is the 
order of the Nevanlinna matrix. We should say explicitly that we do not investigate the limit point case, apart from giving
conditions under which limit point case prevails. 

There are three essentially different cases.
\begin{Itemize}
\item Large diagonal: 
	\[
		\beta_2>\beta_1\quad\text{or}\quad\Big(\beta_1=\beta_2\wedge|y_0|>2x_0\Big)
		.
	\]
\item Small diagonal:
	\[
		\beta_2<\beta_1\quad\text{or}\quad\Big(\beta_1=\beta_2\wedge|y_0|<2x_0\Big)
		.
	\]
\item Critical case: 
	\[
		\beta_2=\beta_1\quad\text{and}\quad |y_0|=2x_0
		.
	\]
\end{Itemize}
We state two theorems. The first settles the cases of large and small diagonal, and the second goes into the critical case.

\begin{Theorem}
\label{U48}
	Let $a_n\in\bb R$ and $b_n>0$ be as in \cref{U50} and let $\ms J$ be the Jacobi matrix with these parameters.
	\begin{Enumerate}
	\item If $\ms J$ has ``large diagonal'', then $\ms J$ is in limit point case. 
	\item If $\ms J$ has ``small diagonal'', then $\ms J$ is in limit circle case if and only if $\beta_1>1$. 
	\item If $\ms J$ has ``small diagonal'' and is in limit circle case, then 
		the Nevanlinna matrix is of order $\frac 1{\beta_1}$ with finite and positive type. 
	\end{Enumerate}
\end{Theorem}

\noindent
In order to handle the critical case, where off-diagonal and diagonal entries of $\ms J$ are comparable with ratio $\pm 2$,
we require one more term in the asymptotic expansion of the parameters:
\begin{equation}
\label{U92}
	b_n=n^{\beta_1}\Big(x_0+\frac{x_1}n+\frac{x_2}{n^2}+\BigO\Big(\frac 1{n^{2+\epsilon}}\Big)\Big),\quad
	a_n=n^{\beta_2}\Big(y_0+\frac{y_1}n+\frac{y_2}{n^2}+\BigO\Big(\frac 1{n^{2+\epsilon}}\Big)\Big)
\end{equation}
where 
\[
	\beta_1,\beta_2\in\bb R,\quad x_0>0,y_0\neq 0,\quad x_1,x_2,y_1,y_2\in\bb R,\quad \epsilon>0
	.
\]

\begin{Theorem}
\label{U49}
\IndexxS{Theorem!Jacobi parameters with power asymptotics}
	Let $a_n\in\bb R$ and $b_n>0$ be as in \cref{U92} and let $\ms J$ be the Jacobi matrix with these parameters.
	Assume that $\ms J$ is in the ``critical case''.
	\begin{Enumerate}
	\item $\ms J$ is in limit circle case, if and only if either
		\[
			\frac 32<\beta_1<2\Big(\frac{x_1}{x_0}-\frac{y_1}{y_0} \Big)
			,
		\]
		or
		\[
			\beta_1 = 2 \Big(\frac{x_1}{x_0}-\frac{y_1}{y_0} \Big) \ \wedge \ 
			1<\frac{x_1}{x_0}-\frac{y_1}{y_0}<\frac 34+
			\frac{2x_2}{x_0}-\frac{2y_2}{y_0}+\frac{y_1}{y_0}-\frac{2x_1y_1}{x_0y_0}+\frac{2y_1^2}{y_0^2}
			.
		\]
	\item Assume that $\ms J$ is in limit circle case.
		\begin{Itemize}
		\item If $\beta_1>2$, then the Nevanlinna matrix is of order $\frac 1{\beta_1}$ with finite and positive type. 
		\item If $\beta_1<2$, then the Nevanlinna matrix is of order $\frac 1{2(\beta_1-1)}$ 
			with finite and positive type. 
		\item If $\beta_1=2$, then the Nevanlinna matrix is of order $\frac 12$ with positive type. 
			It has finite type w.r.t.\ the regularly varying function $r^{\frac 12}\log r$.
		\end{Itemize}
	\end{Enumerate}
\end{Theorem}

\REMARKS{%
\item The above results are obtained by combining work of R.~Pruckner \cite{pruckner:blubb} and J.~Reiffenstein
	\cite{reiffenstein:kachamA-arXiv}, and referring to \cite{pruckner.reiffenstein.woracek:sinqB-arXiv}. 
\item We used the assumptions \cref{U50} and \cref{U92} on the asymptotic expansion of parameters to increase readability, but
	in some cases they can be relaxed. For details we refer to \cite[Theorem~3.4]{reiffenstein:kachamA-arXiv}. 
\item The proofs of \Cref{U48} and \Cref{U49} are quite different. What they have in common is that the starting point is an
	asymptotic analysis of the solutions of the three-term recurrence, which is more complicated in the second theorem than
	in the first. In \Cref{U48} the	bounds come from an explicit analysis of the fundamental solution, while for \Cref{U49} 
	they are obtained by passing to the Hamiltonian parameters. 
\item The proof of \Cref{U48} yields explicit bounds for the type $\tau$ of the Nevanlinna matrix, and by complex analysis
	thus bounds for the counting function of the spectrum of $\ms J$. In the ``critical case'' only a lower bound is given 
	explicitly:
	\begin{Itemize}
	\item If $\ms J$ has parameters \cref{U50}, ``small diagonal'', and is in limit circle case, then 
		\[
			\tau\geq\beta_1\Big(\frac 1{x_0}\Big)^{\frac 1{\beta_1}}
			,\qquad 
			\limsup_{r\to\infty}\frac{n_{\ms J}(r)}{n^{\frac 1{\beta_1}}}\geq
			\frac{\beta_1-1}{\beta_1}\Big(\frac 1{x_0}\Big)^{\frac 1{\beta_1}}
			,
		\]
		where $\IndexN{n_{\ms J}(r)}$ is the number of eigenvalues of absolute value less than $r$ of the Jacobi operator.
		If $\ms J$ has parameters \cref{U92}, is in the ``critical case'' and in limit circle case, the same
		lower bounds hold.
	\item If $\ms J$ has parameters \cref{U50}, ``small diagonal'', and is in limit circle case, then 
		\[
			\tau\leq\frac{\pi}{\sin(\frac\pi{\beta_1})}\Big(\frac a{x_0}\Big)^{\frac 1{\beta_1}}
			,\qquad
			\limsup_{r\to\infty}\frac{n_{\ms J}(r)}{n^{\frac 1{\beta_1}}}
			\leq\frac{e\pi}{\beta_1\sin(\frac\pi{\beta_1})}\Big(\frac a{x_0}\Big)^{\frac 1{\beta_1}}
			,
		\]
		where 
		\[
			a\DE
			\begin{cases}
				1 \CAS \beta_2<\beta_1,
				\\
				\big(1-\frac{y_0^2}{4x_0^2}\big)^{-\frac 12} \CAS \beta_1=\beta_2\wedge|y_0|<2x_0.
			\end{cases}
		\]
	\end{Itemize}
\item From a very recent perspective, \Cref{U48}\,(ii),(iii) are a consequence of Berezanskii's theorem as stated in \Cref{U12}.
	When \Cref{U48} was originally proved in \cite{pruckner:blubb}, this version of Berezanskii's theorem was not yet
	available. Furthermore, Pruckner's work features the type estimates presented above, which \Cref{U12} does not provide. 
\item Consider the ``critical case'' as in \Cref{U49}. Then the convergence exponent of $(b_n)_{n=1}^\infty$ is 
	$\frac 1{\beta_1}$. If $\beta_1<2$, the order of the Nevanlinna matrix is $\frac 1{2(\beta_1-1)}$, and hence 
	strictly larger than this convergence exponent (this yields the explicit example given in the remarks to
	\Cref{U152}).
	Up to our knowledge, \Cref{U49} yields the first examples of Jacobi matrices in limit circle case whose order
	can be computed and is different from the convergence exponent of $(b_n)_{n=1}^\infty$.
\item We do not know whether the results presented above are already close to being sharp when it comes to assumptions on the
	regularity of the parameters, or if the hypothesis that they possess power asymptotics with a certain number of terms
	(or as in \cite{reiffenstein:kachamA-arXiv}) can be weakened significantly.
\item In the limit point case there is a vast literature about spectral asymptotics (but we do not go into this direction). 
	We do not know results that make assertions about the distribution of the spectrum in limit point case without knowing
	actual asymptotics.
}

%**************************************************************************
%***                            Last Change: Mon 31 Mar 2025 11:27
%***   < PART IV >
%***
%**************************************************************************

\clearpage
\PART{IV}{Hamburger Hamiltonians}
\label{U147}

%%%%%%%%%%%%%%%%%%%%%%%%%%%%%%%%%%%%%%%%%%%%%%%%%%%%%%%%%%%%%%%%%%%%%%%%%%%%%%%%%%%%%%%%%

\Intro{%
	In this part we change perspectives towards the study of Hamburger Hamiltonians. That is, we look at Hamburger moment
	problems and Jacobi matrices from the canonical systems point of view. As in Part~III we focus on the indeterminate case, 
	in the present language the limit circle case. 
	It turns out that the parameters of the Hamiltonian are often better suited for determining
	the growth of the monodromy matrix than the moment sequence or the Jacobi parameters. 

	A Hamburger Hamiltonian $H$ is given by the sequences $(l_n)_{n=1}^\infty$ and $(\phi_n)_{n=1}^\infty$ of its lengths
	and angles. On an intuitive level one may say that the growth of the monodromy matrix $W_H$ is given by a combination 
	of three factors: the rate of decay of the lengths $l_j$, the decay of the angle-differences $|\sin(\phi_{j+1}-\phi_j)|$, and 
	the speed of possible convergence of the angles, each measured pointwise or in an averaged sense. 
	Fast decay of lengths and angle-differences and fast convergence of angles means slow growth of the monodromy matrix. 
	We present a variety of results giving upper and lower bounds for $W_H$, each of which gives quantitative meaning to the
	above intuition in its own way. 

	There is a distinct threshold between quickly and slowly decaying lengths and angle-differences, corresponding to the
	order of $W_H$ being less than and larger than $\frac 12$, respectively. In the case of fast decay (and some regularity
	of the data) the growth of $W_H$ is determined by the lengths and angle-differences only. Contrasting this, in the case
	of slow decay the growth of $W_H$ is also sensible to the speed of possible convergence of angles. In particular, there
	are lower and upper bounds for $W_H$ in terms of the lengths and angle-differences, which are distinct even on the level of order and can both be attained.
	}
\begin{center}
	{\large\bf Table of contents}
\end{center}
\begin{flushleft}
	
	\S.\,\ref{U31}\ A scale of lower bounds\ \dotfill\quad\pageref{U31}
	\\[1mm]
	\S.\,\ref{U133}\ Upper bounds in terms of convergence exponents\ \dotfill\quad\pageref{U133}
	\\[1mm]
	\S.\,\ref{U160}\ Upper bounds in terms of tails of convergent series\ \dotfill\quad\pageref{U160}
	\\[1mm]
	\S.\,\ref{U135}\ A Berezanskii-type theorem\ \dotfill\quad\pageref{U135}
	\\[1mm]
	\S.\,\ref{U119}\ Exploiting regular variation I: upper bound from majorants\ \dotfill\quad\pageref{U119}
	\\[1mm]
	\S.\,\ref{U150}\ Exploiting regular variation II: lower bound from minorants\ \dotfill\quad\pageref{U150}
	\\[1mm]
\end{flushleft}
\makeatother
\renewcommand{\thesection}{\arabic{section}}
\renewcommand{\thelemma}{\arabic{section}.\arabic{lemma}}
\makeatletter
\clearpage

%%%%%%%%%%%%%%%%%%%%%%%%%%%%%%%%%%%%%%%%%%%%%%%%%%%%%%%%%%%%%%%%%%%%%%%%%%%%%%%%%%%%%%%%%

%
%
%
\section[{A scale of lower bounds}]{A scale of lower bounds}
\label{U31}

Consider a Hamburger Hamiltonian in limit circle case given by lengths $l_j$ and angles $\phi_j$, and recall the notation $x_n\DE\sum_{j=1}^{n}l_j$ and $W_H=(w_{H,ij})_{i,j=1}^2$.

For each $k\in\bb N$, $k\geq 2$, we define a sequence $(b_n^{(k)})_{n=0}^\infty$ of positive numbers as 
\begin{equation}
\label{U224}
	\IndexN{b_n^{(k)}}\DE\Big[\frac 12\sum_{i,j=n+1}^{n+k} l_il_j \sin^2(\phi_i-\phi_j)\Big]^{-\frac 12}
	.
\end{equation}
We have $\sum_{n=1}^\infty \frac 1{b_n^{(k)}} <\infty$ for all $k$ and $b_n^{(k)}\geq b_n^{(k+1)}$ for all $k$ and $n$. 
In particular, the convergence exponent of the sequence $(b_n^{(k)})_{n=0}^\infty$ is at most $1$ and depends nondecreasingly 
on $k$. 

\begin{Theorem}
\label{U208}
\IndexxS{Theorem!monodromy matrix!lower bound (Hamburger Hamiltonian)}
	Let $l_j>0$ and $\phi_j\in\bb R$ with $\phi_{j+1}-\phi_j\not\equiv 0\mod\pi$, 
	and let $H$ be the Hamburger Hamiltonian with these lengths and angles.
	Assume that $H$ is in limit circle case. 
	Then, for each $k\geq 2$, the order
	of the monodromy matrix of $H$ is not less than the convergence exponent of the sequence $(b_n^{(k)})_{n=0}^\infty$. 
\end{Theorem}

\noindent
In the following theorem we collect two versions of this lower bound that estimate $\log |w_{H,22}(ir)|$ in a more quantitative
way. The first bound holds on a sequence $r_m \to \infty$ whereas the second is true for all large $r$.

\begin{Theorem}
\label{U14}
\IndexxS{Theorem!monodromy matrix!lower bound (Hamburger Hamiltonian)}
	Let $l_j>0$ and $\phi_j\in\bb R$ with $\phi_{j+1}-\phi_j\not\equiv 0\mod\pi$, 
	and let $H$ be the Hamburger Hamiltonian with these lengths and angles.
	Assume that $H$ is in limit circle case. 
	The following statements hold true for any $k \geq 2$.
	\begin{Enumerate}
	\item Let $\ms f, \ms g \DF (0,\infty) \to (0,\infty)$ with $\sum_{n=1}^\infty \ms f(b_n^{(k)})=\infty$ and 
		$\sum_{n=1}^\infty\ms g(1/n)<\infty$. Then there exists an increasing and unbounded sequence 
		$(r_m)_{m=1}^\infty$ of positive numbers with
		\[
			\log |w_{H,22}(ir_m)| \gtrsim \frac{1}{\ms g^- (\ms f(r_m))}, \qquad m \in \bb N,
		\]
		where $\ms g^-(t) \DE \sup \{r>0 \DF \, \ms g(r) <t\}$.
	\item Assume we have a sequence $(f_n)_{n=0}^\infty$ of
		nonnegative numbers such that 
		\[
			b_n^{(k)}\lesssim \frac 1{f_n}
			.
		\]
		Then (for, say, $r\geq 1$) 
		\[
			\log|w_{H,22}(ir)|\gtrsim r	\sum_{n=h(r)}^{\infty} f_n
		\]
		where 
\[
h(r) \DE 1+\max \big\{n \in \bb N \DF \, f_n > r^{-1} \big\}.
\]
	\end{Enumerate}
\end{Theorem}

\REMARKS{%
\item The two statements in \Cref{U14} are related but in general incomparable. Both are taken from work of J.~Reiffenstein
	\cite{reiffenstein:kachamB-arXiv}. The proofs are different in their nature, but both of them use the
	Weyl coefficient approach, cf.\ \Cref{U110,U132}. 
\item The above results should be seen as a scale of bounds emanating from the lower bounds in terms of the Jacobi parameters
	$b_n$ that we discussed in \Cref{U152}. In fact,
	\begin{align}
	\label{U225}
		b_n^{(2)}=\frac{1}{\sqrt{l_{n+2}l_{n+1}}|\sin(\phi_{n+2}-\phi_{n+1})|}=b_n
	\end{align}
	and a computation shows that 
	\begin{align}
		b_n^{(3)}=\frac{b_n b_{n+1}}{\sqrt{a_{n+1}^2+b_n^2+b_{n+1}^2}}.
	\end{align}
	We see that \Cref{U97} matches the particular cases ``$k=2,3$'' in \Cref{U208}. Analogously, \Cref{U98} follows from
	\Cref{U14}. 
\item By comparing \eqref{U224} with \eqref{U218} we have
	\begin{align*}
		b_n^{(k)}=\frac{1}{\sqrt{\det \Omega_H (x_n,x_{n+k})}}.
	\end{align*} 
\item In many situations the correct value for the order $\rho$ of the monodromy matrix is already given by the
	convergence exponent of $(b_n^{(2)})_{n=0}^\infty$. Due to \Cref{U208} and monotonicity in $k$ of the convergence
	exponents of $(b_n^{(k)})_{n=0}^\infty$, in such a situation the convergence exponents all coincide. 
\item An example where not all sequences $(b_n^{(k)})_{n=0}^\infty$ have the same convergence exponent given in 
	\cite[Example~3.4]{reiffenstein:kachamB-arXiv}. Let $\alpha_1>\alpha_0>1$, and set 
	\[
		l_j\DE
		\begin{cases}
			j^{-\alpha_0} \CAS j\text{ even}
			\\
			j^{-\alpha_1} \CAS j\text{ odd}
		\end{cases}
		\qquad\phi_j\DE j\frac\pi4
	\]
	Then the convergence exponent of $(b_n^{(2)})_{n=0}^\infty$ is $\frac 2{\alpha_0+\alpha_1}$ while the one of 
	$(b_n^{(3)})_{n=0}^\infty$ is $\frac 1{\alpha_0}$. The latter gives the correct value: $W_H$ is of order 
	$\frac 1{\alpha_0}$ with positive type. 
}

\section[{Upper bounds in terms of convergence exponents}]{Upper bounds in terms of convergence exponents}
\label{U133}

The upper estimates in this section feature convergence exponents as a measure for the growth of sequences associated with the lengths $l_j$ and angles $\phi_j$ of a Hamburger Hamiltonian. While the influence of the lengths on the bound is straightforward, the angles can only improve the bound if their increments are summable or the angles themselves converge. 

Recall the fact that a Hamburger Hamiltonian is in limit circle case, if and only if $\sum_{j=1}^\infty l_j<\infty$.

We first give a pointwise upper bound for $w_{H,22}(ir)$, i.e., the lower right entry of the monodromy matrix $W_H$ along the imaginary axis. Recall in this context that the growth along the imaginary axis of one entry of $W_H$ dominates $W_H$ itself in the sense of 
\Cref{U194} and \Cref{U197}. 

\begin{Theorem}
\label{U99}
\IndexxS{Theorem!monodromy matrix!upper bound (Hamburger Hamiltonian)}
	Let $l_j>0$ and $\phi_j\in\bb R$ with $\phi_{j+1}-\phi_j\not\equiv 0\mod\pi$, 
	and let $H$ be the Hamburger Hamiltonian with these lengths and angles. 
	Assume that $H$ is in limit circle case. 
	\begin{Enumerate}
	\item Suppose we have $\alpha, \beta \geq 1$ such that
\begin{align*}
(l_j)_{j=1}^\infty \in \ell^{\frac{1}{\alpha}}, \qquad (|\sin (\phi_{j+1}-\phi_j)|))_{j=1}^\infty \in \ell^{\frac{1}{\beta}}.
\end{align*}
Then (for, say, $r \geq 1$)
\begin{align}
\log |w_{H,22}(ir)| \lesssim r^{\frac{1}{\alpha+\beta}}.
\end{align}
	\item Suppose we have $\psi \in \bb R$ and $\alpha, \omega \geq 1$ such that
\begin{align}
\label{U155}
(l_j)_{j=1}^\infty \in \ell^{\frac{1}{\alpha}}, \qquad (l_j \sin^2 (\phi_j-\psi))_{j=1}^\infty \in \ell^{\frac{1}{\omega}}.
\end{align}
Then (for, say, $r \geq 1$)
\begin{align}
\label{U153}
\log |w_{H,22}(ir)| \lesssim r^{\frac{2}{\alpha+\omega}}.
\end{align}
	\end{Enumerate}
\end{Theorem}

Taking infima of possible upper bounds from \Cref{U99} we obtain explicit upper bounds for the order. The suprema defined below
can be written in terms of actual convergence exponents, but we feel that this would reduce readability. 

\begin{Corollary}
\label{U86}
	Let $l_j>0$ and $\phi_j\in\bb R$ with $\phi_{j+1}-\phi_j\not\equiv 0\mod\pi$, 
	and let $H$ be the Hamburger Hamiltonian with these lengths and angles. 
	Assume that $H$ is in limit circle case, let $\rho_H$ be the order of its monodromy matrix, and let
	\begin{align*}
		\alpha_0 &\DE \sup \Big\{\alpha>0 \, \big| \, \sum_{j=1}^{\infty} l_j^{\frac 1\alpha}<\infty \Big\} \geq 1.
	\end{align*}
	\begin{Enumerate}
	\item Assume that $\sum_{j=1}^{\infty} |\sin (\phi_{j+1}-\phi_j)| <\infty$, and set 
		\begin{align*}
			\beta_0 &\DE  \sup \Big\{\beta>0 \, \big| \, 
			\sum_{j=1}^{\infty} |\sin(\phi_{j+1}-\phi_j)|^{\frac 1\beta}<\infty \Big\}.
		\end{align*}
		Then
		\begin{align*}
			\rho_H \leq \frac{1}{\alpha_0+\beta_0}.
		\end{align*}
	\item Fix $\psi \in \bb R$ and set 
		\begin{align*}
			\omega_0 &\DE  \sup \Big\{\omega>0 \, \big| \, 
			\sum_{j=1}^{\infty} \big(l_j \sin^2(\phi_j-\psi)\big)^{\frac 1\omega}<\infty \Big\}.
		\end{align*}
		Then
		\begin{align*}
			\rho_H \leq \frac{2}{\alpha_0+\omega_0}.
		\end{align*}
	\end{Enumerate}
\end{Corollary}

\REMARKS{%
\item These results stem from \cite[Theorems~4.4, 4.7]{reiffenstein:kachamB-arXiv}, and their proofs are consequences of the
	Weyl coefficient approach. 
\item The first upper bound for the order in terms of a convergence exponent was given in \cite[Theorem~4.7]{berg.szwarc:2014}.
	It is formulated in the language of orthogonal polynomials as an upper estimate for $P(z) \DE \sqrt{\sum_{n=0}^\infty
	|p_n(z)|^2}$. It was rewritten to the Hamiltonian setting in \cite[Proposition 2.3]{pruckner.romanov.woracek:jaco}, and
	states that $(l_j)_{j=1}^{\infty} \in \ell^{p}$ implies $\rho_H \leq p$ with finite type (the angles are not
	taken into account). The proof is carried out by an elementary estimate of a product, noting that the monodromy matrix
	is an infinite product of linear polynomial factors. 
\item The upper bound from the previous item can be reobtained from item (ii) of \Cref{U99}, taking $\alpha \DE \omega \DE \frac{1}{p}$. In contrast to \cite{berg.szwarc:2014} we do not get an explicit type estimate from \Cref{U99}.
\item A weaker variant, which does not consider the angles either and which has an additional multiplicative factor $\log r$ on
	the right-hand side of the estimate, can be obtained from Romanov's Theorem I. This is shown in the discussion in
	\cite[\S4.3]{romanov:2017}. 
\item Usually one of the bounds from \Cref{U86} is correct if we can properly take into account the angles, but both are far 
	from the truth if we fail to do so (which is usually the case when the angles do not converge). For instance, 
	take $\delta_l\geq 1$ and consider the following situations.
\begin{Enumerate}
	\item Assume that
	\begin{align*}
			& l_1\DE 1,\qquad l_j\DE
			\begin{cases}
				\frac 1{j^{\delta_l}}\CAS \delta_l>1,
				\\
				\frac 1{j\log^2j}\CAS \delta_l=1,
			\end{cases}
			\quad\text{for }j\geq 2,
			\\
			& |\sin (\phi_{j+1}-\phi_j)| \asymp j^{-\delta_\phi} \quad\text{ with } \delta_\phi>1
			.
		\end{align*}
Applying the first item of \Cref{U86} yields $\rho_H \leq \frac{1}{\delta_l+\delta_\phi}$.	In fact we have $\rho_H = \frac{1}{\delta_l+\delta_\phi}$ due to \Cref{U208} (taking $k=2$). 
	\item With the same lengths as in (i), consider the angles
	\begin{align*}
			\phi_j\DE\sum_{k=1}^j\frac{(-1)^k}{k^{\delta_\phi}},\ j\in\bb N \quad \text{ with } \delta_\phi \leq 1
			.
		\end{align*}
		We first note that the first item of \Cref{U86} is not applicable. Nevertheless, due to $|\phi_{j+2}-\phi_j| \asymp j^{-\delta_\phi-1}$ the limit $\psi \DE \lim_{j \to \infty} \phi_j$ exists and the estimate $|\phi_j-\psi| \lesssim j^{-\delta_\phi}$ holds. Hence, in the second item of \Cref{U86} we have $\alpha_0 = \delta_l$, $\omega_0 = \delta_l+2\delta_\phi$. It follows that
	the order of $W_H$ is not larger than $\frac 1{\delta_l+\delta_\phi}$, so as in (i) we arrive at $\rho_H = \frac{1}{\delta_l+\delta_\phi}$.
	\item With the same lengths as in (i), consider the angles
	\begin{align*}
			\phi_j\DE\sum_{k=1}^j\frac{1}{k^{\delta_\phi}},\ j\in\bb N \quad \text{ with } \delta_\phi \leq 1
			.
		\end{align*}
Again the first item of \Cref{U86} is not applicable, while the second item is applicable but there is no convergence of angles. For any $\psi$ we have $\omega_0=\delta_l$ and thus only $\rho_H \leq 1/\delta_l$. But this upper bound is not optimal: By
	\Cref{U157}, which makes use of the fact that the lengths, angles, and angle-differences in this example are regularly
	varying, we have 
\begin{align*}
\rho_H= \begin{cases}
\frac{1}{\delta_l+\delta_\phi} & \text{if } \delta_l+\delta_\phi>2, \\
\frac{1-\delta_\phi}{\delta_l-\delta_\phi} & \text{if } \delta_l+\delta_\phi \leq 2.
\end{cases}
\end{align*}	
See also \Cref{U173}.
\end{Enumerate}
\item \Cref{U99} (i) and \Cref{U86} (i) ask us to assess the decay of the lengths and the decay of the angle-differences
	separately, while each of the lower bounds in \Cref{U31} needs only one sequence (of terms combining the lengths and
	angle-differences) as input. If both the lengths and the angle-differences are regularly varying and such that
	\Cref{U86} (i) applies, explicit computations enabled by Karamata's theorem show that the lower and upper bounds for the
	order of $W_H$ obtained from \Cref{U208} and \Cref{U86} (i) coincide.
	
%\item One situation where the upper bound from \Cref{U99} does coincide with the lower bound from \Cref{U97}, and hence determines 
%	the order of $W_H$, is that $\phi_n\in\{0,\frac\pi2\}$ and $\frac{l_{n+1}}{l_n}\asymp 1$. 
%	This is obviously a very particular case: a diagonal Hamburger Hamiltonian with non-oscillating lengths. However, when
%	translating to other languages, it gains significance: it corresponds to symmetric moment problems and Jacobi
%	matrices with vanishing diagonal. 
%	Knowledge about this particular class transfers to Stieltjes moment problems and Hamburger Hamiltonians with a
%	particular monotonicity property corresponding to Krein strings. This connection is established using a symmetrisation 
%	procedure (see \cite{berg:1995} for the context of moment problems and \Cref{U131} for the context of canonical
%	systems).
}

\section[{Upper bounds in terms of tails of convergent series}]{Upper bounds in terms of tails of convergent series}
\label{U160}

We present two upper bounds for the lower right entry $w_{H,22}$ of $W_H$ along the imaginary axis, which take into account the
decay of data in an averaged sense. The first result, \Cref{U158}, also assumes summability conditions as in \Cref{U133}.
Recall that the growth along the imaginary axis of one entry of $W_H$ dominates $W_H$ itself in the sense of 
\Cref{U194} and \Cref{U197}. 

\begin{theorem}
\label{U158}
\IndexxS{Theorem!monodromy matrix!upper bound (Hamburger Hamiltonian)}
	Let $l_j>0$ and $\phi_j\in\bb R$ with $\phi_{j+1}-\phi_j\not\equiv 0\mod\pi$, 
	and let $H$ be the Hamburger Hamiltonian with these lengths and angles. 
	Assume we have $\alpha \geq 1$ and $\beta \in (0,1)$ such that
	\begin{align*}
		(l_j)_{j=1}^\infty \in \ell^{\frac{1}{\alpha}}, \qquad 
		(|\sin (\phi_{j+1}-\phi_j)|))_{j=1}^\infty \in \ell^{\frac{1}{\beta}}.
	\end{align*} 
	Fix $\psi \in \mathbb{R}$, let $F\DF\bb N\to(0,\infty)$ be the increasing function 
	\[
		F(n)\DE n^{\frac{1-\beta}{\alpha}} \Big[
		\Big(\!\sum_{j=n+1}^\infty l_j \Big)\cdot\Big(\!\sum_{j=n+1}^\infty l_j \sin^2 (\phi_j-\psi) \Big)
		\Big]^{-\frac {\alpha+1}{2\alpha}}
		,
	\]
	and let $F^-\DF(0,\infty)\to\bb N$ be its left inverse 
	\[
		F^-(r)\DE\min\big\{n\in\bb N\DS F(n) \leq r\big\}
		.
	\]
Then $H$ is in limit circle case and (for, say, $r \geq 1$)
	\[
	\log |w_{H,22}(ir)| \lesssim \big(rF^-(r)^{1-\beta} \big)^{\frac{1}{\alpha+1}}
	.
	\]
\end{theorem}

\noindent
In the following variant the contributions of lengths and angles are fully separated. It requires as an
assumption that the telescoping sum of angle-differences converges absolutely. 

\begin{Theorem}
\label{U159}
	Let $l_j>0$ and $\phi_j\in\bb R$ with $\phi_{j+1}-\phi_j\not\equiv 0\mod\pi$, 
	and let $H$ be the Hamburger Hamiltonian with these lengths and angles. 
	Assume that $H$ is in limit circle case, and that $\sum_{j=0}^\infty|\sin(\phi_{j+1}-\phi_j)|<\infty$. 
	Let $G\DF\bb N\to(0,\infty)$ be the decreasing function 
	\[
		G(n)\DE\frac 1n
		\Big(\!\sum_{j=n+1}^\infty l_j\Big)^{\frac 12}
		\Big(\!\sum_{j=n+1}^\infty \big|\sin(\phi_{j+1}-\phi_j)\big|\,\Big)^{\frac 12}
		,
	\]
	and let $G^-\DF(0,\infty)\to\bb N$ be its left inverse 
	\[
		G^-(r)\DE\min\big\{n\in\bb N\DS G(n) \geq r\big\}
		.
	\]
	Then (for, say, $r\geq 2$)
	\[
		\log|w_{H,22}(ir)|\lesssim G^-\Big(\frac{\log r}{\sqrt r}\Big)\cdot\log r
		.
	\]
\end{Theorem}

\REMARKS{%
\item These results are taken from work of J.~Reiffenstein \cite{reiffenstein:kachamB-arXiv}. 
Their proofs rely on the Weyl coefficient approach, cf. \Cref{U110}, and direct estimates of the integral of $K_H$ over $[0,L]$. For $n \in \bb N$ the contributions of $[0,x_n]$ and of $[x_n,L]$ are estimated in different ways, and choosing $n=F^-(r)$ (in \Cref{U158}) or $n=G^-(r)$ (in \Cref{U159}) balances these contributions. 
\item Consider again the second part of the example from the remarks to \Cref{U133}: For $\delta_l\geq 1,\delta_\phi \geq 0$ set 
	\begin{align*}
		& l_1\DE 1,\qquad l_j\DE
		\begin{cases}
			\frac 1{j^{\delta_l}}\CAS \delta_l>1,
			\\
			\frac 1{j\log^2j}\CAS \delta_l=1,
		\end{cases}
		\quad\text{for }j\geq 2,
		\\
		& \phi_j\DE j^{1-\delta_\phi},\ j\in\bb N
		.
	\end{align*}
	For any $\epsilon>0$ we can apply \Cref{U158} with $\alpha \DE \delta_l-\epsilon$ and $\beta \DE \delta_\phi-\epsilon$,
	yielding
	\[
		\log|w_{H,22}(ir)|\lesssim r^{\frac{\delta_l-\delta_\phi+\epsilon}{\delta_l^2-\delta_\phi+\epsilon (2-\delta_l)}}.
	\]
	Since $\epsilon$ can be taken arbitrarily small, as upper bound for the order of $W_H$ we obtain 
	$\rho_H \leq(\delta_l-\delta_\phi)/(\delta_l^2-\delta_\phi)$. While this upper bound is better than the one we 
	obtained in the remarks to \Cref{U133}, it does not give the correct value of $\rho_H$ which is
	$(1-\delta_\phi)/(\delta_l-\delta_\phi)$, cf. \Cref{U157}.
\item If both the lengths $l_j$ and the angle-differences $|\sin (\phi_{j+1}-\phi_j)|$ are regularly varying, using Karamata's
	theorem we see that the lower and upper bounds for the order of $W_H$ obtained from \Cref{U208} and \Cref{U159}
	coincide. In this context \Cref{U159} is an improvement upon \Cref{U86} (i) since it provides a regularly varying upper
	bound for $\log|w_{H,22}(ir)|$ whose index is the same as the order of $W_H$. Yet an even better estimate is available
	(see the remarks to \Cref{U119}), which addresses regularly varying data specifically.
}

\section[{A Berezanskii-type theorem}]{A Berezanskii-type theorem}
\label{U135}

As a motivation we revisit Berezanskii's theorem discussed in \Cref{U107} from the Hamiltonian viewpoint. 
\Cref{U12} determines the order of the Nevanlinna matrix from the off-diagonal Jacobi parameters. 
There are three assumptions in that theorem. The first is that \IndexS{Carleman's condition}{Carleman's condition} is violated, which in terms of the
Hamiltonian parameters means that
\[
	\big(\sqrt{l_nl_{n+1}}\cdot|\sin(\phi_{n+1}-\phi_n)|\big)_{n=1}^\infty\in\ell^1.
\]
The other two assumptions, about regularity of off-diagonal and smallness of diagonal, can hardly be translated.
Let us turn to the conclusion of \Cref{U12}, which has a perfect translation: 
$(l_n)_{n=1}^\infty\in\ell^1$ and the order $\rho_H$ of the monodromy matrix is equal to the
convergence exponent of the sequence $\big([\sqrt{l_nl_{n+1}}\cdot|\sin(\phi_{n+1}-\phi_n)|]^{-1}\big)_{n=1}^\infty$.

The result we present now is of a similar kind. It gives conditions in terms of lengths and angles, sufficient for the
order of the monodromy matrix to be computed as a certain convergence exponent. 

\begin{Theorem}
\label{U141}
\IndexxS{Theorem!Berezanskii-type}
	Let $l_j>0$ and $\phi_j\in\bb R$ with $\phi_{j+1}-\phi_j\not\equiv 0\mod\pi$, 
	and let $H$ be the Hamburger Hamiltonian with these lengths and angles. 
	Assume that $H$ is in limit circle case, and, for simplicity of presentation, that 
	$\phi_n\not\equiv 0\mod\pi$ for all $n\in\bb N$. 
	Assume that the following three conditions hold.
	\begin{Itemize}
	\item ${\displaystyle \sum_{n=1}^\infty\big(l_n\sin^2\phi_n\big)^{\frac 12}\log n<\infty}$.
	\item The sequence $(l_n\sin^2\phi_n)_{n=1}^\infty$ is nonincreasing and $(|\cot\phi_n-\cot\phi_{n+1}|)_{n=1}^\infty$ 
		is nondecreasing.
	\item ${\displaystyle \sup_{n\geq 1}\big|\cot\phi_n-\cot\phi_{n+1}|<\infty}$.
	\end{Itemize}
	Then the order of the monodromy matrix is equal to the convergence exponent of the sequence 
	$\big([l_n\sin^2\phi_n]^{-1}\big)_{n=1}^\infty$. 
\end{Theorem}

\noindent
A typical case where the assumption on angle-differences is satisfied occurs when the angles perform a walk on the grid 
$\Arccot\bb Z$, since then $|\cot\phi_n-\cot\phi_{n+1}|=1$ for all $n\geq 1$. 

\REMARKS{%
\item \Cref{U141} is taken from \cite{pruckner.woracek:srt}. Its proof is altogether quite unusual. One uses formal
	symmetrisation to get a ``Hamburger Hamiltonian with positive and negative lengths''. In turn, one passes to a Hamburger
	Hamiltonian that is possibly in limit point case, and applies the result \cite[Theorem~1]{kac:1990} of I.S.~Kac about 
	Stieltjes strings. 
\item \Cref{U141} and \Cref{U12} have some overlap. The core of Berezanskii's theorem is the case that $a_n=0$ (adding a ``small
	diagonal'' can be understood as a small perturbation). 
	For vanishing diagonal the translation from Jacobi parameters to Hamiltonian
	parameters is simple, and for orders $<\frac 12$ (to be precise, knowing that 
	$\sum_{n=1}^\infty \sqrt{l_n}\log n<\infty$) 
	the assertion of Berezanskii's theorem (in its original variant) follows from \Cref{U141}.
\item In general \Cref{U141} and Berezanskii's theorem (in the strong variant from \cite{reiffenstein:kachamB-arXiv} stated in 
	\Cref{U12}) are incomparable. \Cref{U141} is weaker in the sense that it is bound to orders $\leq\frac 12$ and requires
	limit circle case as an a priori assumption. On the other hand, it is stronger in the sense that ``smallness and regularity 
	of the diagonal'' is not necessarily present: in \cite[Example~4.8]{pruckner.woracek:srt} an example is given where
	\Cref{U141} applies, while the corresponding Jacobi parameters satisfy $\lim_{n\to\infty}\frac{b_{n+1}}{b_n}=1$ and 
	\[
		\liminf_{n\to\infty}\frac{|a_n|}{b_n}=0,\quad \limsup_{n\to\infty}\frac{|a_n|}{b_n}=2
		.
	\]
}

\section[{Exploiting regular variation I: upper bound from majorants}]{Exploiting regular variation I: 
	upper bound from majorants}
\label{U119}

To improve computability one can work with regularly varying majorants instead of the lengths and angles themselves.
By doing so, it becomes possible to treat even irregularly behaving lengths and angles, potentially at the cost of 
some loss of precision.

\begin{theorem}
\label{U45}
\IndexxS{Theorem!Hamburger Hamiltonian!regular variation}
\IndexxS{Theorem!monodromy matrix!upper bound (Hamburger Hamiltonian)}
	Let $l_j>0$ and $\phi_j\in\bb R$ with $\phi_{j+1}-\phi_j\not\equiv 0\mod\pi$, 
	and let $H$ be the Hamburger Hamiltonian with these lengths and angles. 
	Assume that $H$ is in limit circle case. 

	Let $\psi\in\bb R$ and let $\ms d_l,\ms d_\phi,\ms c_l,\ms c_\phi$ be regularly varying functions that are $\sim$ 
	to some nonincreasing functions, such that $\ms d_l \asymp 1 \asymp \ms d_\phi$ locally, 
	$\ms c_\phi(t)\lesssim\ms c_l(t)$ for sufficiently large $t$, $\Ind\ms d_l+\Ind\ms d_\phi<0$, and 
	$\lim_{t\to\infty}(\ms c_l\ms c_\phi)(t)=0$. Assume that
	\begin{align*}
		l_j &\lesssim\ms d_l(j),& |\sin(\phi_{j+1}-\phi_j)|&\lesssim\ms d_\phi(j)
		,
		\\
		\sum_{j=N+1}^\infty l_j&\lesssim\ms c_l(N),& \sum_{j=N+1}^\infty l_j\sin^2(\phi_j-\psi)&\lesssim\ms c_\phi(N)
		.
	\end{align*}
	Denote
	\[
		\IndexN{\ms D(t)}\DE \frac{1}{(\ms d_l\ms d_\phi)(t)},\ \delta\DE\Ind\ms D,\qquad
		\IndexN{\ms C(t)}\DE \frac{1}{(\ms c_l\ms c_\phi)^{\frac 12}(t)},\ \gamma\DE\Ind\ms C
		,
	\]
	and set 
	\begin{align*}
		& \delta_l \DE -\Ind \ms d_l, && \delta\DE\Ind\ms D = \delta_l+\delta_\phi, 
		\\
		& \delta_\phi \DE -\Ind \ms d_\phi, && \gamma\DE\Ind\ms C.
	\end{align*}
	Let $\IndexN{\ms k(r)}$ be an asymptotic inverse of $\ms D$, and set 
	\[
		\IndexN{\ms h(r)}\DE
		\sup\Big\{t\in[1,\infty)\DS \sup_{1\leq s\leq t}\frac{\ms d_\phi(s)}{r\ms d_l(s)}\leq 1\Big\}
		.
	\]
	Then we have the following bounds for the monodromy matrix $W_H$ and its order $\rho_H$. 
	\\[3mm]
	\scalebox{0.8}{
	\[
		\begin{array}{lll|l|l}
			\multicolumn{3}{l|}{\text{\rm Data satisfies}} 
			& \max_{|z|=r}\log\|W_H(z)\|\text{\rm\ is }\lesssim
			&\rho_H \leq
			\raisebox{-6pt}{\rule{0pt}{17pt}}
			\\
			\hline
			\hline
			\multicolumn{3}{l|}{\ms D(t)\lesssim t\ms C(t)}
			& \frac r{\ms C(\ms f^-(r))}\quad\text{\rm where}
			\raisebox{-9pt}{\rule{0pt}{25pt}}
			\\
			& 
			& 
			& \ms f(t)\DE t\ms C(t)\log\big[\alpha\frac{t\ms C(t)}{\ms D(t)}\big]
			\raisebox{-6pt}{\rule{0pt}{1pt}}
			& \raisebox{-2pt}{\mbox{\LARGE\(\frac{1}{1+\gamma}\)}}
			\\
			& 
			& 
			& \alpha\DE 4\sup_{t\geq 1}\frac{\ms D(t)}{t\ms C(t)}
			\raisebox{-13pt}{\rule{0pt}{1pt}}
			\\
			\hline
			t\ms C(t)\lesssim\ms D(t),
			& \int\limits_1^\infty\ms D(s)^{-\frac 12}\DD s<\infty,
			& 
			& \raisebox{-10pt}{$r^{\frac 12}\int\limits_{\ms k(r)}^\infty\ms D(s)^{-\frac 12}\DD s$}
			\raisebox{-13pt}{\rule{0pt}{30pt}}
			& \raisebox{-12pt}{\mbox{\LARGE\(\frac{1}{\delta}\)}}
			\\[-8pt]
			& \multicolumn{2}{l|}{\parbox{55mm}{\rm ($\gamma>0$ or $\frac{\ms d_\phi}{\ms d_l}$ 
			$\approx$ to nondecreasing)}}
			& 
			\raisebox{-10pt}{\rule{0pt}{1pt}}
			\\
			\cline{2-5}
			& \int\limits_1^\infty\ms D(s)^{-\frac 12}\DD s=\infty,
			& \frac 1{\ms d_l(t)}\lesssim t\ms C(t),
			& \frac r{\ms C(\ms f_1^-(r))}\quad\text{\rm where}
			\raisebox{-13pt}{\rule{0pt}{30pt}}
			& \raisebox{-10pt}{\mbox{\large\(\frac{2-\delta+\gamma}{2-\delta+2\gamma}\)}}
			\\
			& 
			& \parbox{24mm}{\rm $(\delta,\gamma)\neq(2,0)$}
			& \ms f_1(t)\DE[\ms C(t)\int_1^t\ms D(s)^{-\frac 12}\DD s]^2
			\raisebox{-10pt}{\rule{0pt}{1pt}}
			\\
			\cline{3-5}
			& 
			& t\ms C(t)\lesssim\frac 1{\ms d_l(t)},
			& \raisebox{-6pt}{$r^{\frac 12}\int\limits_1^{\ms h(r)}\ms D(s)^{-\frac 12}\DD s$}
			& \raisebox{-6pt}{\mbox{\Large\(\frac{1-\delta_\phi}{\delta_l-\delta_\phi}\)}}
			\raisebox{-12pt}{\rule{0pt}{31pt}}
			\\[-12pt]
			& 
			& \parbox{25.5mm}{\rm $\int\limits_1^\infty\ms d_l(s)\DD s<\infty$, \\
			$\delta_l\!>\!\delta_\phi$}
			& 
			\raisebox{-2pt}{$\quad \,\,\, +r\mkern-5mu\int\limits_{\ms h(r)}^\infty\mkern-5mu\ms d_l(s)\DD s$}
			& 
			\\[7mm]
			\hline
		\end{array}
	\]
	}
\end{theorem}

\REMARKS{%
\item This theorem is taken from the recent work \cite{pruckner.reiffenstein.woracek:sinqB-arXiv} of R.~Pruckner,
	J.~Reiffenstein and H.~Woracek, where it is deduced from a general (much more cumbersome) result, namely
	\cite[Theorem~2.2]{pruckner.reiffenstein.woracek:sinqB-arXiv}.
\item \Cref{U45} still allows small oscillations of lengths and angle-differences: a function that is $\asymp$ to a
	regularly varying function need not necessarily be itself regularly varying.
\item The main assumptions on the data $\ms d_l,\ms d_\phi,\ms c_l,\ms c_\phi$ made in \Cref{U45} are that those functions are
	regularly varying and that $\Ind\ms d_l+\Ind\ms d_\phi<0$. 
	The monotonicity assumptions are only minor restrictions: for example they are automatically fulfilled unless the 
	function under consideration is slowly varying. The assumption that $\ms c_\phi\lesssim\ms c_l$ is no loss in
	generality, since replacing $\ms c_\phi$ by $\min\{\ms c_\phi,\ms c_l\}$ does not affect validity of any of the other
	assumptions. 
\item Intuitively the function $\ms h(r)$ is just the inverse to the quotient $\frac{\ms d_\phi}{\ms d_l}$ 
	(which is of course formally not correct). In fact, the situation is as follows:
	If the quotient $\frac{\ms d_\phi}{\ms d_l}$ is bounded, we have eventually $\ms h(r)=\infty$. 
	If $\Ind\frac{\ms d_\phi}{\ms d_l}>0$, then $\ms h(r)$ is an asymptotic inverse of $\frac{\ms d_\phi}{\ms d_l}$. 
	However, if $\Ind\frac{\ms d_\phi}{\ms d_l}=0$, oscillations are possible, and then it is necessary to define 
	$\ms h(r)$ as the written supremum.
\item An illustrative example, where we discuss in detail the situation that lengths and angles are majorised by
	certain concrete and simple regularly varying functions, is given in \Cref{U149}.
\item \Cref{U45} has an overlap with \Cref{U159}: 
	Assume that we are given only regularly varying functions $\ms d_l,\ms d_\phi$, such that:
	$\ms d_l,\ms d_\phi$ are locally $\asymp 1$, satisfy $\Ind\ms d_l, \Ind \ms d_\phi<-1$, 
	and majorise our data as
	\[
		l_j\lesssim\ms d_l(j),\qquad |\sin(\phi_{j+1}-\phi_j)|\lesssim\ms d_\phi(j)
		.
	\]
	Then we can use the trivial choices $\ms c_l(t)\DE t\ms d_l(t)$ and $\ms c_\phi(t)\DE t^3\ms d_l(t)\ms d_\phi(t)^2$
	in \Cref{U45}. The upper bound for $\max_{|z|=r}\log\|W_H(z)\|$ obtained in this way is $\ms k(r)$, while
	\Cref{U159} gives the upper bound $\ms k (r/(\log r)^2) \log r$. Since $\Ind \ms k<\frac 12$, the bound from
	\Cref{U45} is better, but only by a slowly varying factor.
\item A systematic treatment of the case where we are only given functions $\ms d_l,\ms d_\phi$ is presented in \Cref{U134}.
}

\section[{Exploiting regular variation II: lower bound from minorants}]{Exploiting regular variation II: 
	lower bound from minorants}
\label{U150}

We present two theorems and a corollary. The first theorem is generic, while the second deals with a particular situation. Recall that $w_{H,22}$ denotes the lower right entry of $W_H$.

\begin{Theorem}
\label{U95}
\IndexxS{Theorem!Hamburger Hamiltonian!regular variation}
\IndexxS{Theorem!monodromy matrix!lower bound (Hamburger Hamiltonian)}
	Let $l_j>0$ and $\phi_j\in\bb R$ with $\phi_{j+1}-\phi_j\not\equiv 0\mod\pi$, 
	and let $H$ be the Hamburger Hamiltonian with these lengths and angles.
	Assume that $H$ is in limit circle case. 
	Let $\ms f,\ms g$ be regularly varying functions that are asymptotic inverses of each other, 
	and fix $k\in\bb N$, $k\geq 2$. Then 
	\[
		\Big[\frac 12\sum_{i,j=n+1}^{n+k} l_il_j \sin^2(\phi_i-\phi_j)\Big]^{\frac 12}\gtrsim \frac{1}{\ms f(n)} 
		\ \Longrightarrow\ \log|w_{H,22}(ir)|\gtrsim\ms g(r)
		.
	\]
\end{Theorem}

\noindent
Taking $k=2$ and using notation similar to \Cref{U119} we obtain as a corollary:

\begin{Corollary}
\label{U47}
	Let $l_j>0$ and $\phi_j\in\bb R$ with $\phi_{j+1}-\phi_j\not\equiv 0\mod\pi$, 
	and let $H$ be the Hamburger Hamiltonian with these lengths and angles. 
	Assume that $H$ is in limit circle case. 
	Let $\ms d_l$ and $\ms d_\phi$ be regularly varying and such that
	\[
		l_n\gtrsim\ms d_l(n),\qquad \big|\sin(\phi_{n+1}-\phi_n)\big|\gtrsim\ms d_\phi(n),
	\]
	and let	$\IndexN{\ms k(r)}$ be an asymptotic inverse of $\IndexN{\ms D(t)}\DE \frac 1{(\ms d_l\ms d_\phi)(t)}$. Then (for, say, $r\geq 1$)
	\[
		\log|w_{H,22}(ir)|\gtrsim\ms k(r)
		.
	\]
\end{Corollary}

\begin{Theorem}
\label{U148}
\IndexxS{Theorem!Hamburger Hamiltonian!monotone angles}
\IndexxS{Theorem!monodromy matrix!lower bound (Hamburger Hamiltonian)}
	Let $l_j>0$ and $\phi_j\in\bb R$ with $\phi_{j+1}-\phi_j\not\equiv 0\mod\pi$, 
	and let $H$ be the Hamburger Hamiltonian with these lengths and angles. 
	Assume that $H$ is in limit circle case. 
	Let $\ms d_l$ and $\ms d_\phi$ be regularly varying with $-1<\Ind\ms d_\phi<0$, such that
	\[
		l_n\gtrsim\ms d_l(n),\qquad \big|\phi_{n+1}-\phi_n\big|\asymp\ms d_\phi(n),
	\]
	Further assume that $(\phi_j)_{j=0}^\infty$ is eventually monotone. Then, letting $\IndexN{\ms h(r)}$ be an asymptotic inverse of $\frac{\ms d_\phi(t)}{\ms d_l(t)}$, we have (for, say, $r\geq 1$)
	\[
		\log|w_{H,22}(ir)|\gtrsim r\int_{\ms h(r)}^\infty\ms d_l(t)\DD t
		.
	\]
\end{Theorem}

\REMARKS{%
\item \Cref{U95} is a consequence of \Cref{U14}\,(ii), and Karamata's theorem. \Cref{U148} is 
	\cite[Theorem~1.1]{reiffenstein:kachamA-arXiv}. Both results are not elementary and use
	the Weyl coefficient approach. 
\item \Cref{U47}, although a corollary of \Cref{U95} for $k=2$, also has an elementary proof, for which we refer to
	\cite[Corollary~2.5]{pruckner.woracek:sinqA}.
\item We will see in \Cref{U134} that the lower bounds in \Cref{U47,U148} are sharp. Interestingly, the generic and elementary
	bound from \Cref{U47} is attained quite often. 
\item \Cref{U148} deals with the extreme situation that angles are eventually monotone. A sharp lower bound that takes into
	account oscillations of the angles is at present not available. 
}

%**************************************************************************
%***                            Last Change: Mon 31 Mar 2025 11:28
%***   < PART V >
%***
%**************************************************************************

\clearpage
\PART{V}{Additions and examples}
\label{U128}

%%%%%%%%%%%%%%%%%%%%%%%%%%%%%%%%%%%%%%%%%%%%%%%%%%%%%%%%%%%%%%%%%%%%%%%%%%%%%%%%%%%%%%%%%

\Intro{%
	In this part we collect some selected theorems and examples that supplement the general theory and illustrate phenomena
	related to growth and density. In the first two sections we discuss Hamiltonians with continuous rotation angle. First
	we present some general results, then a concrete (and enlightening) example. In the second pair of sections we turn to
	Hamburger Hamiltonians. Again we first give some general results and then a concrete example. 

	In \Cref{U105} we revisit the setting of Jacobi matrices and present a recent result which is important, 
	yet somewhat peculiar. There for one very particular Hamiltonian not only finiteness of
	$\limsup_{r\to\infty}\frac{n_H(r)}{r^\rho}$ (where $\rho$ is the order) is shown, but the actual value of this limit
	superior is determined. This theorem is out of the realm of the general results presented in Part~II, since those
	determine for intrinsic reasons the limit superior only up to universal multiplicative constants, and thus cannot
	possibly determine its exact value.

	The next two sections contain one perturbation result each, and in the pair of sections following those we discuss 
	symmetry and Krein strings. Exploiting symmetry is a practical tool. The work on strings of I.S.~Kac and M.G.~Krein 
	has played an important role in the development of many results about canonical systems and cannot be appreciated enough. 

	Last but not least, in \Cref{U273}, we touch upon the inverse problem (where we focus on the limit circle case): 
	can one construct Hamiltonians whose monodromy matrix has prescribed growth? 
}
\begin{flushleft}
	\S.\,\ref{U106}\ Hamiltonians with continuous rotation angle\ \dotfill\quad\pageref{U106}
	\\[1mm]
	\S.\,\ref{U122}\ An illustrative example: the chirp function\ \dotfill\quad\pageref{U122}
	\\[1mm]
	\S.\,\ref{U134}\ Hamburger Hamiltonians with regularly varying data\ \dotfill\quad\pageref{U134}
	\\[1mm]
	\S.\,\ref{U149}\ Hamburger Hamiltonians with power-log decay\ \dotfill\quad\pageref{U149}
	\\[1mm]
	\S.\,\ref{U105}\ Bochkov's theorem about type\ \dotfill\quad\pageref{U105}
	\\[1mm]
	\S.\,\ref{U176}\ Modifying the rotation of the Hamiltonian\ \dotfill\quad\pageref{U176}
	\\[1mm]
	\S.\,\ref{U121}\ Cutting out pieces of the Hamiltonian\ \dotfill\quad\pageref{U121}
	\\[1mm]
	\S.\,\ref{U131}\ Exploiting symmetry\ \dotfill\quad\pageref{U131}
	\\[1mm]
	\S.\,\ref{U124}\ Krein strings\ \dotfill\quad\pageref{U124}
	\\[1mm]
	\S.\,\ref{U273}\ Inverse results\ \dotfill\quad\pageref{U273}
	\\[1mm]
\end{flushleft}
\makeatother
\renewcommand{\thesection}{\arabic{section}}
\renewcommand{\thelemma}{\arabic{section}.\arabic{lemma}}
\makeatletter
\clearpage

%%%%%%%%%%%%%%%%%%%%%%%%%%%%%%%%%%%%%%%%%%%%%%%%%%%%%%%%%%%%%%%%%%%%%%%%%%%%%%%%%%%%%%%%%

%
%
%
\section[{Hamiltonians with continuous rotation angle}]{Hamiltonians with continuous rotation angle}
\label{U106}

For continuously rotating Hamiltonians one can show a growth estimate in terms of the quality of continuity of their rotation 
(measured, e.g., by the H\"older exponent). The intuition is that smoothly rotating Hamiltonians have slowly growing fundamental 
matrices. This intuition is legitimate only up to Lipschitz continuity (so one must not think of differentiability or
similar) and for orders between $\frac 12$ and $1$. It breaks down
at growth of speed $r^{\frac 12}$. 

\begin{Theorem}
\label{U18}
\IndexxS{Theorem!monodromy matrix!continuous angle}
	Let $H\in\bb H_{a,b}$ be in limit circle case, and assume that $H$ can be written in the form 
	\begin{equation}
	\label{U23}
		H(t)=\Tr H(t)\cdot\xi_{\phi(t)}\xi_{\phi(t)}^T,\qquad t\in(a,b)\text{ a.e.}
		,
	\end{equation}
	with a continuous function $\phi\DF[a,b]\to\bb R$. 
	Let 
	\[
		\omega(\delta)\DE\sup\big\{|\phi(t)-\phi(s)|\DS t,s\in[a,b],|t-s|\leq\delta\big\}
	\]
	be the \IndexS{modulus of continuity}{modulus of continuity} of $\phi$, let $\Omega$ be the inverse function of $\delta\mapsto\delta\omega(\delta)$, and 
	denote 
	\[
		\Gamma(r)\DE\frac 1{\Omega(\frac 1r)}
		.
	\]
	Moreover, set $l\DE b-a$ and $L\DE\int_a^b\Tr H(t)\DD t$. Then
	\begin{equation}
	\label{U20}
		\log \Big(\max_{|z|=r}\|W_H(z)\|\Big)\leq 3l\cdot\Gamma\Big(\frac Llr\Big)+\BigO\big(\log r\big)
		.
	\end{equation}
\end{Theorem}

\noindent
Observe that $\Gamma$ depends monotonically on $\omega$. This leads to the following corollary.

\begin{Corollary}
\label{U19}
	Consider a Hamiltonian $H$ of the form \cref{U23}, and let $\kappa\in(0,1]$. If $\phi$ is \IndexS{H\"older continuous}{H\"older continuous} with 
	exponent $\kappa$, then 
	\[
		\log \Big(\max_{|z|=r}\|W_H(z)\|\Big)\lesssim r^{\frac 1{1+\kappa}}
		.
	\]
\end{Corollary}

\noindent
As the following example shows the estimate in this corollary is sharp. 

\begin{Example}
\label{U184}
	Let $\alpha\in(0,1)$ and $\beta>\frac 1\alpha$. The Weierstra{\ss} function with parameters $\alpha,\beta$ is 
	\[
		\varphi(t) \DE \sum_{n=0}^\infty \alpha^n\cos(\beta^n\pi t), \qquad t\in\bb R.
	\]
	It is H\"older continuous with exponent $\kappa\DE-\frac{\log\alpha}{\log\beta}$. 

	The growth of the monodromy matrix corresponding to the Hamiltonian 
	$H(t)\DE\xi_\varphi(t)\xi_\varphi(t)^T$, $t\in(0,1)$, can be determined 
	when $\beta$ is an even integer with $\beta>\frac{1+\pi/2}\alpha$ (this assumption is there for technical reasons; 
	it is probably not necessary). It is 
	\[
		\log \Big(\max_{|z|=r}\|W_H(z)\|\Big)\asymp r^{\frac 1{1+\kappa}}
		.
	\]
	The modulus of continuity of $\varphi$ satisfies $\omega(\delta)\asymp\delta^\kappa$, which yields that also 
	$\Gamma(r)\asymp r^{\frac 1{1+\kappa}}$.  
\end{Example}

\noindent
Since $\omega$ is the modulus of continuity of a continuous function, we have $\Gamma(r)\gtrsim r^{\frac 12}$,
and hence the bound from \Cref{U18} can never go below $r^{\frac 12}$. 
The next result emphasises that order $\frac 12$ is indeed an intrinisic threshold. 

\begin{Proposition}
\label{U272}
	Let $H$ be as in \cref{U23}.
	\begin{Enumerate}
	\item Assume that $\phi$ is Lipschitz continuous and not constant. Then 
		\[
			\log \Big(\max_{|z|=r}\|W_H(z)\|\Big)\asymp r^{\frac 12}
			.
		\]
	\item Assume that $\phi$ is monotone and bounded. If the order of $W_H$ is less than $\frac 12$, then $\phi'=0$ a.e.
	\end{Enumerate}
\end{Proposition}

\REMARKS{%
\item \Cref{U18} and \Cref{U19} are established in \cite[\S5]{pruckner.woracek:sinqA} as an application of Romanov's Theorem~I
	(cf.\ \Cref{U43}), and \Cref{U184} is \cite[Theorem~7.4]{langer.reiffenstein.woracek:kacest-arXiv}. 
	Item (i) of \Cref{U272} is \cite[Remark~6.7]{langer.reiffenstein.woracek:kacest-arXiv}, and item (ii) is essentially 
	\cite[Theorem~1.6(c)]{remling.scarbrough:2020}.
\item \Cref{U19} improves upon \cite[Corollary~4(1)]{romanov:2017} where the bound $r^{1-\frac\kappa 2}$ is given. 
\item Also on a finer scale than powers the bound from \Cref{U18} is sharp, at least up to a small gap (and we expect that this
	gap can be removed). 
	Examples are constructed in \cite[\S6]{pruckner.woracek:sinqA}: Let $\ms g$ be a regularly varying function with 
	$\Ind\ms g\in(\frac 12,1)$, let $\ms m$ be regularly varying such that $\int_1^\infty\frac 1{\ms m(t)}\DD t<\infty$, and
	let $\ms n$ be an asymptotic inverse of $\ms m$. 
	Then there exists $H$ as in \Cref{U18}, such that $\Gamma\asymp\ms g$ and 
	\[
		(\ms n\circ\ms g)(r)\lesssim\log \Big(\max_{|z|=r}\|W_H(z)\|\Big)\lesssim\ms g(r)
		.
	\]
	Choosing for example $\ms m(r)\DE r(\log r)^{1+\epsilon}$ with some $\epsilon>0$, then 
	$(\ms n\circ\ms g)(r)\asymp\frac{\ms g(r)}{(\log r)^{1+\epsilon}}$. 
\item The bound from \Cref{U18} may fail badly: there exist Hamiltonians with continuous rotation angle whose monodromy matrix 
	has arbitrary small order. Examples are constructed e.g., in \cite[Example~5.5]{pruckner.woracek:sinqA} or 
	\cite[Theorem~1.6(b)]{remling.scarbrough:2020}.

	An example where the order of $W_H$ is in $(\frac 12,1)$, but is still strictly smaller than what is given by
	\Cref{U19}, is presented in \Cref{U142}. See the remarks to \Cref{U122}. 
\item A discrete variant of \Cref{U18} is given in \cite{reiffenstein:kachamB-arXiv}. It reads as follows:
	Let $l_j>0$ and $\phi_j\in\bb R$, let $H$ be the Hamburger Hamiltonian with these lengths and angles, 
	cf. \Cref{U286}, and recall the notation $x_n\DE\sum_{j=1}^n l_j$. Assume that $H$ is in limit circle case.
	Moreover, let $\kappa\in(0,1]$. If the angles satisfy the discrete H\"older condition 
	\[
		|\phi_{m+1}-\phi_n|\lesssim|x_m-x_n|^\kappa\qquad\text{for }n,m\in\bb N,
	\]
	then 
	\[
		\log \Big(\max_{|z|=r}\|W_H(z)\|\Big)\lesssim r^{\frac 1{1+\kappa}}
		.
	\]
\item At present no work has been undertaken to deal with continuously rotating Hamiltonians in limit point case.
}

\section[{An illustrative example: the chirp function}]{An illustrative example: the chirp function}
\label{U122}

We give an example of Hamiltonians with H\"older continuous and possibly oscillating rotation angle, for which 
the formula \cref{U140} for the growth of $\log|w_{H,22}(ir)|$ can be evaluated (recall here the connection \cref{U77}).

\begin{theorem}
\label{U142}
\IndexxS{Theorem!chirp function}
	Let $(\gamma,\beta)\in\bb R^2 \setminus (0,0)$, let $\varphi\DF(0,1]\to\bb R$ be the chirp signal 
	\[
		\varphi(t)\DE t^\gamma\sin\Big(\frac 1{t^\beta}\Big)
		,
	\]
	and consider the Hamiltonian $H(t)\DE\xi_{\varphi(t)}\xi_{\varphi(t)}^T$, $t\in(0,1)$. Then 
	\begin{equation}
	\label{U167}
		\log\Big(\max_{|z|=r}\|W_H(z)\|\Big)\asymp 
		\begin{cases}
			r^{\frac{\beta}{\beta+\gamma+1}} \CAS \beta>\gamma+1\wedge\gamma\geq 0,
			\\
			r^{\frac{\beta-\gamma}{\beta-\gamma+1}} \CAS \beta>\gamma+1\wedge\gamma <0,
			\\
			r^{\frac 12}\log r \CAS \beta=\gamma+1,
			\\
			r^{\frac 12} \CAS \beta<\gamma+1.
		\end{cases}
	\end{equation}
\end{theorem}

\medskip
\begin{center}
	\begin{tikzpicture}[x=1.2pt,y=1.2pt,scale=0.8,font=\fontsize{8}{8}]
		\draw[dashed,->] (10,0)--(234,0);
		\draw[dashed] (130,-100)--(130,20);
		\draw[->] (130,20)--(130,122);
		\draw     (20,-90)--(130,20);
		\draw (130,20) -- (220,110) node[midway, sloped, below]{ ${\displaystyle r^{\frac 12}\log r}$};
		
		\draw (122,120) node {${\displaystyle \beta}$};
		\draw (234,-8) node {${\displaystyle \gamma}$};
		
		\draw (160,-30) node {\large $r^{\frac 12}$};
		\fill[pattern=dots, opacity=0.1] (20,-90)--(220,110)--(220,-90)--(20,-90);
		
		\draw (80,40) node {\Large $r^{\frac{\beta-\gamma}{\beta-\gamma+1}}$};
		\fill[pattern=crosshatch dots, opacity=0.25] (20,-90)--(130,20)--(130,110)--(20,110)--(20,-90);
		
		\draw (155,80) node {\Large $r^{\frac{\beta}{\beta+\gamma+1}}$};
		\fill[pattern=crosshatch, opacity=0.15] (130,20)--(130,110)--(220,110)--(130,20);
	\end{tikzpicture}
\end{center}

\REMARKS{%
\item This result is \cite[Theorem~6.9]{langer.reiffenstein.woracek:kacest-arXiv}.
	To prove \Cref{U142} one estimates $\kappa_H(r)$ (see \Cref{U143}) from below with some expression, and 
	estimates the integral on the right side of \cref{U140} with the same expression from above. 
	Applying \Cref{U144} finishes the argument.
\item For the Hamiltonians $H$ as in the theorem we have $\log|w_{H,22}(ir)|\asymp\kappa_H(r)$. 
	Thus they provide examples where the lower bound in \Cref{U144} is attained. 
\item In this example the speed of growth of $\log|w_{H,22}(ir)|$ is determined by the behaviour of $H$ at a single point
	(namely, the left endpoint $0$). In fact, let $c\in(0,1)$. Since $\varphi$ is Lipschitz continuous on $(c,1)$, the monodromy matrix of the restriction of $H$ to $(c,1)$ satisfies
	\[
		\log\big|w_{H|_{(c,1)},22}(ir)\big|\asymp r^{\frac 12}
		,
	\]
	cf.\ \Cref{U272}. Now note that the exponent appearing in the first and second case on the right side of \cref{U167} 
	is $>\frac 12$. 
\item It is interesting to analyse the dependency on $(\beta,\gamma)$ of the exponent appearing on the right side of
	\cref{U167}. 
	\begin{Itemize}
	\item If the frequency is small relative to decay/increase, we are in the last case of \cref{U167}, and the growth is
		as slow as it possible can be (recall the previous remark).
	\item If $\varphi(t)$ oscillates with relatively high frequency, the growth ranges over all possible powers $r^\alpha$
		with $\alpha\in(\frac 12,1)$. 
	\item We have to distinguish whether oscillations are damped or amplified. In the first case the exponent $\alpha$
		depends only on the relation between frequency and damping, in the second case frequency and amplification act
		as independend parameters. 
	\end{Itemize}
\item \Cref{U142} illustrates a limitation of \Cref{U19}. 
	Consider parameters $\beta\geq\gamma>0$. In this case $\varphi(t)$ is \IndexS{H\"older continuous}{H\"older continuous} with H\"older exponent 
	$\frac{\gamma}{\beta+1}$. \Cref{U19} implies $\rho_H\leq(1+\frac\gamma{\beta+1})^{-1}=\frac{\beta+1}{\beta+\gamma+1}$, 
	which is not the correct value. With a reparametrisation trick one can improve this to 
	$\rho_H\leq\frac{\beta}{\beta+\gamma}$, cf.\ \cite[Example 5.6]{pruckner.woracek:sinqA}. 
	However, this is still not the correct value for $\rho_H$ (unless $\beta=\gamma$).
}

\section[{Hamburger Hamiltonians with regularly varying data}]{Hamburger Hamiltonians with
	regularly varying data}
\label{U134}

In this section we consider Hamburger Hamiltonians in limit circle case whose lengths $l_j$ and angle-differences 
$|\sin(\phi_{j+1}-\phi_j)|$ are bounded by, or even comparable to, functions $\ms d_l$ and $\ms d_\phi$ that satisfy:
\begin{Itemize}
\item $\ms d_l$ and $\ms d_\phi$ are regularly varying, 
\item $\ms d_l$ and $\ms d_\phi$ are locally $\asymp 1$,
\item $\ms d_\phi$ is $\sim$ to some nonincreasing function. 
\end{Itemize}
The second and third requirements are made for technical reasons; the second is no restriction, and the third is a minor
restriction. Denote
\[
	\delta_l\DE-\Ind\ms d_l,\qquad\delta_\phi\DE-\Ind\ms d_\phi,
\]
then $\delta_l\geq 1$ with $\int_1^\infty\ms d_l(t)\DD t<\infty$ since the lengths are summable, and 
$\delta_\phi\geq 0$ since the sine of angle-differences is bounded. 

For data with 
\begin{equation}
\label{U154}
	l_j\asymp\ms d_l(n),\qquad \big|\sin(\phi_{j+1}-\phi_j)\big|\asymp\ms d_\phi(n)
	,
\end{equation}
a fairly complete picture can be given, nicely illustrating the way how lengths and angles 
influence the growth of $W_H$. However, it must be emphasised that \cref{U154} does not contain enough information 
to always determine the growth (e.g.\ the order) of the monodromy matrix.

One may say that there are two regions of essentially different behaviour, one corresponding to slow growth 
and the other to faster growth. The border line is the threshold ``$\rho_H=\frac 12$'' (where $\rho_H$ is the
order of $W_H$). 
Finally, there is one exceptional point where not much can be said.
\begin{itemize}
\item[\textsf{(I)}] Region of simple behaviour (small orders): $\delta_l+\delta_\phi>2$.
\item[\textsf{(II)}] \IndexS{Critical triangle}{critical triangle} (large orders): $\delta_l+\delta_\phi<2$.
\item[\textsf{(III)}] Switching line ($\rho_H=\frac 12$): $\delta_l+\delta_\phi=2$ with 
	$(\delta_l,\delta_\phi)\neq(1,1)$.
\item[\textsf{(IV)}] Exceptional point: $(\delta_l,\delta_\phi)=(1,1)$.
\end{itemize}

\subsubsection*{Overview on the scale of order}
On the rough scale of order only the indices $\delta_l$ and $\delta_\phi$ of the functions $\ms d_l$ and $\ms d_\phi$,
respectively, play a role. Assuming \cref{U154} we have the following picture:
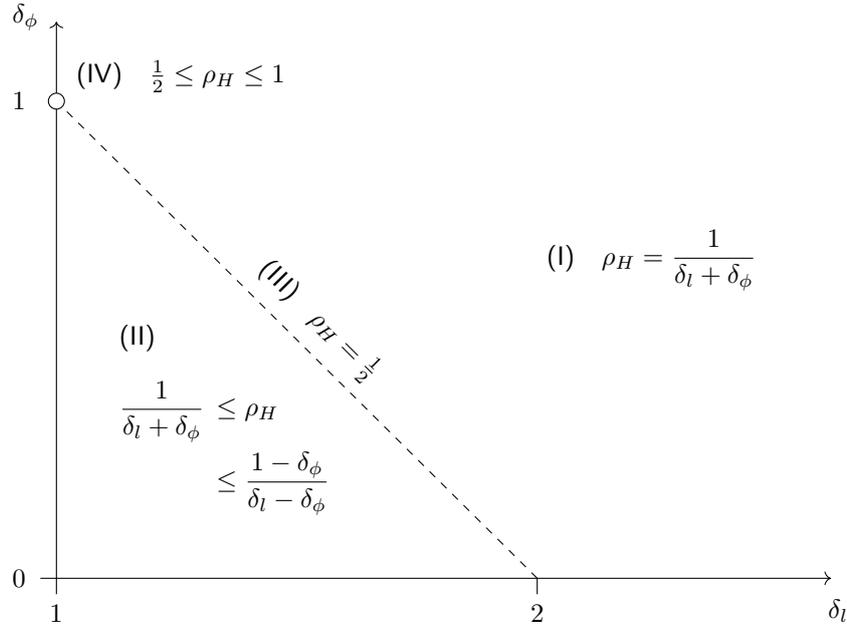
\begin{center}
\begin{tikzpicture}[x=1pt,y=1pt,scale=2,font=\fontsize{10}{10}]
	\draw[->] (7,10)--(155,10);
	\draw (160,0) node[anchor=south east] {$\delta_l$};
	\draw (10,7)--(10,98.5);
	\draw[->] (10,101.5)--(10,115);
	\draw (0,120) node[anchor=north west] {$\delta_\phi$};

	\draw (100,7)--(100,10);
	\draw (0,10) node[anchor=west] {$0$};
	\draw (0,100) node[anchor=west] {$1$};
	\draw (10,0) node[anchor=south] {$1$};
	\draw (100,0) node[anchor=south] {$2$};

	\draw (10,100) circle (1.5);
	\draw[dashed] (100,10)--(11,99);

	\draw (100,70) node[anchor=west] {\textsf{(I)}\quad $\displaystyle \rho_H=\frac 1{\delta_l+\delta_\phi}$};
	\draw (20,60) node[anchor=north west] {%
		$\begin{aligned}
			& \textsf{(II)} && \phantom{}
			\\[6pt]
			& \frac 1{\delta_l+\delta_\phi}\mkern-15mu && \leq\rho_H
			\\ 
			& \phantom{} && \leq\frac{1-\delta_\phi}{\delta_l-\delta_\phi}
		\end{aligned}$%
		};
	\draw (46,70) node[anchor=north west] {%
		\begin{rotate}{-45} \textsf{(III)}\quad $\rho_H=\frac 12$\end{rotate}%
		};
	\draw (12,100) node[anchor=south west] {\textsf{(IV)}\quad $\frac 12\leq\rho_H\leq 1$};
\end{tikzpicture}
\end{center}
The formula for $\rho_H$ in the region of simple behaviour should be seen as a lower bound and an upper bound which coincide in
that region. In this way we may describe the picture as follows:
\begin{Itemize}
\item The lower bound ``$\frac 1{\delta_l+\delta_\phi}\leq \rho_H$'' holds throughout the whole picture. 
\item The upper bound is different in different regions. It coincides with the lower bound in \textsf{(I)}, and 
	splits from it when passing through the switching line. 
\item In the exceptional point \textsf{(IV)}, the upper bound jumps to the trivial bound ``$\rho_H\leq 1$''. 
\end{Itemize}
The bounds in the critical triangle are sharp. 

\begin{Example}[\textsf{(II)}]
\label{U173}
	Let $\delta_l\in[1,2)$ and $\delta_\phi\in[0,2-\delta_l)$.
	\begin{Enumerate}
	\item Set 
		\begin{align*}
			& l_1\DE 1,\qquad l_j\DE
			\begin{cases}
				\frac 1{j^{\delta_l}}\CAS \delta_l>1,
				\\
				\frac 1{j\log^2j}\CAS \delta_l=1,
			\end{cases}
			\quad\text{for }j\geq 2,
			\\
			& \phi_j\DE\sum_{k=1}^j\frac{(-1)^k}{k^{\delta_\phi}},\ j\in\bb N
			.
		\end{align*}
		Then $\rho_H=\frac 1{\delta_l+\delta_\phi}$.
	\item Take the same lengths as in (i), but use
		\[
			\phi_j\DE\sum_{k=1}^j\frac 1{k^{\delta_\phi}},\ j\in\bb N
			,
		\]
		instead. Then $\rho_H=\frac{1-\delta_\phi}{\delta_l-\delta_\phi}$.
	\end{Enumerate}
\end{Example}

\noindent
In the exceptional point all values in the interval $[\frac 12,1]$ can be attained.

\begin{Example}[\textsf{IV}]
\label{U171}
	Let $\nu\in(1,2)$, and consider a Hamburger Hamiltonian $H$ with 
	\begin{equation}
	\label{U172}
		l_j\asymp \frac 1{j\log^\nu j},\quad |\phi_{j+1}-\phi_j|\asymp\frac 1j
		.
	\end{equation}
	Then 
	\[
		\log\Big(\max_{|z|=r}\|W_H(z)\|\Big)\asymp r^{\frac 1\nu}
		.
	\]
	If \cref{U172} holds with ``$\nu=2$'', then $\rho_H=\frac 12$. If 
	\[
		l_j\asymp \frac 1{j(\log j)(\log\log j)^\alpha},\quad |\phi_{j+1}-\phi_j|\asymp\frac 1j
	\]
	with some $\alpha>1$, then $\rho_H=1$. 
\end{Example}

\subsubsection*{Results in detail}
Now we take into account the actual functions $\ms d_l,\ms d_\phi$ and not only their indices. 
First, we state the lower bound which again holds throughout.

\begin{Proposition}[lower bound]
\label{U168}
	Assume that lengths and angles satisfy 
	\[
		l_j\gtrsim\ms d_l(n),\qquad \big|\sin(\phi_{j+1}-\phi_j)\big|\gtrsim\ms d_\phi(n)
		.
	\]
	Let $\IndexN{\ms k(r)}$ be an asymptotic inverse of the function $[\ms d_l(t)\ms d_\phi(t)]^{-1}$. Then 
	\[
		\ms k(r)\lesssim\log|w_{H,22}(ir)|
		.
	\]
	We have $\Ind\ms k(r)=\frac 1{\delta_l+\delta_\phi}$.
\end{Proposition}

\noindent
We come to upper bounds.
In the now considered finer picture the reason for the different behaviour in different regions becomes much clearer. 
\begin{Itemize}
\item Apart from the exceptional point \textsf{(IV)} one has to distinguish between a region of convergence, where 
	$\int_1^\infty[\ms d_l(t)\ms d_\phi(t)]^{\frac 12}\DD t<\infty$, and a region of divergence where this integral
	is infinite. 
	
	Region \textsf{(I)} belongs to the region of convergence, and \textsf{(II)} to the region of
	divergence. In \textsf{(I)} upper and lower bounds coincide, while in \textsf{(II)} they split.  
\item The role played by the switching line \textsf{(III)} is much more subtle than in the previous rough picture. 
	Already in the part of that line belonging to the region of convergence the upper and lower bounds start to split 
	(even before entering the region of divergence or the critical triangle). 
	Within the switching line the gap between the bounds is slowly varying.
\item The case that $\ms d\phi\asymp 1$ (jumping angles) is also exceptional: 
	the upper and lower bounds coincide unless $\delta_l=2$ where they differ by a slowly varying factor.
\item In the exceptional point \textsf{(IV)} we have, in general, 
	no better information than the trivial upper bound $\lesssim r$. 
\end{Itemize}
The exact formulae read as follows.

\begin{proposition}[upper bound; region of convergence]
\label{U169}
	Assume that lengths and angles satisfy 
	\begin{equation}
	\label{U219}
		l_j\lesssim\ms d_l(n),\qquad \big|\sin(\phi_{j+1}-\phi_j)\big|\lesssim\ms d_\phi(n)
		,
	\end{equation}
	and that
	\[
		\int_1^\infty\big[\ms d_l(t)\ms d_\phi(t)\big]^{\frac 12}\DD t<\infty
		,\qquad(\delta_l,\delta_\phi)\neq(1,1)
		.
	\]
	Then 
	\[
		\log\Big(\max_{|z|=r}\|W_H(z)\|\Big)\lesssim 
		r^{\frac 12}\int_{\ms k(r)}^\infty\big[\ms d_l(t)\ms d_\phi(t)\big]^{\frac 12}\DD t
		.
	\]
\end{proposition}

\begin{proposition}[upper bounds; region of divergence]
\label{U170}
	Assume that lengths and angles satisfy \cref{U219} and that 
	\[
		\int_1^\infty\big[\ms d_l(t)\ms d_\phi(t)\big]^{\frac 12}\DD t=\infty
		,\qquad(\delta_l,\delta_\phi)\neq(1,1)
		.
	\]
	Let $\IndexN{\ms h(r)}$ be an asymptotic inverse of $\frac{\ms d_\phi(t)}{\ms d_l(t)}$, which exists since in the present case 
	$\delta_l>\delta_\phi$. 
	Then 
	\[
		\log\Big(\max_{|z|=r}\|W_H(z)\|\Big)\lesssim 
		\begin{cases}
			r\int_{\ms h(r)}^\infty\ms d_l(t)\DD t \CAS \delta_l=1,
			\\[1mm]
			\ms h(r)\ms d_\phi(\ms h(r)) \CAS \delta_l>1\wedge\delta_l+\delta_\phi<2,
			\\[1mm]
			r^{\frac 12}\int_1^{\ms h(r)}[\ms d_l(t)\ms d_\phi(t)]^{\frac 12}\DD t \CAS \delta_l+\delta_\phi=2.
		\end{cases}
	\]
\end{proposition}

\begin{Proposition}[comparison of upper and lower bounds]
\label{U156}
	Denote by $\ms b$ the upper bound exhibited in \Cref{U169,U170} if $(\delta_l,\delta_\phi) \neq (1,1)$, and $\ms b(r)
	\DE r$ if $(\delta_l,\delta_\phi)=(1,1)$. We compare $\ms b(r)$ to the lower bound $\ms k(r)$ from \Cref{U168}.
	\begin{itemize}
	\item[\sf{(I)}] 
		Assume that $\delta_l+\delta_\phi>2$. Then $\ms k(r)\asymp\ms b(r)$. 
	\item[\sf{(II)}] Assume that $\delta_l+\delta_\phi<2$. 
		\begin{Itemize}
		\item If $\delta_\phi>0$, then $\frac{\ms b(r)}{\ms k(r)}$ has positive index. 
		\item If $\delta_\phi=0$ and $\ms d_\phi(t)\ll 1$, then $\frac{\ms b(r)}{\ms k(r)}$ is slowly varying with 
			$\frac{\ms b(r)}{\ms k(r)}\ll 1$.
		\item If $\delta_\phi=0$, $\ms d_\phi(t)\asymp 1$, and $\delta_l>1$, then $\ms k(r)\asymp\ms b(r)$. 
		\item If $\delta_\phi=0$, $\ms d_\phi(t)\asymp 1$, and $\delta_l=1$, then $\frac{\ms b(r)}{\ms k(r)}$ is 
			slowly varying with $\frac{\ms b(r)}{\ms k(r)}\ll 1$.
		\end{Itemize}
	\item[\sf{(III)}] Assume that $\delta_l+\delta_\phi=2$ and $(\delta_l,\delta_\phi)\neq(1,1)$. Then 
		$\frac{\ms b(r)}{\ms k(r)}$ is slowly varying with $\frac{\ms b(r)}{\ms k(r)}\ll 1$.
	\item[\sf{(IV)}] Assume that $\delta_l=\delta_\phi=1$, then $\frac{\ms b(r)}{\ms k(r)}$ has index $\frac 12$. 
	\end{itemize}
\end{Proposition}

\subsubsection*{About sharpness}
In the critical triangle the size of $\log\max_{|z|=r}\|W_H(z)\|$ depends on the distribution of signs of the angle 
differences $\sin(\phi_{j+1}-\phi_j)$. The next result exhibits two extreme cases: eventually constant signs vs.\ nearly 
alternating signs. 

\begin{Proposition}[Region \textsf{(II)}: extreme cases]
\label{U157}
	Assume that lengths and angles satisfy \cref{U154} and that $\delta_l+\delta_\phi<2$.
	\begin{Itemize}
	\item If $\delta_\phi>0$, $\limsup_{j\to\infty}|\phi_{j+1}-\phi_j|<\pi$, and $(\phi_j)_{j=0}^\infty$ is 
		eventually monotone, then 
		\[
			\log\Big(\max_{|z|=r}\|W_H(z)\|\Big)\asymp
			r\int_{\ms h(r)}^\infty \ms d_l(t)\DD t
			.
		\]
	\item If $\delta_l>1$ and there exists $\psi\in\bb R$ such that 
		$|\sin(\phi_j-\psi)|\lesssim|\sin(\phi_{j+1}-\phi_j)|$ for $j \in \bb N$, then 
		\[
			\log\Big(\max_{|z|=r}\|W_H(z)\|\Big)\asymp\ms k(r)
			.
		\]
	\end{Itemize}
\end{Proposition}

\REMARKS{%
\item The lower bound in \Cref{U168} is just \Cref{U47}. The upper bounds in \Cref{U169,U170} follow from \Cref{U45} by making
	appropriate choices for the additional data $\ms c_l,\ms c_\phi$. In most cases one can use 
	\[
		\ms c_l(t)=\ms c_\phi(t)\DE \int_t^\infty\ms d_l(s)\DD s
		,
	\]
	in one boundary case we use 
	\[
		\ms c_l(t)\DE \int_t^\infty\ms d_l(s)\DD s,\quad
		\ms c_\phi(t)\DE \int_t^\infty\ms d_l(s)\DD s\cdot\Big(\int_t^\infty\ms d_\phi(s)\DD s\Big)^2
		.
	\]
	For details see \cite[Corollary~4.7]{pruckner.reiffenstein.woracek:sinqB-arXiv}.
	\Cref{U157} is contained in \cite[Theorem~1.2]{reiffenstein:kachamA-arXiv}.
	\Cref{U173} follows from \Cref{U157}, and \Cref{U171} is \cite[Proposition~3.4]{reiffenstein:kachamA-arXiv} combined with 
	the monotonicity result \Cref{U181}\,(i) for the boundary cases $\nu\in\{1,2\}$. 
	Finally, \Cref{U156} is a simple application of Karamata's theorem. 
\item The bound in the second case of \Cref{U170} can be rewritten in different forms. 
	First, by the definition of $\ms h$, we have 
	\[
		\ms h(r)\ms d_\phi(\ms h(r))=r\ms h(r)\ms d_l(\ms h(r))
		,
	\]
	and second, by Karamata's theorem, we have 
	\[
		r\ms h(r)\ms d_l(\ms h(r))=
		\begin{cases}
			r\int_{\ms h(r)}^\infty\ms d_l(t)\DD t \CAS \delta_l>1,
			\\[2mm]
			r^{\frac 12}\int_1^{\ms h(r)}[\ms d_l(t)\ms d_\phi(t)]^{\frac 12}\DD t \CAS \delta_l+\delta_\phi<2.
		\end{cases}
	\]
\item The intuition behind the requirement in the second item of \Cref{U157} is that $\phi_j$ should converge to $\psi$ such
	that the signs of $\phi_j-\psi$ are roughly alternating. In fact, if these signs are exactly alternating, then we have
	$|\phi_j-\psi|\leq |\phi_{j+1}-\phi_j|$ for all $j$. Contrasting this, in the first item the series of angles 
	diverges as fast as it possibly can.
}

\section[{Hamburger Hamiltonians with power-log decay}]{Hamburger Hamiltonians with power-log decay}
\label{U149}

We consider a very concrete scale of examples to illustrate the results from \Cref{U119,U150,U134}.
In order to give an efficient presentation, we use the following notation. 

\begin{Definition}
\label{U151}
	Let $f$ be a function defined on $(0,\infty)$ or $\bb N$ (or similar) and taking values in $(0,\infty)$. Further, let 
	$(\delta,\alpha)\in\bb R\times\bb R$. Then we write
	\[
		f\circeq(\delta,\alpha) \quad\DI\quad f(t)\asymp t^{-\delta}(\log t)^{-\alpha}\ \text{for $t$ suff.\ large}
	\]
	and analogously ``$f\circleq(\delta,\alpha)$'' or ``$f\circgeq(\delta,\alpha)$'', if the above holds with 
	``$\lesssim$'' or ``$\gtrsim$'', respectively. 
	
	Moreover, let $\preceq$ denote the lexicographical order on $\bb R\times\bb R$. 
\end{Definition}

\noindent 
We consider a Hamburger Hamiltonian in limit circle case, and compare its lengths $l_n$, its angle-differences 
$|\sin(\phi_{n+1}-\phi_n)|$, and the speed of possible convergence of its angles $|\sin(\phi_n-\psi)|$ towards some limit angle 
$\psi$, to power-log functions:
\[
	l_n\text{ vs.\ }(\delta_l,\alpha_l),\quad |\sin(\phi_{n+1}-\phi_n)|\text{ vs.\ }(\delta_\phi,\alpha_\phi),
	\quad |\sin(\phi_n-\psi)|\text{ vs.\ }(\lambda,\kappa)
\]
Since $(l_n)_{n=1}^\infty\in\ell^1$, $(|\sin(\phi_{n+1}-\phi_n)|)_{n=1}^\infty$ is bounded, and $|\sin(\phi_n-\psi)|$
compares to the telescoping series of angle-differences, it is natural to put the following assumptions on the parameters 
$\delta_l,\delta_\phi,\lambda,\alpha_l,\alpha_\phi,\kappa$:
\begin{align}
	\label{U266}
	&\,\, (\delta_l,\alpha_l)\succ(1,1),\quad (\delta_\phi,\alpha_\phi)\succeq(0,0),
	\\[1mm]
	\label{U226}
	& \left.\begin{cases}
		(0,0) \CAS (\delta_\phi,\alpha_\phi)\preceq(1,1),
		\\
		(0,\alpha_\phi-1) \CAS \delta_\phi=1,\alpha_\phi>1,
		\\
		(\delta_\phi-1,\alpha_\phi) \CAS \delta_\phi>1,
	\end{cases}
	\right\}
	\preceq(\lambda,\kappa)\preceq(\delta_\phi,\alpha_\phi)
	.
\end{align}

\noindent
Throughout this section we keep this notation and the standing assumptions \cref{U266} and \cref{U226}.

We always have the lower bound from \Cref{U47}.

\begin{Proposition}[lower bound]
\label{U174}
	Assume that lengths and angles satisfy 
	\[
		l_n\circgeq(\delta_l,\alpha_l),\quad |\sin(\phi_{n+1}-\phi_n)|\circgeq(\delta_\phi,\alpha_\phi)
		.
	\]
	Then 
	\[
		-\bigg(\frac 1{\delta_l+\delta_\phi},-\frac{\alpha_l+\alpha_\phi}{\delta_l+\delta_\phi}\bigg)
		\circleq\log|w_{H,22}(ir)|
		.
	\]
\end{Proposition}

\noindent
For the discussion of upper bounds we
distinguish four parameter regions in the $(\delta_l,\delta_\phi,\lambda)$-space. This is parallel to the case distinction 
\textsf{(I)}--\textsf{(IV)} which appeared in \Cref{U134}, yet slightly refined.
\begin{itemize}
\item[\textsf{(I)}] $\delta_l+\delta_\phi>2$.
\item[\textsf{(II)}] $\delta_l+\delta_\phi<2$.
\item[\textsf{(III)}] $\delta_l+\delta_\phi=2$ with $(\delta_l,\delta_\phi,\lambda)\neq(1,1,0)$.
\item[\textsf{(IV)}] $(\delta_l,\delta_\phi,\lambda)=(1,1,0)$.
\end{itemize}
Unless we are in region \textsf{(I)}, where upper and lower bounds coincide, there will occur further case distinctions. 
Contrasting \Cref{U134} we can also treat the exceptional point \textsf{(IV)}.

\begin{proposition}[upper bound; \textsf{I}]
\label{U162}
	Assume that 
	\[
		l_n\circleq(\delta_l,\alpha_l),\quad |\sin(\phi_{n+1}-\phi_n)|\circleq(\delta_\phi,\alpha_\phi)
		,
	\]
	and that $\delta_l+\delta_\phi>2$. Then 
	\begin{equation}
	\label{U166}
		\log\Big(\max_{|z|=r}\|W_H(z)\|\Big)\circleq
		-\bigg(\frac 1{\delta_l+\delta_\phi},-\frac{\alpha_l+\alpha_\phi}{\delta_l+\delta_\phi}\bigg)
		.
	\end{equation}
\end{proposition}

\begin{proposition}[upper bound; \textsf{II}]
\label{U163}
	Assume that 
	\begin{equation}
	\label{U161}
		l_n\circleq(\delta_l,\alpha_l),\quad |\sin(\phi_{n+1}-\phi_n)|\circleq(\delta_\phi,\alpha_\phi),
		\quad |\sin(\phi_n-\psi)|\circleq(\lambda,\kappa)
		,
	\end{equation}
	and that $\delta_l+\delta_\phi<2$. Then 
	\[
		\log\Big(\max_{|z|=r}\|W_H(z)\|\Big)\circleq
		-\bigg(\frac{1-\delta_\phi+\lambda}{\delta_l-\delta_\phi+2\lambda},
		-\Gamma_{\rm II}\big(\delta_l,\delta_\phi,\lambda,\alpha_l,\alpha_\phi,\kappa\big)\bigg)
		,
	\]
	where
	\[
		\Gamma_{\rm II}\DE
		\begin{cases}
			\frac{\kappa(2-\delta_l-\delta_\phi)+\alpha_l(1-\delta_\phi+\lambda)+\alpha_\phi(\delta_l-1+\lambda)}{
			\delta_l-\delta_\phi+2\lambda} \CAS \delta_l>1,
			\\
			\frac{\kappa(1-\delta_\phi)+\alpha_l(1-\delta_\phi+\lambda)+\alpha_\phi\lambda-\frac 12(1-\delta_\phi)}{
			1-\delta_\phi+2\lambda} \CAS \delta_l=1,\lambda>0,
			\\
			\alpha_l+\kappa-1 \CAS \delta_l=1,\lambda=0,\kappa\geq 1,
			\\
			\alpha_l-1 \CAS \delta_l=1,\lambda=0,\kappa<1.
		\end{cases}
	\]
\end{proposition}

\begin{proposition}[upper bound; \textsf{III}]
\label{U164}
	Assume \cref{U161} and that $\delta_l+\delta_\phi=2$ and $(\delta_l,\delta_\phi,\lambda)\neq(1,1,0)$. Then 
	\[
		\log\Big(\max_{|z|=r}\|W_H(z)\|\Big)\circleq
		-\bigg(\frac 12,-\Gamma_{\rm III}\big(\delta_l,\delta_\phi,\lambda,\alpha_l,\alpha_\phi\big)\bigg)
		,
	\]
	where
	\[
		\Gamma_{\rm III}\DE
		\begin{cases}
			\frac{\alpha_l+\alpha_\phi-2}{2(1-\delta_\phi+\lambda)} \CAS \alpha_l+\alpha_\phi<2, \delta_l>1, 
			\\[2mm]
			\frac{\alpha_l+\alpha_\phi-2}2 \CASO
			.
		\end{cases}
	\]
\end{proposition}

\begin{proposition}[upper bound; \textsf{IV}]
\label{U165}
	Assume \cref{U161} and that $(\delta_l,\delta_\phi,\lambda)=(1,1,0)$. Then 
	\[
		\log\Big(\max_{|z|=r}\|W_H(z)\|\Big)\circleq
		\begin{cases}
			-\big(\frac 12,-\frac{\alpha_l+\alpha_\phi-2}2\big) \CAS \alpha_l+\alpha_\phi>2,
			\\
			-\big(\frac 12,1\big) \CAS \alpha_l+\alpha_\phi=2,
			\\
			-\big(\frac{1-\alpha_\phi+\kappa}{\alpha_l-\alpha_\phi+2\kappa},0\big) \CAS \alpha_l+\alpha_\phi<2.
		\end{cases}
	\]
\end{proposition}

\begin{Proposition}[Region \textsf{(II)}: extreme cases]
\label{U175}
	Assume that 
	\[
		l_n\circeq(\delta_l,\alpha_l),\quad |\sin(\phi_{n+1}-\phi_n)|\circeq(\delta_\phi,\alpha_\phi),
	\]
	and that $\delta_l+\delta_\phi<2$.
	\begin{Itemize}
	\item If $\delta_\phi>0$, $\limsup_{j\to\infty}|\phi_{j+1}-\phi_j|<\pi$ and $(\phi_j)_{j=0}^\infty$ is 
		eventually monotone, then 
		\[
			\log\Big(\max_{|z|=r}\|W_H(z)\|\Big)\circeq
			-\bigg(\frac{1-\delta_\phi}{\delta_l-\delta_\phi},
			-\frac{\alpha_l (1-\delta_\phi)+\alpha_\phi (\delta_l-1)}{\delta_l-\delta_\phi}\bigg)
			.
		\]
	\item If $\delta_l>1$ and there exists $\psi\in\bb R$ such that 
		$|\sin(\phi_j-\psi)|\lesssim|\sin(\phi_{j+1}-\phi_j)|$, then 
		\[
			\log\Big(\max_{|z|=r}\|W_H(z)\|\Big)\circeq
			-\bigg(\frac 1{\delta_l+\delta_\phi},-\frac{\alpha_l+\alpha_\phi}{\delta_l+\delta_\phi}\bigg)
			.
		\]
	\end{Itemize}
\end{Proposition}

\REMARKS{%
\item \Cref{U162,U163,U164}, and the case in \Cref{U165} where $\alpha_l+\alpha_\phi>2$, follow
	from \Cref{U45}. The case $\alpha_l+\alpha_\phi\leq 2$ in \Cref{U165} follows from 
	\cite[Examples~5.8 and 5.9]{pruckner.reiffenstein.woracek:sinqB-arXiv}; in the extreme case $\lambda=\delta_\phi>0$
	an even better bound is available there. 
\item The lower bound \Cref{U174} is \Cref{U47}. The lower bound in
	\Cref{U175}\,(i) comes from \Cref{U148}.
\item For every choice of parameters satisfying \cref{U266} and \cref{U226} there is a Hamburger Hamiltonian satisfying 
	\cref{U161} with equality. An example is obtained as follows. Set $l_1\DE 1,\phi_1\DE 0$ and
	\[
		l_n\DE\frac 1{n^{\delta_l}\log^{\alpha_l}n},\quad 
		\phi_n\DE\sum_{j=1}^{n-1}\frac{\epsilon_n}{n^{\delta_\phi}\log^{\alpha_\phi}n}
		\qquad\text{for }n\geq 2, 
	\]
	where $\epsilon_n\in\{-1,1\}$ are chosen such that 
	\[
		(\lambda,\kappa)\succ(0,0)\ \Rightarrow\ \psi\DE\lim_{j\to\infty}\phi_j\text{ exists }\wedge
		|\sin(\phi_n-\psi)|\circeq(\lambda,\kappa).
	\]
}

\section[{Bochkov's theorem about type}]{Bochkov's theorem about type}
\label{U105}

Up to our knowledge there is only one example of an indeterminate moment problem where not only the order but also the type 
of the Nevanlinna matrix can be computed without knowing explicit expressions or explicit asymptotics for its entries. 

This example is given in terms of the Jacobi parameters, namely, taking $a_n\DE 0$ and $b_n\DE(n+1)^p$ where $p>1$. 
The corresponding moment problem is indeterminate by Berezanskii's Theorem, cf.\ \Cref{U12}.

\begin{Theorem}
\label{U9}
\IndexxS{Theorem!Bochkov}
	Let $p>1$ and consider the moment problem corresponding to the Jacobi parameters 
	\[
		a_n\DE 0,\ b_n\DE(n+1)^p\quad\text{for }n\in\bb N
		.
	\]
	This problem is indeterminate and the entries of the Nevanlinna matrix are of order $\frac 1p$ with type 
	\begin{equation}
	\label{U11}
		p\int\limits_0^1\frac{\DD t}{(1-t^{2p})^{\frac 1p}}
		.
	\end{equation}
\end{Theorem}

\REMARKS{%
\item The values for order and type were conjectured by G.~Valent in \cite[Conjecture~1]{valent:1999} on the basis of two
	explicit examples. His starting point was the study of birth-and-death processes with polynomial rates. An explicit
	translation to the presently stated form can be found in \cite{berg.szwarc:2017} (in particular, Remark~1.10 of that
	paper). 
\item The statement about order in \Cref{U9} was first established in \cite[Corollary~6]{romanov:2017} by a straightforward 
	application of Romanov's Theorem~I, cf.\ \Cref{U43}. A different proof can be found in \cite[\S6]{berg.szwarc:2017}. 
	In the latter paper also bounds for the type are given. 
\item The fact that the type is equal to the value \cref{U11} is shown by I.~Bochkov in \cite{bochkov:2021}. 
	He analyses the Taylor coefficients of the power series appearing in the Berg-Szwarc Theorem 
	(\Cref{U28} and specifically \cite[Theorem~1.11]{berg.szwarc:2017}). This analysis relies on a
	combinatorial trick which seems to be quite specific for the situation that the Jacobi parameters are powers (though, as
	far as we know, no serious attempts were made to generalise the approach). 
}

\section[{Modifying the rotation of the Hamiltonian}]{Modifying the rotation of the Hamiltonian}
\label{U176}

It is a very intuitive fact that the density of the eigenvalues should not increase when the rotation
of the Hamiltonian becomes slower, and we give two results that instanciate this intuition.
As in \Cref{U110} we denote the entries of the Hamiltonian $H$ as follows,
\[
	H(t)=\begin{pmatrix} h_1(t) & h_3(t)\\ h_3(t) & h_2(t)\end{pmatrix},\qquad t\in(a,b)\text{ a.e.},
\]
and set $\omega_{H,j}(s,t) \DE \int_s^t h_j(u) \DD u$ for $j=1,2,3$. 
As usual, we use $\det\Omega_H(s,t)$ as a measure for the speed of rotation, and the growth of the eigenvalue counting function $n_H(r)$ is
expressed via its Stieltjes transform.

\begin{Theorem}
\label{U177}
\IndexxS{Theorem!comparing two Hamiltonians}
	Let $H,\wt H\in\bb H_{a,b}$ be definite, assume that $\int_a^b h_1(t) \DD t<\infty$, and that $H, \wt H$ are either both in the limit circle or both in the limit point case. Assume that the spectra of $A_H$ and $A_{\wt H}$ are discrete.
\begin{Enumerate}
\item If there are $\gamma_-,\gamma_+>0$ such that
	\begin{alignat}{2}
		\label{U178}
		& \det\Omega_{\wt H}(s,t) \le \gamma_+^2\det\Omega_H(s,t),\qquad && a<s<t<b,
		\\[1ex]
		\label{U179}
		& \omega_{\wt H,2}(a,t)\le \gamma_+\omega_{H,2}(a,t),\qquad && t\in(a,b), 
		\\[1ex]
		\label{U180}
		& \gamma_-h_1(t) \le \tilde h_1(t) \le \gamma_+h_1(t),\qquad && t\in(a,b)\text{ a.e.}
	\end{alignat}
	then (for say, $r\geq 1$)
	\[
		\int_0^\infty\frac 1{t+r^2}\cdot\frac{n_{\wt H}(\sqrt t)}t\DD t\lesssim
		\int_0^\infty\frac 1{t+r^2}\cdot\frac{n_H(\sqrt t)}t\DD t
	\]
	The constant implicit in the relation ``\/$\lesssim$'' depends on $\gamma_-,\gamma_+$ but not on $H,\wt H$. 
\item Assume that $H,\wt H$ are both in limit circle case and that \eqref{U178} holds for some $\gamma_+>0$. Then (for say, $r\geq 1$)
\[
		\int_0^\infty\frac 1{t+r^2}\cdot\frac{n_{\wt H}(\sqrt t)}t\DD t\lesssim \log r \cdot
		\int_0^\infty\frac 1{t+r^2}\cdot\frac{n_H(\sqrt t)}t\DD t,
	\]
where the constant implicit in the relation ``\/$\lesssim$'' depends on $\gamma_+$ but not on $H,\wt H$.
	\end{Enumerate}
\end{Theorem}

\REMARKS{%
\item This result is shown in \cite[Theorem~7.1]{langer.reiffenstein.woracek:kacest-arXiv}.
\item Assume that $\wt H(t)\le H(t)$ for a.e.\ $t\in(a,b)$. Then \cref{U178}, \cref{U179}, and the second inequality in 
	\cref{U180} are automatically satisfied. Thus, it is enough to assume additionally that $\tilde h_1\lesssim h_1$ 
	in order to ensure that \Cref{U177} (i) is applicable. 
}

\section[{Cutting out pieces of the Hamiltonian}]{Cutting out pieces of the Hamiltonian}
\label{U121}

It is a very intuitive fact that the growth of the monodromy matrix should not increase when one cuts out pieces of the
Hamiltonian, and this intuition can indeed be established. We give two variants of this fact.

First, we have to make precise what we mean by ``cutting out pieces'' of a Hamiltonian. To this end, assume we have 
$H\in\bb H_{a,b}$ (which is in limit circle case) and let $\Delta\subseteq[a,b]$ be measurable with nonzero Lebesgue measure. 
Then we set
\begin{align*}
	\lambda(t) &\DE \int_a^t\mathds{1}_{\Delta}(u)\DD u,\quad t\in[a,b], \qquad\qquad
	\tilde a\DE 0,\ \tilde b\DE\lambda(b),
	\\[1ex]
	\chi(s) &\DE \min\bigl\{t\in[a,b]\DS \lambda(t)\ge s\bigr\},\quad s\in[\tilde a,\tilde b],
	\\[1ex]
	\wt H(s) &\DE H(\chi(s)),\quad s\in[\tilde a,\tilde b].
\end{align*}
It is clear that $\wt H\in\bb H_{\tilde a,\tilde b}$, and we think of $\wt H$ as the Hamiltonian that emerged from $H$ by
cutting out the complement of $\Delta$ from $[a,b]$. 

\begin{Theorem}
\label{U181}
\IndexxS{Theorem!cutting out pieces}
	Let $H\in\bb H_{a,b}$ be in limit circle case, let $\Delta\subseteq[a,b]$ be measurable with nonzero Lebesgue measure, 
	and let $\wt H\in\bb H_{\tilde a,\tilde b}$ be the Hamiltonian that emerges by cutting out the complement of $\Delta$. 
	Assume that $\wt H$ is definite. 
	\begin{Enumerate}
	\item We have 
		\[
			\int_0^\infty\frac 1{t+r^2}\cdot\frac{n_{\wt H}(\sqrt t)}t\DD t\lesssim
			\log r\cdot\int_0^\infty\frac 1{t+r^2}\cdot\frac{n_H(\sqrt t)}t\DD t
		\]
		for sufficiently large $r$.
	\item Assume additionally that for every $H$-indivisible interval $(c,d)\subseteq(a,b)$ either 
		$(c,d)\cap\Delta$ or $(c,d)\setminus\Delta$ has measure zero. Then 
		\[
			n_{\wt H}(r)\leq n_H(r)\quad\text{for all }r>0
			.
		\]
	\end{Enumerate}
\end{Theorem}

\REMARKS{%
\item Item (i) in \Cref{U181} is \cite[Theorem~5.10]{langer.reiffenstein.woracek:kacest-arXiv}, and 
	item (ii) follows from \cite[Theorem~3.4]{pruckner.woracek:sinqA}. Interestingly, the proofs of these two results are
	very different. The first is purely analytic and the second is operator theoretic. 
\item We comment on the proof of \Cref{U181}\,(i). 
	It is straightforward to show that $\kappa_H(r)\geq\kappa_{\wt H}(r)$ for all $r>0$
	where $\kappa_H$ and $\kappa_{\wt H}$ are as in \Cref{U143} with some constant $c>0$. Namely, 
	one proceeds executing the algorithm in \Cref{U143} and notes that 
	\[
		\det\Omega_{\tilde H}\big(\lambda(\sigma_{j-1}^{(r)}),\lambda(\sigma_j^{(r)})\big)\leq
		\det\Omega_H\big(\sigma_{j-1}^{(r)},\sigma_j^{(r)}\big)\leq\frac 1{r^2}
		,
	\]
	which leads to $\kappa_H(r)\geq\kappa_{\tilde H}(r)$. 
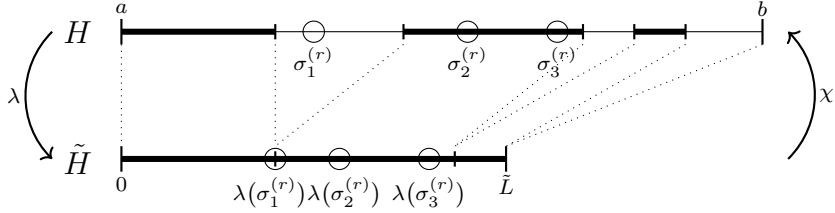
	\begin{center}
		\begin{tikzpicture}[x=1.2pt,y=1.2pt,scale=0.8,font=\fontsize{8}{8}]
			\draw (10,30)--(260,30);
			\draw[line width=.8mm] (10,30)--(70,30);
			\draw[line width=.8mm] (120,30)--(190,30);
			\draw[line width=.8mm] (210,30)--(230,30);
			\draw[thick] (10,25)--(10,35);
			\draw[thick] (70,27)--(70,33);
			\draw[thick] (120,27)--(120,33);
			\draw[thick] (190,27)--(190,33);
			\draw[thick] (210,27)--(210,33);
			\draw[thick] (230,27)--(230,33);
			\draw[thick] (260,25)--(260,35);
			\draw (-7,30) node {\large $H$};
			\draw (10,39) node {${\displaystyle a}$};
			\draw (260,39) node {${\displaystyle b}$};
			
			\draw (85,30) circle (5pt);
			\draw (145,30) circle (5pt);
			\draw (180,30) circle (5pt);
			
			\draw (85,18) node {${\sigma_1^{(r)}}$};
			\draw (145,18) node {${\sigma_2^{(r)}}$};
			\draw (180,18) node {${\sigma_3^{(r)}}$};
			
			\draw[dotted] (10,25)--(10,-15);
			\draw[dotted] (70,25)--(70,-15);
			\draw[dotted] (120,25)--(70,-15);
			\draw[dotted] (190,25)--(140,-15);
			\draw[dotted] (210,25)--(140,-15);
			\draw[dotted] (230,25)--(160,-15);
			\draw[dotted] (260,25)--(160,-15);
			
			\draw[line width=.8mm] (10,-20)--(160,-20);
			\draw[thick] (10,-25)--(10,-15);
			\draw[thick] (70,-23)--(70,-17);
			\draw[thick] (140,-23)--(140,-17);
			\draw[thick] (160,-25)--(160,-15);
			\draw (-7,-20) node {\large $\tilde H$};
			\draw (10,-30) node {${\displaystyle 0}$};
			\draw (-32,5) node {$\lambda$};
			\draw (285,5) node {$\chi$};
			
			\draw (70,-20) circle (5pt);
			\draw (95,-20) circle (5pt);
			\draw (130,-20) circle (5pt);
			
			\draw (68,-33) node {${\lambda \big(\sigma_1^{(r)}\big)}$};
			\draw (97,-33) node {${\lambda \big(\sigma_2^{(r)}\big)}$};
			\draw (130,-33) node {${\lambda \big(\sigma_3^{(r)}\big)}$};		
			
			\draw (160,-30) node {${\displaystyle \tilde L}$};
			
			\draw[->, thick] (-17,30) to[out=225, in=135] (-17,-20);
			\draw[->, thick] (270,-20) to[out=45, in=315] (270,30);
		\end{tikzpicture}
	\end{center}
\item We comment on the proof of \Cref{U181}\,(ii). Under the stated additional assumption one can relate the operator models
	associated with $H$ and $\wt H$. First, the map $V$ acting as 
	\[
		V\DF f\mapsto (f\circ\lambda)\cdot\mathds{1}_{\Delta}
	\]
	induces an isometry of $L^2(\tilde H)$ into $L^2(H)$. Now consider the multiplication operator 
	$M_{\mathds{1}_{\Delta}}$ with $\mathds{1}_{\Delta}$, and the Volterra operator $R_H$ defined as the inverse of the
	restriction of the maximal operator with boundary condition ``$f(a)=0$''. Then 
	$\Ran M_{\mathds{1}_{\Delta}}R_HV\subseteq\Ran V$ and $R_{\tilde H}=V^{-1}M_{\mathds{1}_{\Delta}}R_HV$. 
}
	\begin{center}
	\begin{tikzcd}[column sep=normal]
		L^2(H) \arrow{r}{R_H} & L^2(H) \arrow{r}{M_{\mathds{1}_{\Delta}}} 
		& \Ran(M_{\mathds{1}_{\Delta}}R_HV) \arrow{d}{\subseteq}
		\\
		L^2(\tilde H) \arrow{r}[swap]{R_{\tilde H}} \arrow{u}{V}
		& L^2(\tilde H) \arrow[r, bend left=15, "V"] \ar[r, phantom, "\cong"] 
		& \Ran V \ar[l, bend left=15, "V^{-1}"]
	\end{tikzcd}
	\end{center}

\section[{Exploiting symmetry}]{Exploiting symmetry}
\label{U131}

We call a Hamiltonian $H=\smmatrix{h_1}{h_3}{h_3}{h_2}$ \IndexS{diagonal}{Hamiltonian!diagonal} if $h_3=0$ a.e. 
If $H$ is diagonal many formulae
simplify because the canonical system splits and can easily be rewritten as a scalar second order equation of familiar form. 
On the level of solutions of the system (and Weyl functions, spectrum, etc.) this is expressed by symmetry properties. 
Diagonal Hamiltonians can be transformed to non-diagonal ones of a very particular kind via a splitting procedure, and this
process can be reversed.

First we state the mentioned symmetry properties.

\begin{Theorem}
\label{U299}
	Let $H\in\bb H_{a,b}$. 
	\begin{Enumerate}
	\item Assume that $H$ is in limit circle case, and consider the monodromy matrix $W_H=(w_{H,ij})_{i,j=1}^2$ of $H$. 
		Then $H$ is diagonal if and only if the functions $w_{H,11},w_{H,22}$ are even and $w_{H,12},w_{H,21}$ are odd. 
	\item Assume that $H$ is in limit point case, and consider the Weyl coefficient $q_H$ of $H$. 
		Then $H$ is diagonal if and only if the function $q_H$ is odd.
	\end{Enumerate}
\end{Theorem}

\noindent
Second we explain the splitting procedure: a diagonal Hamiltonian gives rise to two non-diagonal ones.

\begin{Definition}
\label{U298}
	Let $H\in\bb H_{a,b}$ be diagonal. Set 
	\[
		\check m(t)\DE\int_a^t h_1(s)\DD s,\quad \hat m(t)\DE\int_a^t h_2(s)\DD s
		,
	\]
	and let $\check\rho$ and $\hat\rho$ be the left-continuous right inverses of $\check m$ and $\hat m$, respectively. 
	\begin{center}
		\begin{tikzpicture}[x=1.2pt,y=1.2pt,scale=0.8,font=\fontsize{8}{8}]
			\draw (20,0) node[anchor=east] {$[a,b]$};
			\draw (40,0) node[anchor=west] {$[0,\int_a^bh_1(x)\,dx]$};
			\draw[->] (20,2) to[out=30, in=120] (40,2);
			\draw[->] (40,-2) to[out=210, in=300] (20,-2);
			\draw (30,12) node {$\check m$};
			\draw (30,-12) node {$\check\rho$};
			\draw (30,-25) node {$\check m\circ\check\rho=\Id$};
			\draw (170,0) node[anchor=east] {$[a,b]$};
			\draw (190,0) node[anchor=west] {$[0,\int_a^bh_2(x)\,dx]$};
			\draw[->] (170,2) to[out=30, in=120] (190,2);
			\draw[->] (190,-2) to[out=210, in=300] (170,-2);
			\draw (180,12) node {$\hat m$};
			\draw (180,-12) node {$\hat\rho$};
			\draw (180,-25) node {$\hat m\circ\hat\rho=\Id$};
		\end{tikzpicture}
	\end{center}
	\begin{Enumerate}
	\item Assume that $H$ is in limit circle case. Then we define 
		$\check H\DF(0,\int_a^b h_2(s)\DD s)\to\bb R^{2\times 2}$ and 
		$\hat H\DF(0,\int_a^b h_1(s)\DD s)\to\bb R^{2\times 2}$ as 
		\begin{align}
		\label{U297}
			\check H(t)\DE &\, 
			\begin{pmatrix}
				(\check m\circ\hat\rho)(t)^2 & (\check m\circ\hat\rho)(t)
				\\
				(\check m\circ\hat\rho)(t) & 1
			\end{pmatrix}
			\\
		\label{U296}
			\hat H(t)\DE &\, 
			\begin{pmatrix}
				1 & -(\hat m\circ\check\rho)(t)
				\\
				-(\hat m\circ\check\rho)(t) & (\hat m\circ\check\rho)(t)^2
			\end{pmatrix}
		\end{align}
	\item Assume that $H$ is in limit point case. Then we define $\check H,\hat H\DF(0,\infty)\to\bb R^{2\times 2}$ by the
		formulae \cref{U297} and \cref{U296}, respectively, and by 
		\begin{align*}
			\check H(t)\DE &\, \begin{pmatrix} 1 & 0 \\ 0 & 0 \end{pmatrix},\quad t>\int_a^b h_1(s)\DD s
			,
			\\
			\hat H(t)\DE &\, \begin{pmatrix} 0 & 0 \\ 0 & 1 \end{pmatrix},\quad t>\int_a^b h_2(s)\DD s
			.
		\end{align*}
		Note that these additional definitions are nonvoid for at most one of $\check H$ and $\hat H$, since 
		$\int_a^b(h_1(s)+h_2(s))\DD s=\infty$. 
	\end{Enumerate}
\end{Definition}

\noindent
Clearly, the functions $\check H$ and $\hat H$ are Hamiltonians on their respective domains. They are both in limit circle case
if $H$ is, and they are both in limit point case if $H$ is. 

The monodromy matrices of $H$ and $\check H,\hat H$ and their Weyl coefficients, respectively, are related by simple explicit
formulae. They perfectly express the nature of symmetry for a diagonal Hamiltonian. 

\begin{Theorem}
\label{U295}
	Let $H$ be a diagonal Hamiltonian. 
	\begin{Enumerate}
	\item Assume that $H$ is in limit circle case. Then 
		\begin{align*}
			W_{\check H}(t,z^2)= &\, 
			\smmatrix{\tfrac 1z}001 W_H(t,z)\smmatrix z001\smmatrix 1{w_{H,12}'(t,0)}01
			,
			\\
			W_{\hat H}(t,z^2)= &\, 
			\smmatrix 100{\tfrac 1z} W_H(t,z)\smmatrix 100z\smmatrix 10{-w_{H,21}'(t,0)}1
			.
		\end{align*}
	\item Assume that $H$ is in limit point case. Then 
		\[
			zq_{\check H}(z^2)=q_H(z)=\frac 1z q_{\hat H}(z^2)
			.
		\]
	\end{Enumerate}
\end{Theorem}

\noindent
The relevance of this result for the present setting is that zeroes and growth of $w_{H,22}$ or poles of $q_H$, respectively, 
are related to zeroes and growth or poles of the corresponding functions for $\check H$ and $\hat H$ just by taking square roots. 
Hence, any knowledge about one of $H,\check H,\hat H$ translates immediately to the other two.

\REMARKS{%
\item The proof of \Cref{U299} is contained in \cite{kaltenbaeck.winkler.woracek:bimmel}. \Cref{U295} is established by
	a straightforward computation.
\item As we presented the matters, they look really ad hoc. But this is not at all so; there is a very clean reason behind. 
	It rests on another characterisation of diagonality in terms of the model space $L^2(H)$ which reads as follows: 
	\begin{Itemize}
	\item $H$ is diagonal if and only if the involution
		\[
			\iota\DF\binom{f_1}{f_2}\mapsto\binom{-f_1}{f_2}
		\]
		maps $L^2(H)$ isometrically into itself. 
	\end{Itemize}
	The proof is by elementary computation. 

	Having an isometric involution on a Hilbert space automatically 
	leads to an orthogonal decomposition of the space. Namely, we have the orthogonal
	projections $\check P\DE\frac 12(\Id+\iota)$ and $\hat P\DE\frac 12(\Id-\iota)$, and the corresponding decomposition 
	\[
		L^2(H)=\Ran\check P\oplus\Ran\hat P
		.
	\]
	It turns out that the spaces $\Ran\check P$ and $\Ran\hat P$ can be identified with the model spaces of some 
	Hamiltonians, and these are $\check H$ and $\hat H$: 
	\[
		\Ran\check P\cong L^2(\check H),\quad \Ran\hat P\cong L^2(\hat H)
		.
	\]
	The isomorphisms are known explicitly and the model operators decompose accordingly, cf.\ 
	\cite[\S 2.f]{winkler.woracek:del} and \cite{kaltenbaeck.winkler.woracek:bimmel}. 
\item On the level of Nevanlinna functions, the splitting of an odd function $q$ is also classical and occurred in the context
	of strings, cf.\ \cite{kac.krein:1968a}.
\item The splitting concept is common to all kinds of objects related with canonical systems, and can be extended to a
	sign-indefinite setting. This is in parts still work in progress, see 
	\cite{kaltenbaeck.winkler.woracek:dbsym,kaltenbaeck.winkler.woracek:nksym,kaltenbaeck.winkler.woracek:string,
	winkler.woracek:evoddb,winkler.woracek:varsym}.
}

\section[{Krein strings}]{Krein strings}
\label{U124}

A \IndexS{string}{string} is a pair, we write it as $S[L,\mu]$, where $L\in[0,\infty]$ and $\mu$ is a positive Borel measure 
on $\bb R$ (compact sets are assumed to have finite measure) with $\Supp\mu\subseteq[0,L]$, $\sup\Supp\mu=L$, and $\mu(\{L\})=0$.
The number $L$ is called the \IndexS{length}{string!length} of the string, and the distribution function 
\[
	m(x):=\mu\big((-\infty,x)\big),\qquad x\in [-\infty,\infty]
\]
its \IndexS{mass-function}{string!mass-function}.

When Fourier's method is applied to the partial differential equation that describes the vibrations of a string on the interval 
$[0,L)$ with free left endpoint $0$ whose mass on the interval $[0,x)$ is $m(x)$, then the equation 
\[
	y'(x)+\int_{[0,x]}zy(u)\DD\mu(u)=0
\]
arises, which is called the \IndexS{string equation}{string equation}. 
For each $z\in\bb C$ there exist unique solutions $\varphi(x,z)$ and $\psi(x,z)$ of the string equation 
satisfying the initial conditions
\[
	\varphi(0,z)=1,\ \varphi'(0-,z)=0,\qquad \psi(0,z)=0,\ \psi'(0-,z)=1
	.
\]
The limit
\begin{equation}
\label{U279}
	q_S(z):=\lim_{x\to L}\frac{\psi(x,z)}{\varphi(x,z)}
\end{equation}
exists locally uniformly on $\bb C\setminus(-\infty,0]$, and is called the 
\IndexS{principal Titchmarsh--Weyl coefficient}{principal Titchmarsh--Weyl coefficient} of the string $S[L,\mu]$. 
An important theorem of M.G.~Krein states that the assignment $S[L,\mu]\mapsto q_S$ is a bijection between all strings and all
\IndexS{Stieltjes class}{Stieltjes class} functions, i.e.,
functions $q$ that are analytic in $\bb C\setminus(-\infty,0]$ and admit a representation of the form
\[
	q(z):=b+\int_{[0,\infty)}\frac{d\sigma(t)}{t-z}
	,
\]
where $b\geq 0$ and $\sigma$ is a positive Borel measure with $\int_{[0,\infty)}\frac{d\sigma(t)}{1+t}<\infty$.
The measure occurring in the representation of the function $q_S$ is called the 
\IndexS{principal spectral measure}{principal spectral measure} of the string $S[L,\mu]$. 

A string gives rise to a diagonal Hamiltonian.

\begin{Definition}
\label{U276}
	Let $S[L,\mu]$ be a string. Denote by $\lambda$ the Lebesgue measure, and set $\tau(x)\DE x+m(x)$. 
	Then $\tau$ is an increasing function of $[0,L]$ onto a subset of $[0,\infty]$ whose complement consists of an at most
	countable union of intervals. 
	
	We define 
	\begin{equation}
	\label{U274}
		H_d(t)\DE
		\begin{cases}
			\begin{pmatrix}
				\frac{\DD\lambda}{\DD(\lambda+\mu)}(\tau^{-1}(t)) & 0
				\\
				0 & \frac{\DD\mu}{\DD(\lambda+\mu)}(\tau^{-1}(t))
			\end{pmatrix}
			\CAS t\in\Ran\tau,
			\\[4mm]
			\begin{pmatrix}
				0 & 0
				\\
				0 & 1
			\end{pmatrix}
			\CAS t\in(0,\infty)\setminus\Ran\tau.
		\end{cases}
	\end{equation}
\end{Definition}

\noindent
Clearly, $H_d$ is a trace-normalised Hamiltonian on the interval $(0,\infty)$. Its relation to the string equation is as
follows.

\begin{Theorem}
\label{U278}
\IndexxS{Theorem!string vs. Hamiltonian}
	Let $S[L,\mu]$ be a string, and let $H_d$ be as in \cref{U274}. 
	The fundamental solution $W_d(t,z)$ of the Hamiltonian $H_d$ can be expressed in terms of the functions 
	$\varphi,\psi$ as
	\begin{equation}
	\label{U277}
		W_d(t,z)=
		\begin{pmatrix}
			\psi'(\tau^{-1}(t),z^2) & z\psi(\tau^{-1}(t),z^2)
			\\[2mm]
			\frac 1z\varphi'(\tau^{-1}(t),z^2) & \varphi(\tau^{-1}(t),z^2)
		\end{pmatrix}
		,\qquad t\in\Ran\tau
		.
	\end{equation}
	The complement of $\Ran\tau$ consists of indivisible intervals, and hence the values of $W_d(t,z)$ for 
	$t\in(0,\infty)\setminus\Ran\tau$ can be obtained from \cref{U277} by linear interpolation. 

	The Weyl function $q_{H_d}$ is related to the principal Titchmarsh--Weyl coefficient of the string as 
	\[
		q_{H_d}(z)=zq_S(z^2)
		.
	\]
\end{Theorem}

\noindent
The assignment $S[L,\mu]\mapsto H_d$ can be reversed: using the notation from \Cref{U298}, we have 
\begin{equation}
\label{U275}
	m(x)=(\hat m\circ\check\rho)(x)
	.
\end{equation}
Together, \cref{U274} and \cref{U275} yield a bijective correspondence between strings and (trace normed) diagonal Hamiltonians. 

\REMARKS{%
\item A standard reference for strings is \cite{kac.krein:1968,kac:1994}, or also \cite{dym.mckean:1976}.
	The (operator) theory of Krein strings historically preceeds the theory of canonical systems. In fact, the first 
	initiated a lot of research about the latter. The connection between strings and Hamiltonians was long known as ``common
	knowledge''. An explicit presentation, which includes a detailed study of the corresponding operator models, is given in 
	\cite{kaltenbaeck.winkler.woracek:bimmel}. 
\item The spectrum of a string is encoded in its principal Titchmarsh--Weyl coefficient in the same way as the spectrum of a
	canonical system is encoded in its Weyl function. In particular, $S[L,\mu]$ has discrete spectrum if and only if 
	the corresponding diagonal Hamiltonian $H_d$ has. Moreover, if the spectrum is discrete, the eigenvalues of $H_d$ are 
	obtained from those of $S[L,\mu]$ by taking all square roots. 
\item In the present context we mention in particular the work of I.S.~Kac and M.G.~Krein about the distribution of eigenvalues
	of a string with discrete spectrum 
	\cite{krein:1952a,kac.krein:1958,kac:1959,kac:1962,kac:1973,kac:1978,kac:1986,kac:1990}.
	From the nowadays viewpoint those results appear as particular cases of general results about canonical systems, but
	this is neither fully true nor fair to say: a translation from one setting to the other by means of explicit
	computation is in some cases not known, and the work about strings -- including proof methods -- was a great
	inspiration for the general case. 
}

\section[{Inverse results}]{Inverse results}
\label{U273}

In this section we discuss the problem:
\begin{Itemize}
\item Given a regularly varying function $\ms g$, can we construct a Hamiltonian $H$ that is in limit circle case and
	satisfies
	\begin{equation}
	\label{U271}
		\log\Big(\max_{|z|=r}\|W_H(z)\|\Big)\asymp\ms g(r)
		\quad\text{?}
	\end{equation}
\end{Itemize}
The case $\ms g(r)\asymp r$ is of course trivial: by the Krein-de~Branges formula we can take any Hamiltonian 
whose determinant does not vanish identically (and only such). 
The same is true in the case that $\ms g(r)\asymp\log r$, for which 
\cref{U271} is satisfied for any Hamiltonian consisting only of a finite number of indivisible intervals (and only for those). 
Hence, we may focus on functions $\ms g$ with 
\begin{equation}
\label{U270}
	\log r=\Smallo(\ms g(r)),\quad\ms g(r)=\Smallo(r)
	.
\end{equation}
The first is necessary for elementary reasons (the Cauchy estimates), concerning the second recall \Cref{U196} (ii). 

Existence of Hamiltonian with \cref{U271} can always be granted.

\begin{Theorem}
\label{U269}
\IndexxS{Theorem!monodromy matrix!inverse result on growth}
	Let $\ms g$ be a regularly varying function satisfying \cref{U270}. Then there exists a Hamiltonian $H$ which is in limit
	circle case and such that \cref{U271} holds. 

	If $\int_1^\infty\frac 1{\ms f(t)}\DD t<\infty$, where $\ms f$ is an asymptotic inverse of $\ms g$, then 
	$H$ can be chosen to be a Hamburger Hamiltonian. 
\end{Theorem}

\noindent
\Cref{U269} is not constructive regarding our initially posed problem: in the proof one entry of $W_H$ is constructed, but not
$H$ itself. In the following theorem we present large classes of regularly varying functions, where one can easily write down a
Hamiltonian satisfying \cref{U271}. Solutions to the problem can be of very different kind; we exhibit a discrete and a 
continuous solution. 

\begin{Theorem}
\label{U268}
\IndexxS{Theorem!monodromy matrix!prescribed growth}
	Let $\ms g$ be increasing and smoothly varying with positive index, and let $\ms f$ be the inverse function of $\ms g$. 
	\begin{Enumerate}
	\item Assume that $\Ind\ms g\in(0,\frac 12)\cup(\frac 12,1)$, and set 
		\[
			l_n\DE\frac 1{\ms f(n)},\quad \phi_n\DE n\frac\pi2
			.
		\]
		Then the Hamburger Hamiltonian with lengths $l_n$ and angles $\phi_n$ satisfies \cref{U271}.
	\item Assume that $\Ind\ms g\in(\frac 12,1)$, and let $\phi$ be the inverse function of $\frac r{\ms f(r)}$. Then 
		the Hamiltonian $H(t)\DE\xi_{\phi(t)}\xi_{\phi(t)}^T$, defined on some interval $(0,b)$, satisfies \cref{U271}. 
	\end{Enumerate}
\end{Theorem}

\REMARKS{%
\item \Cref{U269} is taken from \cite[Theorem~3.1,Remark~3.2(iii)]{baranov.woracek:smsub}. Alternative approaches could proceed 
	via \cite[Theorem~5.6]{berg.pedersen:1995} or (at least for small growth) via \cite[\S5]{remling.scarbrough:2020}.
	Item (i) of \Cref{U268} is obtained from \Cref{U156}, item (ii) is from 
	\cite[Theorem~6.13]{langer.reiffenstein.woracek:kacest-arXiv}. 
\item We stated two situations in \Cref{U268} that are perfectly simple and explicit. There are also other cases of functions 
	$\ms g$ for which explicit constructions are possible. For example in 
	\cite[Theorem~6.13]{langer.reiffenstein.woracek:kacest-arXiv} also some $\ms g$ with $\Ind\ms g=\frac 12$ are treated.
	A full (constructive) solution of the problem is not known. 
\item The assumption in \Cref{U268} that $\ms g$ is smoothly varying and increasing is no loss in generality since 
	$\Ind\ms g>0$, cf.\ \Cref{U200}.
}

%**************************************************************************
%***                            Last Change: Mon 31 Mar 2025 12:07
%***   < APPENDIX >
%***
%**************************************************************************

\clearpage
\AUXILIARY{APPENDIX}{Auxiliary notions}{Appendix}
\label{U129}

%%%%%%%%%%%%%%%%%%%%%%%%%%%%%%%%%%%%%%%%%%%%%%%%%%%%%%%%%%%%%%%%%%%%%%%%%%%%%%%%%%%%%%%%%

\Intro{%
	In this short appendix we collect some auxiliary notions and results that are used in the core Parts~I-V, but 
	might not be common knowledge to all readers. We include: regular variation in the sense of J.~Karamata, 
	entire functions of Cartwright class, and $J$-inner matrix functions.
}
\begin{center}
	{\large\bf Table of contents}
\end{center}
\begin{flushleft}
	\S\,\ref{U123}.\ Regularly varying functions\ \dotfill\quad\pageref{U123}
	\\[1mm]
	\S\,\ref{U136}.\ Entire functions of Cartwright class\ \dotfill\quad\pageref{U136}
	\\[1mm]
	\S\,\ref{U130}.\ $J$-inner matrix functions\ \dotfill\quad\pageref{U130}
	\\[1mm]
\end{flushleft}
\makeatother
\renewcommand{\thesection}{\Alph{section}}
\renewcommand{\thelemma}{\Alph{section}.\arabic{lemma}}
\makeatletter
\setcounter{section}{3}
\clearpage

%%%%%%%%%%%%%%%%%%%%%%%%%%%%%%%%%%%%%%%%%%%%%%%%%%%%%%%%%%%%%%%%%%%%%%%%%%%%%%%%%%%%%%%%%

%
%
%
\section[{Regularly varying functions}]{Regularly varying functions}
\label{U123}

In complex analysis the growth of the maximum modulus $\max_{|z|=r}|F(z)|$ of an entire function $F$ is compared to functions of
the form $\exp(\ms g(r))$. The most classical comparison functions are powers $\ms g(r)=r^\rho$, and this leads to the notions
of order and type. Let us recall that for an entire function $F$ one defines the \IndexS{order}{order} of $F$ as 
\begin{align*}
	\IndexN{\rho(F)}\DE &\, 
	\limsup_{r\to\infty}\frac{\log\log(\max_{|z|=r}|F(z)|)}{\log r}
	\\
	= &\, 
	\inf\big\{\rho>0\DS\exists c,c'>0\DQ\forall z\in\bb C\DP |F(z)|\leq ce^{c'r^\rho}\big\}\in[0,\infty]
	.
\end{align*}
If $\rho(F)<\infty$, the \IndexS{type of $F$ w.r.t.\ to its order}{type!w.r.t.\ order} is 
\begin{align*}
	\IndexN{\tau(F)}\DE &\,
	\limsup_{r\to\infty}\frac{\log(\max_{|z|=r}|F(z)|)}{r^{\rho(F)}}
	\\
	= &\,
	\inf\big\{\tau>0\DS\exists c>0\DQ\forall z\in\bb C\DP |F(z)|\leq ce^{\tau r^{\rho(F)}}\big\}\in[0,\infty]
	.
\end{align*}
A refined comparison scale was introduced already at a very early stage by E.~Lindel\"of \cite{lindeloef:1905} who considered 
comparison functions behaving for $r\to\infty$ like
\begin{equation}
\label{U83}
	\ms g(r)\DE r^\alpha\cdot\bigl(\log r\bigr)^{\beta_1}\cdot\bigl(\log\log r\bigr)^{\beta_2}
	\cdot\ldots\cdot
	\bigl(\underbrace{\log\cdots\log}_{\text{\footnotesize$m$\textsuperscript{th} iterate}}r\bigr)^{\beta_m},
\end{equation}
where $\alpha>0$ and $\beta_1,\ldots,\beta_m\in\bb R$. We refer to functions of this form as 
\IndexS{Lindel\"of comparison functions}{Lindel\"of comparison function}.

Functions that are nowadays commonly used as comparison functions are regularly varying functions in Karamata sense.
An up-to-date standard reference is \cite[Chapter~7]{bingham.goldie.teugels:1989}; for other levels of generality see 
\cite{seneta:1976,levin:1980,rubel:1996}, and historically \cite{karamata:1930,karamata:1931,karamata:1931a}. 

\begin{Definition}
\label{U60}
	A function $\ms g\DF[1,\infty)\to(0,\infty)$ is called \IndexS{regularly varying}{regular variation} at $\infty$ with 
	\IndexS{index}{regular variation!index} $\alpha\in\bb R$, if it is measurable and 
	\begin{equation}
	\label{U216}
		\forall \lambda\in(0,\infty)\DP \lim_{r\to\infty}\frac{\ms g(\lambda r)}{\ms g(r)}=\lambda^\alpha
		.
	\end{equation}
	We write $\Ind\ms g$ for the index of a regularly varying function $\ms g$. 
	A regularly varying function with index $0$ is also called \IndexS{slowly varying}{slow variation}.
\end{Definition}

\noindent
Observe that a function of the form \cref{U83} is regularly varying with index $\alpha$.

The notion of type admits an immediate generalisation to arbitrary comparison functions. 

\begin{Definition}
\label{U215}
	Let $\ms g$ be a regularly varying function. For an entire function $F$ we define the 
	\IndexS{type of $F$ w.r.t.\ to the comparison function $\ms g$}{type!w.r.t.\ comparison function} as 
	\begin{align*}
		\IndexN{\tau_{\ms g}(F)}\DE &\, 
		\limsup_{r\to\infty}\frac{\log(\max_{|z|=r}|F(z)|)}{\ms g(r)}
		\\
		= &\,
		\inf\big\{\tau>0\DS\exists c>0\DQ\forall z\in\bb C\DP |F(z)|\leq ce^{\tau\ms g(r)}\big\}\in[0,\infty]
		.
	\end{align*}
\end{Definition}

\noindent
We note the particular case that $\ms g(r)=r$: the number $\tau_r(F)$ is called the 
\IndexS{exponential type}{exponential type} of the function $F$. 

In the sequel we state a number of fundamental theorems on regularly varying functions. 
The basis for many of these results are the following two theorems: the 
\IndexS{uniform convergence theorem}{regular variation!uniform convergence theorem} and the 
\IndexS{representation theorem}{regular variation!representation theorem}.

\begin{Theorem}
\label{U79}
\IndexxS{Theorem!regular variation!uniform convergence}
	Assume that $\ms g$ is regularly varying with index $\alpha$. Then the limit \cref{U216} is attained uniformly for
	\[
		\begin{cases}
			\lambda\in[a,b] \text{ for all } 0<a<b<\infty \CAS \alpha=0,
			\\
			\lambda\in(0,b] \text{ for all } 0<b<\infty \CAS \alpha>0,
			\\
			\lambda\in[a,\infty) \text{ for all } 0<a<\infty \CAS \alpha<0.
		\end{cases}
	\]
\end{Theorem}

\begin{Theorem}
\label{U61}
\IndexxS{Theorem!regular variation!representation theorem}
	Let $\alpha \in \bb R$. A function $\ms g\DF[1,\infty)\to(0,\infty)$ is regularly varying with index $\alpha$ 
	if and only if it has a representation of the form 
	\[
		\ms g(r)=r^\alpha \cdot c(r) \exp \bigg(\int_1^r \epsilon (u) \frac{\DD u}{u} \bigg), \qquad r \in [1,\infty),
	\] 
	where $c,\epsilon$ are measurable, $\lim_{r \to \infty} c(r)=c \in (0,\infty)$, and $\lim_{r \to \infty} \epsilon(r)=0$. 

	If $\ms g$ is slowly varying (i.e., $\alpha=0$) and eventually nondecreasing (nonincreasing), then $\epsilon$ may be 
	taken eventually nonnegative (nonpositive). 
\end{Theorem}

\noindent
It is a legitimate intuition that regularly varying functions fill in the scale of powers, and that a regularly varying 
function with index $\alpha$ behaves roughly like the power $r^\alpha$. The next two results
express this intuition very clearly. 

The first is a variant of the \IndexS{Potter bounds}{regular variation!Potter bounds}.
Here we use the following notation to asymptotically compare two functions. 
Assume $f,g$ are defined on a ray like $[1,\infty)$ and assume positive values. Then we write 
\begin{align*}
	& f\IndexN{\sim}g\ \DI\ \lim_{r\to\infty}\frac{f(r)}{g(r)}=1
	,
	\\
	& f\IndexN{\ll}g\ \DI\ \lim_{r\to\infty}\frac{f(r)}{g(r)}=0
	.
\end{align*}

\begin{Theorem}
\label{U62}
\IndexxS{Theorem!regular variation!Potter bounds}
	Let $\alpha\in\bb R$ and let $\ms g$ be regularly varying with $\Ind\ms g=\alpha$. 
	\begin{Enumerate}
	\item ${\displaystyle
		\forall \epsilon>0\DP r^{\alpha-\epsilon}\ll\ms g(r)\ll r^{\alpha+\epsilon}
		}$
		,		
	\item $\lim_{r \to \infty} \frac{\log \ms g(r)}{\log r}=\alpha$,
	\item For all $\epsilon>0$ the quotients $\frac{\ms g(r)}{r^{\alpha-\epsilon}}$ and 
		$\frac{r^{\alpha+\epsilon}}{\ms g(r)}$ are $\sim$ to an eventually increasing function.
	\end{Enumerate}
\end{Theorem}

\noindent
The second is \IndexS{Karamata theorem}{Karamata theorem!asymptotic integration}%
\IndexxS{regular variation!asymptotic integration} about asymptotic integration. 

\begin{Theorem}
\label{U63}
\IndexxS{Theorem!regular variation!asymptotic integration}
	Let $\ms g$ be regularly varying with index $\alpha\in\bb R$. 
	\begin{Enumerate}
	\item Assume that $\alpha\geq -1$. Then the function $x\mapsto\int_1^x \ms g(t)\DD t$ is regularly
		varying with index $\alpha+1$, and 
		\[
			\lim_{x\to\infty}\bigg(
			\raisebox{3pt}{$x\ms g(x)$}\Big/\,\raisebox{-2pt}{$\int\limits_1^x \ms g(t)\DD t$}
			\bigg)=\alpha+1
			.
		\]
	\item Assume that $\alpha\leq -1$ and $\int_1^\infty \ms g(t)\DD t<\infty$. Then the function 
		$x\mapsto\int_x^\infty \ms g(t)\DD t$ is regularly varying with index $\alpha+1$, and 
		\[
			\lim_{x\to\infty}\bigg(
			\raisebox{3pt}{$x\ms g(x)$}\Big/\,\raisebox{-2pt}{$\int\limits_x^\infty \ms g(t)\DD t$}
			\bigg)=-(\alpha+1)
			.
		\]
	\end{Enumerate}
\end{Theorem}

\noindent
Regularly varying functions $\ms g$ are used to quantify growth for $r\to\infty$, and hence the values of $\ms g(r)$ for small
$r$ are irrelevant. This allows to change $\ms g$ on any finite interval without changing the essence of results, and this
freedom can be used to assume that $\ms g$ has some additional properties. For example any regularly varying function 
can be smoothened.

\begin{Theorem}
\label{U200}
\IndexxS{Theorem!regular variation!smooth variation}
	Let $\ms g$ be regularly varying. Then there exists a function $\ms f$ which is infinitely differentiable and such that 
	$\ms g\sim\ms f$. 
\end{Theorem}

\noindent
Note that in this theorem $\ms f$ is automatically regularly varying with the same index as $\ms g$. 

Let $\ms g,\ms f$ be regularly varying. We say that $\ms f$ is an 
\IndexS{asymptotic inverse}{regular variation!asymptotic inverse} \IndexxS{asymptotic inverse} of $\ms g$ if 
\begin{equation}
\label{U65}
	(\ms g\circ\ms f)(x)\sim(\ms f\circ\ms g)(x)\sim x
	.
\end{equation}
If an asymptotic inverse exists, it is determined uniquely up to $\sim$.
The following result states that an asymptotic inverse exists provided $\ms g$ has positive index.

\begin{Theorem}
\label{U64}
\IndexxS{Theorem!regular variation!asymptotic inversion}
	Let $\ms g$ be regularly varying with index $\alpha>0$, and set 
	\[
		\ms g^-(x) \DE \sup \big\{t\in[1,\infty)\DS \ms g(t)<x\big\}
		.
	\]
	Then $\ms g^-$ is regularly varying with index $\frac{1}{\alpha}$, and is an asymptotic inverse of $\ms g$. 
\end{Theorem}

\noindent
We recall a practical formula for computing asymptotic inverses of functions of a certain form.

\begin{Remark}
\label{U66}
	Let $\ms g$ be a regularly varying function of the form
	\[
		\ms g(r)=r^\rho\big[\ms h(\log r)\big]
	\]
	with $\rho>0$ and $\ms h$ regularly varying. Then 
	\[
		\ms f(r)\DE \rho^{\frac{\Ind\ms h}\rho}\cdot \Big[\frac r{\ms h(\log r)}\Big]^{\frac 1\rho}
	\]
	is an asymptotic inverses of $\ms g$. 
\end{Remark}

\noindent
J.~Karamata also established a result that characterises regular variation of the 
\IndexS{Stieltjes transform}{Stieltjes transform}\IndexxS{Karamata theorem!Stieltjes transform}%
\IndexxS{regular variation!Stieltjes transform} and gives precise information about its size.
We state a formulation taken from \cite{langer.woracek:kara} which includes a boundary case that is often excluded in the
literature.

\begin{Theorem}
\label{U205}
\IndexxS{Theorem!regular variation!Stieltjes transform}
	Let $\mu$ be a measure on $[0,\infty)$, which is not the zero measure and satisfies 
	$\int_{[0,\infty)}(1+t)^{-1}\DD\mu(t)<\infty$, and set 
	\[
		\ms S[\mu](x)\DE \int_{[0,\infty)}\frac{\DD\mu(t)}{t+x},\qquad x>0.
	\]
	Then the following two statements are equivalent:
	\begin{Enumerate}
	\item The distribution function $t\mapsto\mu([0,t))$ is regularly varying
		with index $\alpha$;
	\item $\ms S[\mu]$ is regularly varying with index $\alpha-1$.
	\end{Enumerate}
	If \textup{(i)} and \textup{(ii)} hold, then $\alpha\in[0,1]$ and
	\begin{equation}\label{U206}
		\ms S[\mu](x) \sim C_\alpha\int_x^\infty \frac{\mu([0,t))}{t^2}\DD t,
		\qquad x\to\infty.
	\end{equation}
	with
	\[
		C_\alpha \DE
		\begin{cases}
			\frac{\pi\alpha(1-\alpha)}{\sin(\pi\alpha)} \CAS \alpha\in(0,1),
			\\[1ex]
			1 \CAS \alpha\in\{0,1\}.
		\end{cases}
	\]
	The integral in \cref{U206} is finite for every $x>0$.
\end{Theorem}

\section[{Entire functions of Cartwright class}]{Entire functions of Cartwright class}
\label{U136}

We recall a notion from complex analysis.

\begin{Definition}
\label{U82}
	An entire function $F$ is said to be of \IndexS{Cartwright class}{Cartwright class}, if it is of finite exponential type 
	and satisfies 
	\[
		\int_{\bb R}\frac{\log^+|F(x)|}{1+x^2}\DD x<\infty
		.
	\]
\end{Definition}

\noindent
The following characterisation of Cartwright class goes back to M.G.~Krein \cite{krein:1947} 
(see also \cite[Theorem~6.17]{rosenblum.rovnyak:1994}). 

\begin{Theorem}
\label{U202}
\IndexxS{Theorem!Cartwright class!characterisation}
	An entire function is of Cartwright class if and only if both restrictions $f|_{\bb C^+}$ and $f|_{\bb C^-}$ 
	can be represented as a quotient of two bounded analytic functions in the respective domain $\bb C^+$ or $\bb C^-$. 
\end{Theorem}

\noindent
Functions of Cartwright class have very particular properties when it comes to growth and zero-distribution, see 
\cite{koosis:1998} or also \cite{levin:1980,boas:1954}. We state what is needed in the context of this paper. 

\begin{Theorem}
\label{U196}
\IndexxS{Theorem!Cartwright class!growth and zeroes}
	Let $F$ be an entire function of Cartwright class that has only real zeroes, and satisfies $F=F^\#$ and $F(0)=1$. 
	Let $\tau$ be the exponential type of $F$ and $w_1,w_2,\ldots$ be the zeroes of $F$ listed according to their
	multiplicities. 
	\begin{Enumerate}
	\item The limit $\lim_{R\to\infty}\sum_{|w_n|\leq R}\frac 1{w_n}$ exists, and we have 
		\[
			F(z)=\lim_{R\to\infty}\prod_{|w_n|\leq R}\Big(1-\frac z{w_n}\Big)
			.
		\]
	\item We have $\lim\limits_{r\to\infty}\frac 1r\log\big(\max_{|z|=r}|F(re^{i\vartheta})|\big)=\tau$, and 
		\[
			\forall\vartheta\in(0,\pi)\cup(\pi,2\pi)\DP
			\lim_{r\to\infty}\frac{\log|F(re^{i\vartheta})|}r=\tau\cdot|\sin\vartheta|
			.
		\]
		Along the real axis we have $\limsup\limits_{r\to\infty}\frac{\log|F(\pm r)|}r=0$.
	\item Denote $n_F^+(r)\DE\#\{n\DS w_n\in(0,r)\}$ and $n_F^-\DE\#\{n\DS w_n\in(-r,0)\}$. Then 
		\[
			\lim_{r\to\infty}\frac{n_F^+(r)}r=\lim_{r\to\infty}\frac{n_F^-(r)}r=\frac 1\pi\cdot\tau
			.
		\]
	\item Denote $n_F(r)\DE\#\{n\DS w_n\in(-r,r)\}$. Then 
		\begin{align}
			& \forall r>0\DP \log|F(ir)|=\frac{r^2}2\int_0^\infty\frac 1{t+r^2}\cdot\frac{n_F(\sqrt t)}t\DD t
			\label{U195}
			\\
			& \forall z\in\bb C\DP 
			\log|F(z)|\leq\int_0^{|z|}\frac{n_F(t)}t\DD t+|z|\int_{|z|}^\infty\frac{n_F(t)}{t^2}\DD t
			\nonumber
		\end{align}
		The second relation yields a meaningful result only if $\int_0^\infty\frac{n_F(t)}{t^2}\DD t<\infty$. 
	\end{Enumerate}
\end{Theorem}

\noindent
We note that the relation \cref{U195} implies that 
\[
	\forall r>0\DP n_F(r)\leq \frac 2{\log 2}\cdot\log|F(ir)|
	.
\]
Applying Karamata's theorems to evaluate the asymptotic behaviour of Stieltjes transform and integrals leads to the following 
result (see \cite[Remark~4.16 and Lemma~4.18]{langer.reiffenstein.woracek:kacest-arXiv}). 
It says that the growth of $F(z)$ is governed by its behaviour along the imaginary axis.

\begin{Corollary}
\label{U194}
	Let $F$ be as in \Cref{U196}, and let $\ms g$ be regularly varying with $\lim_{r\to\infty}\ms g(r)=\infty$. 
	\begin{Enumerate}
	\item Assume that $\Ind\ms g\in(0,1)$. Then 
		\begin{align*}
			\big(1-\Ind\ms g\big)\Ind\ms g\cdot &\, \limsup_{r\to\infty}\frac{\log(\max_{|z|=r}|F(z)|)}{\ms g(r)}
			\\
			&\, \leq\limsup_{r\to\infty}\frac{n_F(r)}{\ms g(r)}
			\leq\frac 2{\log 2}\cdot\limsup_{r\to\infty}\frac{\log|F(ir)|}{\ms g(r)}
			.
		\end{align*}
	\item Assume that $\Ind\ms g=0$. Then 
		\begin{equation}
		\label{U193}
		\begin{aligned}
			\limsup_{r\to\infty} &\, \Big(\,\raisebox{4pt}{$\log\big(\max\limits_{|z|=r}|F(z)|\big)$}
			\Big/\raisebox{-4pt}{$\int\limits_1^r\ms g(t)\frac{\DD t}t$}\Big)
			\\
			&\mkern50mu \leq\limsup_{r\to\infty}\frac{n_F(r)}{\ms g(r)}
			\leq\frac 2{\log 2}\limsup_{r\to\infty}\frac{\log|F(ir)|}{\ms g(r)}
			.
		\end{aligned}
		\end{equation}
	\item Assume that $\Ind\ms g=1$ and $\int_1^\infty\frac{\ms g(t)}{t^2}\DD t<\infty$. Then 
		\begin{equation}
		\label{U203}
		\begin{aligned}
			\limsup_{r\to\infty} &\, \Big(\,\raisebox{4pt}{$\log\big(\max\limits_{|z|=r}|F(z)|\big)$}
			\Big/\raisebox{-4pt}{$r\int\limits_r^\infty\frac{\ms g(t)}t\frac{\DD t}t$}\Big)
			\\
			&\mkern50mu \leq\limsup_{r\to\infty}\frac{n_F(r)}{\ms g(r)}
			\leq\frac 2{\log 2}\limsup_{r\to\infty}\frac{\log|F(ir)|}{\ms g(r)}
			.
		\end{aligned}
		\end{equation}
	\end{Enumerate}
\end{Corollary}

\noindent
In \cref{U193} only one inequality gives a meaningful result:
if the left-hand side is positive, the right-hand side is infinite, and if the right-hand side is finite the left-hand side is
zero. This follows from a general fact about entire functions of order zero, cf.\ \cite[Appendix]{berg.pedersen:2007}: 

\begin{Theorem}
\label{U217}
\IndexxS{Theorem!growth of order zero}
	Let $\ms g$ be slowly varying with $\log r=\Smallo(\ms g(r))$, and let $F$ be an entire function with 
	\[
		\limsup_{r\to\infty}\frac{\log(\max_{|z|=r}|F(z)|)}{\ms g(r)}<\infty
		.
	\]
	Then 
	\[
		\limsup_{r\to\infty}\frac{\log(\max_{|z|=r}|F(z)|)}{\ms g(r)}=
		\limsup_{r\to\infty}\frac{\log(\min_{|z|=r}|F(z)|)}{\ms g(r)}
		.
	\]
\end{Theorem}

\noindent
To illustrate the size of the gap between the comparison functions on the left-hand and right-hand sides of 
\cref{U193} and \cref{U203}, we consider Lindel\"of comparison functions. 
In this example all integrals can be computed explicitly. 

\begin{Example}
\label{U204}
	Let $\ms g$ be defined (for sufficiently large $r$) as
	\[
		\ms g(r)=r^\delta \prod_{k=1}^N\big(\log^{[k]}r\big)^{\beta_k}
		,
	\]
	where $\delta\in\{0,1\}$, $N\in\bb N$, and $\beta_1>0$ if $\delta=0$ while $\beta_n<-1$ if $\delta=1$. Then 
	\[
		\begin{cases}
			\int\limits_1^r\ms g(t)\frac{\DD t}t\asymp \ms g(r)\cdot\log r \CAS \delta=0,
			\\
			r\int\limits_r^\infty\frac{\ms g(t)}t\frac{\DD t}t \asymp\ms g(r)\cdot\log r \CAS \delta=1.
		\end{cases}
	\]
\end{Example}

\section[{$J$-inner matrix functions}]{$J$-inner matrix functions}
\label{U130}

In this section we discuss monodromy matrices as stand-alone objects from the function-theoretic viewpoint. 
Recall the notation 
\[
	\IndexN{f^\#}(z)\DE \ov{f(\ov z)},\quad J\DE\begin{pmatrix} 0 & -1 \\ 1 & 0 \end{pmatrix}
	.
\]

\begin{Definition}
\label{U199}
	Let $W\DF\bb C\to\bb C^{2\times 2}$ be a matrix-valued function, and write $W=(w_{ij})_{i,j=1}^2$. 
	We say that $W$ is \IndexS{$J$-inner}{$J$-inner}, if 
	\begin{Enumerate}
	\item the entries $w_{ij}$ are entire and $w_{ij}^\#=w_{ij}$, 
	\item for all $z\in\bb C^+$ we have 
		\[
			\frac 1i\big(W(z)JW(z)^*-J\big)\geq 0
			.
		\]
	\end{Enumerate}
\end{Definition}

\noindent
Assume that $W$ is $J$-inner. For $z\in\bb C^-$ we have $\frac 1i\big(W(z)JW(z)^*-J\big)\leq 0$. Let $z\in\bb R$, then
$W(z)JW(z)^*-J=0$, and since $W(z)\in\bb R^{2\times 2}$ thus $W(z)\in\SL(2,\bb R)$. 
It also follows that $\det W(z)=1$ for all $z\in\bb C$. In particular, the functions $w_{11},w_{12}$ have no common zeroes, and
the same holds for each of the pairs $w_{21},w_{22}$ and $w_{11},w_{21}$ and $w_{12},w_{22}$. 

\begin{Theorem}
\label{U198}
\IndexxS{Theorem!$J$-inner matrix functions}
	Let $W$ be a $J$-inner matrix function.
	\begin{Enumerate}
	\item Each of the quotients
		\begin{equation}
		\label{U192}
			\frac{w_{11}}{w_{21}},\ \frac{w_{12}}{w_{22}},\ \frac{w_{12}}{w_{11}},\ \frac{w_{22}}{w_{21}}
			,
		\end{equation}
		is a Nevanlinna function (cf.\ \Cref{U294}). 
	\item Each entry $w_{ij}$, $i,j\in\{1,2\}$, belongs to the Cartwright class.
	\end{Enumerate}
\end{Theorem}

\noindent
Item (i) follows by elementary manipulations exploiting positivity of $\frac 1i\big(W(z)JW(z)^*-J\big)$. 
Item (ii) is a deeper result due to M.G.~Krein in \cite[Teorema~2]{krein:1951}. 
For us this is an important fact, since it makes available the function theoretic machinery from \Cref{U136}.

We note that a $J$-inner matrix function gives rise to a whole family of Nevanlinna functions, not only those in \cref{U192}. 

In the next proposition we give some properties of meromorphic functions of Nevanlinna class $\mc N_0$.

\begin{Proposition}
\label{U197}
	Let $A,B$ be entire functions with $A=A^\#$, $B=B^\#$, that have no common zeroes, and assume that $\frac BA$ is a
	Nevanlinna function. 
	\begin{Enumerate}
	\item Each of the functions $A,B$ has only real and simple zeroes, and the zeroes of $A$ and $B$ interlace. 
	\item Let $\ms f\DF[0,\infty)\to(0,\infty)$ be a function with $\log r=\Smallo(\ms f(r))$. Then 
		\[
			\forall \vartheta\in(0,\pi)\cup(\pi,2\pi)\DP
			\limsup_{r\to\infty}\frac{\log|A(re^{i\vartheta})|}{\ms f(r)}=
			\limsup_{r\to\infty}\frac{\log|B(re^{i\vartheta})|}{\ms f(r)}
			.
		\]
		If $\ms f$ is regularly varying this also holds for $\vartheta\in\{0,\pi\}$. If 
		$\ms f(r)\sim\sup_{|s-r|\leq 1}\ms f(s)$, then 
		\[
			\limsup_{r\to\infty}\frac{\log(\max_{|z|=r}|A(z)|)}{\ms f(r)}=
			\limsup_{r\to\infty}\frac{\log(\max_{|z|=r}|B(z)|)}{\ms f(r)}
			.
		\]
	\end{Enumerate}
\end{Proposition}

\noindent
Item (i) can be found, e.g., in \cite[Theorem~VII.1]{levin:1980}. Item (ii) follows from the proof
of \cite[Proposition~2.3]{baranov.woracek:smsub}, and by continuity in $\vartheta$ for the case $\vartheta\in\{0,\pi\}$, cf.\ 
\cite[\S16]{levin:1980}. 

Let us make one corollary explicit.

\begin{Corollary}
\label{U201}
	Let $W$ be a $J$-inner matrix function. Then all entries of $W$ have the same order and type. 
\end{Corollary}

%---------
%   FINISH
%---------

\clearpage
\addcontentsline{toc}{section}{\textcolor{Sepia}{List of notation}}
\printindex[not]
\addcontentsline{toc}{section}{\textcolor{Sepia}{Subject index}}
\printindex[sub]

\clearpage
\addcontentsline{toc}{section}{\textcolor{Sepia}{References}}
\printbibliography

{\footnotesize
\begin{flushleft}
	J.~Reiffenstein \\
	Department of Mathematics\\
	Stockholms universitet\\
	106 91 Stockholm\\
	SWEDEN\\
	email: \texttt{jakob.reiffenstein@math.su.se}\\[5mm]
\end{flushleft}
\begin{flushleft}
	H.\,Woracek\\
	Institute for Analysis and Scientific Computing\\
	Vienna University of Technology\\
	Wiedner Hauptstra{\ss}e\ 8--10/101\\
	1040 Wien\\
	AUSTRIA\\
	email: \texttt{harald.woracek@tuwien.ac.at}\\[5mm]
\end{flushleft}
}

\end{document}